\documentclass[reqno,11pt]{amsart}
\usepackage{amssymb, amsmath,latexsym,amsfonts,amsbsy, amsthm,mathtools,cases}
\usepackage{rotating}
\usepackage{enumitem}
\usepackage{tikz}
\usepackage{calc}
\usepackage{tikz-cd}
\usepackage{mathabx}
\usepackage{tikz-3dplot}
\usepackage{hyperref}
\usetikzlibrary{calc}
\usetikzlibrary{decorations.markings,arrows}
\usetikzlibrary{shapes.geometric}
\usepackage{color}
\usepackage{hyperref}
\usepackage{extarrows}
\usepackage{graphicx}
\allowdisplaybreaks[4]
\newcommand\mynumber[1]{%
  \refstepcounter{equation}
  \tag{\theequation}
  \label{#1}
}

\let\pa=\partial
\let\al=\alpha

\let\f=\frac

\let\G=\Gamma

\let\ve=\varepsilon
\let\pa=\partial

\newcommand{\beq}{\begin{equation}}
\newcommand{\eeq}{\end{equation}}
\newcommand{\ben}{\begin{eqnarray}}
\newcommand{\een}{\end{eqnarray}}
\newcommand{\beno}{\begin{eqnarray*}}
\newcommand{\eeno}{\end{eqnarray*}}

\renewcommand{\theequation}{\thesection.\arabic{equation}}

\newtheorem{theorem}{Theorem}[section]

\newtheorem{lemma}[theorem]{Lemma}

\newtheorem{Theorem}{Theorem}[section]

\newtheorem{Lemma}[Theorem]{Lemma}
\newtheorem{Corollary}[Theorem]{Corollary}
\newtheorem{Remark}[Theorem]{Remark}

\newcommand{\ud}{\mathrm{d}}

\newcommand{\nn}{\mathbf{n}}

\newcommand{\II}{\mathbf{I}}

\newcommand{\BR}{\mathbb{R}}

\newcommand{\HQ}{\mathcal{H}_{+}}
\newcommand{\Ker}{\text{Ker}\ \mathcal{H}_{+}}
\newcommand{\Kerb}{(\text{Ker}\ \mathcal{H}_{+})^\bot}

\newcommand{\tr}{\mathrm{Tr}}

\begin{document}
\title[Asymptotic limit for isotropic-nematic dynamics]{On the asymptotic limit for the dynamic isotropic-nematic phase transition with anisotropic elasticity}

\author[H. Dong]{Huan Dong}
\address{School of Mathematical Sciences, Zhejiang University, Hangzhou 310058, China}
\email{huandong\_@zju.edu.cn}

\author[S. Ren]{Siqi Ren}
\address{School of Mathematical Sciences, Zhejiang University of Technology, Hangzhou 310032, China}
\email{sirrenmath@zjut.edu.cn}
\author[W. Wang]{Wei Wang}
\address{School of Mathematical Sciences, Zhejiang University, Hangzhou 310058, China}
\email{wangw07@zju.edu.cn}

\renewcommand{\theequation}{\thesection.\arabic{equation}}
\setcounter{equation}{0}


\begin{abstract}
We consider the isotropic--nematic phase transition dynamics for liquid crystals based on the Landau--de Gennes $Q$-tensor theory.
When the anisotropic elastic constant $L_2$ is negative, we rigorously justify the limit from the Landau--de Gennes flow to a sharp interface system:
the interface which separates the two bulk regions evolves by mean curvature flow; in the isotropic region, the limit alignment tensor \(Q(x,t)\equiv 0\); in the nematic region, the limit $Q(x,t)$ can be characterized by the gradient flow of the Oseen--Frank energy with a strong normal anchoring condition on the interface. This result rigorously verifies a claim by de Gennes \cite{gennes} regarding the surface tension strength of isotropic--nematic interfaces in dynamical settings.
\end{abstract}

\date{\today}
\maketitle

\numberwithin{equation}{section}


\section{Introduction}
Liquid crystal is a state of matter between liquid and solid, in which molecules tend to align along a preferred direction. There are many phases in liquid crystals, such as isotropic, nematic, and smectic et. al. In this paper, we study the dynamics of isotropic--nematic phase transition in the framework of Landau--de Gennes theory \cite{gennes}. Within this model, the state of nematic liquid crystals is described by a symmetric, traceless $3\times3$ matrix $Q$, which is usually called a $Q$-tensor. Physically, $Q$ can be interpreted as the second moment of the orientational distribution function $f$:
$$Q(x)=\int_{\mathbb{S}^2}\Big(m\otimes m-\frac{1}{3}I\Big)f(x,m)dm,$$
where $f(x,m)$ describes the number density of molecules oriented along direction $m$ at material point $x$.
We consider the Landau--de Gennes gradient flow for $x\in \Omega\subset\mathbb{R}^3$ and $t>0$:
\begin{align}\label{equation:main}
 \partial_t Q^\varepsilon(x,t) = \mathcal{L} Q^\varepsilon - \frac{1}{\varepsilon^2} f(Q^\varepsilon),
\end{align}
where
\begin{align}
(\mathcal{L}Q)_{ij}&:=\Delta Q_{ij}+\frac{L_2}{2}\Big(\partial_{ik}Q_{kj}+\partial_{jk}Q_{ki}-\frac{2}{3}\delta_{ij}\partial_{kl}Q_{kl}\Big), \notag\\
f(Q)&:=aQ-b\Big(Q^2-\frac{1}{3}|Q|^2I\Big)+c|Q|^2Q. \label{f}
\end{align}
Here, the parameter $\ve>0$ is small; $a, b, c$ are constants depending on liquid crystal materials; $L_1$, $L_2$ are elastic constants, and in particular, $L_2$ characterizes the anisotropic elasticity of the liquid crystal material. In addition, the Einstein summation convention is used.

The equation \eqref{equation:main} can be regarded as the gradient flow of the well-known Landau--de Gennes energy:
\begin{align}\label{energymain}
\mathcal{F}^\ve(Q,\nabla Q)=&\frac{1}{\ve^2}\int_\Omega \underbrace{\frac{a}{2}\operatorname{Tr}Q^2-\frac{b}{3}\operatorname{Tr}Q^3
+\frac{c}{4}(\operatorname{Tr}Q^2)^2}_{
\triangleq F_b(Q)}dx+\int_\Omega\underbrace{\frac{1}{2}\Big(|\nabla Q|^2+L_2\pa_jQ_{ij}\pa_kQ_{ik}\Big)}_{
\triangleq F_e(\nabla Q)}dx,
\end{align}
where $F_b$ and $F_e$ are called the bulk  and elastic energies, respectively. The bulk energy $F_b$ has two minimal wells at
\begin{align*}
\{0\} \text{ and }\mathcal{N}:=\Big\{s_+(\nn\nn-\frac13\II):\quad \nn\in\mathbb{S}^2,\  s_+=\frac{b+\sqrt{b^2-24ac}}{4c}\Big\},
\end{align*}
which correspond to the isotropic phase and nematic phase in liquid crystals, respectively. We are interested in the phase transition case, i.e., $F_b(0)=F_b(s_+(\nn\nn-\frac13\II))$. In this case, we have
\begin{align*}
  b^2=27ac,\text{ and then }a=\frac{cs_+^2}{3},\quad b=3cs_+.
\end{align*}
In this case $F_b(Q)\ge F_b(0)=0$ for all $Q$. Moreover, if $Q$ has two equal eigenvalues (in this case, it is called {\it uniaxial}), then there exist $s\in\mathbb{R}$ and $n\in\mathbb{S}^2$ such that $Q=s(n\otimes n-\frac{1}{3}I)$, and
 $$F_b(Q)=\frac{c}{9}s^2(s_+-s)^2.$$

The analysis of the isotropic--nematic interface problem was initiated by de Gennes \cite{gennes} in 1971. Under the assumption that $Q$ is diagonal and uniaxial (i.e., Q has two equal eigenvalues), he analyzed a simplified one-dimensional problem and then calculated the surface tension coefficients:
\begin{align}\label{gennesclaim}
    c(L_2)=c_0\begin{cases}
        \sqrt{L_1+\f23L_2},\ L_2<0;\\
        \sqrt{L_1+\f16L_2},\ L_2>0,
    \end{cases}\text{ where }c_0=\sqrt{\f34}\int_0^{s_+}\sqrt{2F_b(s)}ds.
\end{align}
Moreover, near the interface, normal conformation is more favored (i.e., energy is smaller) for
$L_2<0$, while tangential alignment is more favored for $L_2>0$. Since then, several numerical results \cite{Popa1996, Popa1997, KBAM} suggest that de Gennes's calculation is correct for the case $L_2<0$, while the uniaxial assumption near the interface is not valid for $L_2>0$.
In the companion paper by Lin, Zhang and the third author \cite{linwangzhang}, de Gennes' formal derivation \eqref{gennesclaim}  for $L_2<0$ on the isotropic--nematic tension strength, with the (normal) anchoring condition is rigorously confirmed in the static setting. The present paper is devoted to establishing a rigorous mathematical proof for $L_2<0$ in the dynamical setting.

We also remark that the isotropic case ($L_2=0$) is investigated in \cite{FWZZ2018} by Fei-Zhang-Zhang and the third author. It is proved that,  as $\ve\to 0$, the domain $\Omega$ can be divided into two subregions $\Omega_t^-$ and $\Omega_t^+$ such that
\begin{align*}
Q^\ve(x,t)\to 0, \text{ for }x\in\Omega_t^-; \quad Q^\ve(x,t)\to s_+\Big( n(x,t)\otimes n(x,t)-\frac13I\Big)\in \mathcal{N}, \text{ for }x\in\Omega_t^+.
\end{align*}
Moreover, the interface $\Gamma_t=\partial \Omega_t^-\cap \partial \Omega_t^+$ and the unit vector field $n(x,t)$ indeed satisfy the following limit system (assuming $L_1=1$):
\begin{subnumcases}{\label{limit1}}
\partial_t n = \Delta n + |\nabla n|^2 n, \quad  \text{ in } \Omega_t^+,\label{harmonic} \\
\partial_\nu n = 0,  \ \ \ \quad \quad \quad\qquad \text{ on } \Gamma_t,\label{numanncondition} \\
V = \kappa, \ \ \ \ \ \qquad \quad\qquad \text{ on }  \Gamma_t.
\end{subnumcases}
Here $\nu$ and $\kappa$ are the normal vector and mean curvature of the interface $\Gamma_t$ respectively, and $V$ is the normal velocity of $\Gamma_t$.

\subsection{Statement of our main results}
Without loss of generality, we assume
\begin{align*}
a=1,\ b=9,\ c=3,\text{ so that } s_+=1.
\end{align*}
Moreover, we assume $L_1=1$, and denote $L:=L_2<0$ as well as $\beta_L=\sqrt{1+\frac23L}$. We consider $\Omega=\mathbb{T}^3$ to neglect the  technical issues caused by bounded domains.

As $\ve\to 0$, the domain $\Omega$ can also be divided into two subregions $\Omega_t^-$ and $\Omega_t^+$ such that
\begin{align*}
Q^\ve(x,t)\to 0, \text{ for }x\in\Omega_t^-; \quad Q^\ve(x,t)\to n(x,t)\otimes n(x,t)-\frac13I\in \mathcal{N}, \text{ for }x\in\Omega_t^+.
\end{align*}
The limit system for $(\Gamma_t, n(x,t))$ in this case reads as
\begin{subnumcases}{\label{limit}}
(\partial_t n-h)\times n=0, \qquad  \text{ in } \Omega_t^+, \label{limita}\\
n=\nu, \quad\qquad\qquad\qquad\ \text{ on } \Gamma_t, \label{drichleta}\\
V=\beta_L^2 \kappa,\;\;\quad\quad\quad\quad\quad\ \text{ on } \Gamma_t, \label{v=sigmakappa}
\end{subnumcases}
where
$h=-\frac{\delta E(n,\nabla n)}{\delta n}$,
with $E(n,\nabla n)$ the Oseen--Frank energy defined by
\begin{equation}\label{Oseen-Frank energy}
E(n,\nabla n)=\frac12\Big\{|\nabla n|^2+\frac{L}{2}\Big((\nabla\cdot n)^2+|n\times (\nabla \times n)|^2\Big)\Big\}.
\end{equation}
Note that  in this case, the surface tension strength $\beta_L$ coincides with  \eqref{gennesclaim} by de Gennes. In addition, the boundary condition on the interface  $\Gamma$ is just the normal anchoring condition which also agrees with  de Gennes' claim.

 The first result  is the existence of approximate solutions to \eqref{equation:main} valid to an arbitrary order in $\ve$ employing the idea of
``quasi-minimal connecting orbits'' developed by Fei-Lin-Wang-Zhang \cite{FW}.
\begin{theorem}[Construction of approximate solutions]\label{th1}
Let $L\in(-\f32,0)$. Assume that $(\Gamma_t,n)$ is a smooth solution of the limit system \eqref{limit} on $[0,T]$. Then for any $K\in\mathbb{Z}^+$, there exists an approximate solution $Q^{K}$  such that
\begin{equation}\label{equm}
\pa_t Q^{K}=\mathcal{L}Q^K-\frac{1}{\ve^2}f(Q^{K})+\mathcal{R}^K,
\end{equation}
where $\mathcal{R}^K\sim O(\ve^{K-1}),$ and
\begin{align*}
|Q^K(x,t)|\to 0, \text{ for }x\in\Omega_t^-; \quad \Big|Q^K(x,t) -\Big(n(x,t)\otimes n(x,t)-\frac13I\Big)\Big|\to 0, \text{ for }x\in\Omega_t^+.
\end{align*}
\end{theorem}
 We also have the spectral lower bound for the  linear operator around $Q^K$.
\begin{theorem}[Spectral lower bound estimate]\label{th:uA}
Let $L\in(-\f32,0)$ and $Q^K$ be the approximate solution constructed in Theorem \ref{th1}. Then, for any $Q$-tensor $Q\in H^1(\Omega)$ and $t\in [0,T]$, we have
\begin{align}\label{mainineq}
\int_{\Omega}|\nabla Q|^2dx+L\int_{\Omega}|\nabla\cdot Q|^2dx
+\frac{1}{\ve^2}\int_{\Omega}\mathcal{H}_{Q^K}Q:Qdx\geq
-C\int_{\Omega}|Q|^2dx,
\end{align}
where $C$ is a positive constant independent of $\ve$, and $\mathcal{H}_{Q^K}Q$ is the linearized operator defined as
\begin{equation}
    \mathcal{H}_{Q^K}Q=Q+9\Big(\f23I(Q:Q^K)-QQ^K-Q^KQ\Big)
    +3\Big(Q|Q^K|^2+2Q^K(Q:Q^K)\Big).\label{hq0estimta-K}
\end{equation}
\end{theorem}
With the help of the above two theorems, we can obtain the following nonlinear stability of the  approximate solutions constructed above.
\begin{theorem}[Nonlinear stability]\label{th:error}
Let $L\in(-\f32,0)$, $k=9$, $K=k+1$ and $Q^K$ be an approximate solution constructed as in Theorem \ref{th1}. Let $Q^\ve$ be a solution to \eqref{equation:main}. If we have initial condition:
\begin{align*}
\mathcal{E}\big(Q^\ve(\cdot,0)-Q^K(\cdot,0)\big)\leq C\ve^{2k},
\end{align*}
 where
 $$\mathcal{E}(Q)=\sum_{i=0}^{2}\ve^{6i}\int_{\Omega}\|\pa^iQ\|^2dx,$$
 then there exists $\ve_0>0$ such that for any $ \ve \leq \ve_0,$ it holds
\begin{align*}
\mathcal{E}(Q^\ve(\cdot,t)-Q^K(\cdot,t))\leq C\ve^{2k},\ \ \text{for}\ \forall \ t\in [0,T].
\end{align*}
\end{theorem}
We give some remarks on the system \eqref{limit}.
\begin{Remark}
(i) When $L=0$, the energy \eqref{Oseen-Frank energy} reduces to
$\frac12|\nabla n|^2$ and $h=\Delta n$. Then the equation \eqref{limita} is just the well-known harmonic map heat flow:
$$\partial_t n = \Delta n + |\nabla n|^2 n, \ \  x \in \Omega_t^+.$$
In this case, the rigorous limit is justified in \cite{FWZZ2018}.

(ii) In \eqref{drichleta}, the alignment vector $n$ satisfies the Dirichlet boundary condition $n=\nu$, compared with the  Neumann boundary condition $\pa_\nu n=0$ derived in \cite{FWZZ2018} for $L=0$.

(iii) The surface tension $\beta_L$ coincides with \eqref{gennesclaim} for $L_2<0$. Therefore, our results provide rigorous justification of the claim \eqref{gennesclaim}, in the dynamical settings, made by de Gennes \cite{gennes}.
\end{Remark}

\subsection{Background and related results}
There are many results on the asymptotic limits of phase transition problems; an incomplete list of references includes  \cite{Fonseca, deS, C1992, ABC1994, fischer2014, Laux18, Moser90}. In most of these results, the order parameter unknown is a scalar-valued function, whereas in the current isotropic--nematic transition problem for liquid crystals, the order parameter unknown $Q$ is a high-dimensional vector (or a matrix).

The one-dimensional case for the isotropic--nematic interfaces problem has been investigated in many works, see \cite{gennes, PWZZ, Popa1996, Popa1997} for examples. For higher dimensions with  $L_2=0$, the asymptotic limit  was studied in \cite{symliquid} for the radially symmetric case, while the general case is considered by Fei-Wang-Zhang-Zhang in \cite{FWZZ2018} which shows the rigorous limit of the system \eqref{equation:main} to \eqref{limit1}. Moreover,  Laux-Liu \cite{ll2021} provided another proof by using the relative entropy method.

Up to the best of authors' knowledge, there are few analytical results on the high-dimensional isotropic--nematic interfaces problem for non-zero anisotropic elasticity ($L_2\neq0$). Golovaty-Novack-Sternberg-Venkatraman \cite{gsiam,garma} conducted the asymptotic analysis of variational problems based on two dimensional vector models under certain scalings.
Lin-Wang \cite{lw23} analyzed the asymptotic limit for minimizers of the full Ericksen's model, which is also a vector model. In particular,  different boundary conditions (homeotropic, tangential, or free) correspond to different anisotropic elasticities from their analysis. Liu studied the asymptotic limit of an anisotropic Ginzburg-Landau type dynamic model by the relative entropy method \cite{liu24}.

For the $Q$-tensor model, formal derivation of strong anchoring condition $n=\nu$ is conducted in \cite{FWZZ2015}. Park-Wang-Zhang-Zhang \cite{PWZZ} proved that uniaxial solutions with the homeotropic anchoring is stable when $L_2<0$ and becomes unstable for $L_2>0$, which gives a first intuition on the dependence of anchoring alignment on the sign of $L_2$.  The present work might be the first result on the rigorous analysis  of isotropic--nematic phase transition dynamics with anisotropic elasticity.

We also remark that, the isotropic--nematic phase transition problem with isotropic elasticity ($L_2=0$) can be viewed as a special case of the Keller-Rubinstein-Sternberg problem proposed in \cite{RSKf,RSKr}. For this problem, Lin-Pan-Wang \cite{LPW2012} studied the asymptotic behavior  as $\ve\to 0$ for minimizers of the energy functionals. Later, Lin-Wang \cite{lw2019} proved the partial regularity for the minimizers of the limit problem.
For the dynamical problem, Bronsard-Stoth \cite{Bron} studied a special $\mathbb{R}^2$-valued model where the bulk potential vanishes on two concentric circles. They also assumed the solution is radially symmetric.
In \cite{FW}, Fei-Lin-Wang-Zhang  proved the asymptotic limit as $\ve \to 0$ for a typical matrix-valued case. In particular, two different kinds of boundary conditions on the moving interface are derived.  Lately, Liu \cite{liucpam} studied the dynamical Keller-Rubinstein-Sternberg problem by the
relative entropy method.

There are several different methods to study asymptotic limits of dynamical phase transition problems. The first method is the matched asymptotic expansion,
which is first proposed by de Mottoni-Schatzmann in \cite{deS} to study the classical Allen-Cahn dynamics. This method has been further developed by Alikakos-Bates-Chen \cite{ABC1994} to study the asymptotic behavior of the Cahn-Hilliard equation. Thereafter, this method was applied to other phase transition problems, e.g., \cite{Abelsliu2018, FWZZ2015, FWZZ2018, FW} and references therein. Another approach is the geometric measure theory method based on varifold flows, which is introduced in \cite{il} to study the classical Allen-Cahn model and further applied by \cite{Chen96, HT, PP}. This method relies on a crucial monotonicity formula utilizing the ``discrepancy function" which seems difficult to apply to vectorial phase transition problems. The third method is the relative entropy (or modified energy) method, developed by Fischer-Laux-Simon \cite{FLS2020} recently (motivated by \cite{jerrard, fischerj, Lin98}). This method does not rely on the comparison principle, and has also been applied to different phase transition problems, e.g., \cite{ll2021, AFM, liu24, liucpam}.

In this paper, we adopt the classical matched asymptotic expansion method, since both of the latter methods encounter certain essential obstacles when applied to high-dimensional two-phase interface problems.

\subsection{Main difficulties and key ideas.}

The main framework of the proof in this paper follows \cite{ABC1994, FW}: first, for a given smooth solution $(\Gamma_t, n(x,t))$ to the limit system \eqref{limit} on time interval $[0, T]$, we construct approximate solutions to \eqref{equation:main} using the asymptotic expansion method; second, we prove a lower bound estimate for the spectrum of the linearized operator near the approximate solutions; and finally, we prove the stability of the approximate solutions. In the first two steps, we encounter new difficulties.

In the first step,  as in \cite{ABC1994, FW}, we adopt different expansions for two distinct regions: the region away from the interface and the region near the interface, which are called the outer expansion and the inner expansion, respectively. A key point in this step is the solvability of a linearized matrix-valued ODE arising in the inner expansion; see the equation \eqref{Qm}. By an appropriate decomposition as in \cite{FWZZ2018, FW}, we reduce this system to three scalar-valued ODEs, given by equations \eqref{sm0-gamma}-\eqref{sm3-gamma}. The solvability of the first and third equations has already been discussed in \cite{FWZZ2018}, but the second equation is new due to the anisotropic elasticity $L\neq 0$. We establish its solvability by constructing fundamental solutions of the associated homogeneous ODE with prescribed asymptotic behavior at $\pm\infty$, proving the non-vanishing of the relevant Wronskians, and then deriving a Green's--function representation via variation of parameters. This representation yields the existence and uniqueness of a bounded solution, together with the required exponential-decay estimates.

Another key point in the first step is to prove the linearity of the expansion system at each order, which ensures that the system is solvable on the whole interval [0,T]; otherwise, the solution could only be expected to exist on a smaller time interval. This is not a trivial matter and requires exploiting many structural properties of the problem; see Lemma 3.9 and its proof.

There are also new difficulties in the second step, namely the lower bound estimate for the spectrum of the linearized operator. Since this operator is matrix-valued, we need to decompose it as in \cite{FWZZ2018, FW} and reduce the domain of the estimates to a one-dimensional interval along the normal direction near the interface. For the case $L=0$, the proof of these estimates relies strongly on the natural boundary conditions $\partial_\nu n=0$, but in the present situation this condition is no longer valid. We make use of a key div-curl decomposition from \cite{linwangzhang}:
 \begin{align} \label{key-div-curl}
|\nabla \cdot Q|^2= \frac{2}{3} |\nabla Q|^2 - \frac{1}{6} |\mathcal{T}(Q)|^2 - \frac{2}{3}\Big(\partial_k Q_{li} \partial_l Q_{ki} - \partial_k Q_{ki} \partial_l Q_{li}\Big),
\end{align}
where the operator $\mathcal{T}(Q)$ is defined as
\begin{equation*}
\mathcal{T}(Q) = \big( T_{ij}(Q) \big)_{1 \leq i, j \leq 3}, \quad T_{ij}(Q) = \varepsilon^{ikl} \partial_k Q_{lj} + \varepsilon^{jkl} \partial_k Q_{li}.
\end{equation*}
Since the term $\int_{\Omega}(\partial_k Q_{li} \partial_l Q_{ki} - \partial_k Q_{ki} \partial_l Q_{li})dx$ vanishes  due to the Stokes' theorem, one uses \eqref{key-div-curl} to convert \eqref{mainineq} into the following:
\begin{align}\label{mainineqL}
\int_{\Omega}&\Big(1+\frac{2L}{3}\Big)|\nabla Q|^2dx+\frac{1}{\ve^2}\int_{\Omega}\mathcal{H}_{Q^K}Q:Qdx-\frac{L}{6}\int_{\Omega}|\mathcal{T}(Q)|^2dx\geq -C\int_{\Omega}|Q|^2dx.
\end{align}
Moreover, through a delicate analysis, we find that, with the help of good estimates for $|\mathcal{T}(Q)|^2$, the normal anchoring condition exactly allows us to close the entire estimate in this case.

The rest of the paper is organized as follows: In Sections 2 and 3, we derive the systems of each order in outer and inner expansion respectively. Then we present the procedure for solving the system of each order by induction in Section 4. The approximate solutions are constructed by gluing the two expansions in Section 5, which completes the proof of Theorem \ref{th1}. In Section 6, we prove the uniform spectral lower bound estimate for the linearized operator, which yields Theorem \ref{th:uA}. In Section 7, we show the nonlinear stability of the approximate solutions by the standard energy method and then conclude Theorem \ref{th:error}.

\subsection{Notations}
\begin{itemize}
\item Let $d(x,\Gamma_t)$ be the signed distance from $x$ to $\Gamma_t$, and $\nu=\nabla d|_{\Gamma_t}$ be the unit outer normal of $\Omega_t^-$. We denote
\begin{align*}
 &   \Gamma(\delta)=\{(x,t)\in\Omega\times[0,T]: |d(x,\Gamma_t)|<\delta \},\\
 & U_\pm =\{(x,t)\in\Omega\times[0,T]: d(x,\Gamma_t)\gtrless 0\}.
\end{align*}
We also write $\Gamma=\{(x,t)\in\Omega\times[0,T]:x\in\Gamma_t\}$ for simplicity.
\item For any two vectors $m=(m_1,m_2,m_3)$ and $n=(n_1,n_2,n_3)$, we denote
$$mn=m\otimes n=[m_in_j]_{1\le i,j\le 3}.$$
\item The matrices in this paper are $3\times 3$ tensors  that are symmetric and traceless, unless otherwise stated.
\item $\{n,l,m\}$ represents a smooth orthogonal vectors basis  coordinate frame for $\mathbb{R}^3$.
  \item For any two matrices $Q$ and $R$, we denote $$Q:R=\tr(QR)=Q_{ij}R_{ji},\ \  \tr Q^3=Q_{ij}Q_{jk}Q_{ki}, \ \ |Q|^2=Q:Q.$$
    Considering the derivatives, we denote $$Q_{ij,j}=\pa_{x_j}Q_{ij},\ \ \nabla\cdot Q=\pa_{x_i}Q_{ij}.$$
\item
Define the multi-linear operators $B(\cdot,\cdot)$ and $C(\cdot,\cdot,\cdot)$ as follows
\begin{align}
\begin{split}
B(Q,P)&:=9\Big(\frac{2}{3}I(Q:P)-QP-PQ\Big),\label{taylorB}\\
C(Q,P,R)&:=3\Big(Q(P:R)+R(Q:P)+P(Q:R)\Big).
\end{split}
\end{align}
For a given matrix $P$, the linear operator $\mathcal{H}_{P}(\cdot)$ takes as follows
\begin{align}
\mathcal{H}_{P}Q:=&Q+9\Big(\frac{2}{3}I(Q:P)-QP-PQ\Big)
+3\Big(Q|P|^2+2P(Q:P)\Big)\notag\\
=&Q+B(Q,P)+C(Q,P,P)\label{taylor1}.
\end{align}
   \item Let $\alpha>0$.
   By $f(z,x,t) = O(e^{-\alpha|z|})$ as $z \to \pm\infty$ we mean that there exist constants $C > 0$ and $k \ge 0$ such that
       \begin{align}
       |f(z,x,t)|\leq C|z|^ke^{-\al|z|},\;\; \;\text{as}\;z\to\pm\infty.
\end{align}
\item The Levi-Civita symbol $\ve^{ijk}$ is defined by
 \begin{equation}\label{levisivita}
\ve^{ijk}=\left\{
\begin{array}{ll}
1,&\text{if }(i,j,k)\text{ is }(1,2,3),\ (2,3,1)\ \text{or }(3,1,2);\\
-1,&\text{if }(i,j,k)\text{ is }(3,2,1),\ (1,3,2)\ \text{or }(2,1,3);\\
0,&\text{if } i=j,\ j=k\text{ or }i=k.
\end{array}
\right.
\end{equation}
In particular, for any vectors $a,b,c\in\mathbb{R}^3$, we have
\begin{equation*}
\ve^{ijk}a_jb_k=(a\times b)_i,\quad \text{and}\quad \ve^{ijk}a_ib_jc_k=a\cdot(b\times c).
\end{equation*}
Moreover, we have
\begin{equation*}
\ve^{ikl}\ve^{imn}=\delta_{km}\delta_{ln}-\delta_{kn}\delta_{lm}, \quad \text{and }\quad\ve^{ikl}\ve^{jmn}=
\left|\begin{array}{cccc}
\delta_{ji}&\delta_{mi}&\delta_{ni}\\
\delta_{jk}&\delta_{mk}&\delta_{nk}\\
\delta_{jl}&\delta_{ml}&\delta_{nl}
\end{array}\right|.
\end{equation*}
\item  Denote
\begin{align}\label{def:s(z)}
s(z)=\frac{\exp(\gamma z)}{1+\exp(\gamma z)},\ \gamma=\frac{1}{\sqrt{1+\frac{2L}{3}}}>0.
\end{align}
Then $s(z)$ is the solution of
 \begin{equation*}
\begin{cases}
\begin{array}{ll}
(1+\frac{2L}{3})s''+f(s)=0,\\
 s(-\infty)=0,\ s(+\infty)=1,
\end{array}
\end{cases}
\end{equation*}
with $f(s)=-2s^3+3s^2-s$.
\item We denote
\begin{align*}
u(z)\sim \bar{u}(z)\text{ as }z\to+\infty\;\;\text{if}\;\;\lim_{z\to+\infty} \f{u(z)}{\bar{u}(z)}=1\;\;\text{and}\;\;\lim_{z\to+\infty} \f{u'(z)}{\bar{u}'(z)}=1.
\end{align*}
The notation $u(z)\sim \bar{u}(z)\text{ as }z\to-\infty$ is defined similarly.
\end{itemize}
\section{Outer expansion}\label{section2}
In this section, we perform an outer expansion in $U_\pm$ rather than $U_\pm \backslash \Gamma(\delta/2)$.
For $(x,t)\in U_\pm$, we assume that the solution of \eqref{equation:main} has the form
\begin{equation}\label{outerex}
Q^\ve(x,t)=Q_0^\pm(x,t)+\ve Q_1^\pm(x,t)+\ve^2Q_2^\pm(x,t)+\cdots,\quad \
\end{equation}
Then we expand \eqref{f} as follows
\begin{equation}\label{fex}
f(Q^\ve)=f(Q_0^\pm)+\ve\mathcal{H}_{\pm}Q_1^\pm+\sum_{k\geq 2}\ve^k\Big(\mathcal{H}_{\pm}Q_k^\pm+B_{k-1}^\pm+C_{k-1}^\pm\Big),
\end{equation}
where
\begin{align}
&\mathcal{H}_{\pm}Q=\mathcal{H}_{Q_0^\pm}Q,\label{H} \\
&B_{k-1}^\pm=\frac{1}{2}\sum_{i=1}^{k-1}B(Q_{i}^\pm,Q_{k-i}^\pm),\label{B}\\
&C_{k-1}^\pm=\sum_{i=1}^{k-1}C(Q_i^\pm,Q_{k-i}^\pm,Q_0^\pm)+3\sum_{\substack{i+j\leq k-1\\1\leq i,j\leq k-1}}Q_{k-i-j}^\pm(Q_i^\pm:Q^\pm_j).\notag
\end{align}
By substituting \eqref{outerex} and \eqref{fex} into the equation \eqref{equation:main} and collecting  terms of the same order, we have
\begin{align}
O(\ve^{-2})\ \text{term}:&\ f(Q_0^\pm)=0,\label{oleading}\\
O(\ve^{-1})\ \text{term}:&\ \mathcal{H}_{\pm}Q_1^\pm=0,\label{1order}\\
O(\ve^{k-2})\ \text{term}\ (k\geq 2):&\ \mathcal{H}_{\pm}Q_{k}^\pm=-\partial_tQ_{k-2}^\pm+\mathcal{L}Q_{k-2}^\pm-B_{k-1}^\pm-C_{k-1}^\pm,\label{lorder}
\end{align}
where
$$(\mathcal{L}Q)_{ij}:=\Delta Q_{ij}+\frac{L_2}{2}\Big(\partial_{ik}Q_{kj}+\partial_{jk}Q_{ki}-\frac{2}{3}\delta_{ij}\partial_{kl}Q_{kl}\Big).$$
\subsection{Solving $Q_0^+$ and $Q_k^-(k\geq 0)$}
For the $O(\ve^{-2})$ order \eqref{oleading}, we find that $Q_0^\pm$ are two critical points of the bulk energy $F$.
We can take
\begin{equation}\label{Q0}
Q_0^-=0,\ \  Q_0^+=n(x,t)n(x,t)-\frac{1}{3}I,\ \ n\in \mathbb{S}^2.
\end{equation}
 Therefore, the linear operator \eqref{H} takes in the explicit form:
 \begin{equation}\label{Hpm}
 \mathcal{H}_+Q=9\Big(Q-nn Q-Q nn+\frac{2}{3}(nn: Q)I\Big)+6\Big(nn-\frac{1}{3}I\Big)(nn:Q).
 \end{equation}

 Our purpose is to construct $Q_{k}^{\pm}$ ($k\geq 1$) in $U_{\pm}$.  Indeed, since $Q^-=0$ is a solution of \eqref{equation:main} and $Q_0^-=0$, we take
\begin{align}\label{def:Q_m^-}
Q_k^-=0\;\;\text{in}\; U_-,\;\; k\geq 0.
 \end{align}

It remains to construct  $Q_k^+$ ($k\geq 1$) in $U_+$.
 \subsection{Decomposition of $Q_k^+\ (k\geq 1)$}\label{decom2.2}
To begin with, we give a lemma for the linear operator $\mathcal{H}_+$ mainly based on \cite{WZZ}, where the orthogonality \eqref{orth:E} for the basis can be verified directly.
\begin{Lemma}\label{ker}
For $\mathcal{H}_+$ defined in \eqref{Hpm}, it holds that
\begin{itemize}
\item[(1)]
$\text{Ker}\ \mathcal{H}_+=\{nn^\bot +n^\bot n|n^\bot\in \mathbb{V}_n\}$, with $\mathbb{V}_n\triangleq \{n^\bot |n^\bot \cdot n=0\}$.
\item[(2)]
$\mathcal{H}_+$ is a one-to-one mapping on the orthogonal complement $(\text{Ker}\ \mathcal{H}_+)^\bot$ where
$$(\text{Ker}\ \mathcal{H}_+)^\bot\triangleq \{Q|(Q:nn^\bot+n^\bot n)=0,\ \forall n^\bot \in \mathbb{V}_n\}.$$
The inverse operator takes
\begin{equation*}
  \mathcal{H}_{+}^{-1}Q=\frac{1}{9}\Big(Q-nn\cdot Q-Q\cdot nn+\frac{2}{3}(nn:Q)I\Big)+\frac{14}{9}\Big(nn-\frac{1}{3}I\Big)(nn:Q).
\end{equation*}
\end{itemize}
\end{Lemma}

Assume that $l(x,t),\ m(x,t)$ are smooth unit vector fields in $U_+$ such that $l(x,t)\cdot m(x,t)=l(x,t)\cdot n(x,t)=m(x,t)\cdot n(x,t)=0$.
We introduce the orthogonal basis $\{E^i,\;0\leq i\leq 4\}$ of the space of symmetric and traceless $3\times 3$ matrices as follows
\begin{align}
\begin{split}
&E^0(x,t)=\Big(nn-\frac{1}{3}I\Big),\ \ E^1(x,t)=\Big(nl+ln\Big),\ \ E^2(x,t)=\Big(nm+mn\Big),\label{def:E-al}\\
&E^3(x,t)=\Big(ml+lm\Big),\ \ E^4(x,t)=\Big(ll-mm\Big),
\end{split}
\end{align}
with
 \begin{equation}\label{orth:E}
 E^\alpha:E^\beta=\left\{
\begin{array}{ll}
\frac{2}{3},\ \alpha=\beta=0;\\
2,\ \alpha=\beta\neq 0;\\
0,\ \alpha\neq \beta.
\end{array}
\right.
\end{equation}
It follows by Lemma \ref{ker} that span $\{E^0,E^3,E^4\}\subseteq \Kerb$ and span $\{E^1,E^2\}\subseteq \Ker$. Thus, we introduce the decomposition:
\begin{align}
\begin{split}
\label{decom}
Q_k^+(x,t)&=Q^\bot_k(x,t)+Q^\top_k(x,t),\\
Q^\bot_k(x,t)&=q_{k,0}(x,t)E^0(x,t)+q_{k,3}(x,t)E^3(x,t)+q_{k,4}(x,t)E^4(x,t)\in \Kerb,\\
Q^\top_k(x,t)&=q_{k,1}(x,t)E^1(x,t)+q_{k,2}(x,t)E^2(x,t)\in\Ker,
\end{split}
\end{align}
  where the coefficients $q_{k,i}(x,t)$ for $i=0,\cdots, 4$ will be determined later.
  \subsection{Solving the  evolution of alignment  vector field $n$ and $Q_1^+(Q_1^\top)$, $Q_2^\bot$}
 Fix $Q_0^+$ and $Q_1^+$, and recall the definition of operators $B$, $C$ in \eqref{taylorB}. Then we define two linear operators:
\begin{equation}\label{def:calB-calC}
\mathcal{B}_{+}Q=B(Q_1^+,Q),\ \ \mathcal{C}_+Q=C(Q_1^+,Q_0^+,Q).
\end{equation}
\begin{Lemma}\cite[Lemma 2.3]{FWZZ2018}\label{Kerb}
For any  $Q\in \Ker$, we have $\mathcal{B}_+Q,\ \mathcal{C}_+Q\in \Kerb$.
\end{Lemma}
Moreover, we have
\begin{align}
&\mathcal{C}_+E^0=2Q_1^+,\quad\mathcal{C}_+E^1=6q_{1,1}E^0,\quad\mathcal{C}_+E^2=6q_{1,2}E^0,\quad\mathcal{C}_+E^3=\mathcal{C}_+E^4=0,\notag\\
  &\mathcal{B}_+E^0=\mathcal{B}_+E^3=\mathcal{B}_+E^4=0,\quad  \mathcal{B}_+E^1=-9q_{1,1}E^0-9q_{1,2}E^3-9q_{1,1}E^4,\label{bccalcu}\\
  &\mathcal{B}_+E^2=-9q_{1,2}E^0-9q_{1,1}E^3+9q_{1,2}E^4.\notag
  \end{align}
Then we arrive at that for any $Q\in (\Ker)^\top$, we have $\mathcal{B}_+Q=0.$
We can write
$$\mathcal{B}_+Q_k^+=\mathcal{B}_+Q_k^\top.$$

For $k=1$, it follows by $\mathcal{H}_{+}Q_1^+=0$ that
\begin{equation}\label{Q1}
Q_1^+=Q_1^\top=q_{1,1}E^1+q_{1,2}E^2,
\end{equation}
where $q_{1,1}$ and $q_{1,2}$  are determined later.

For $k=2,$ we deduce from \eqref{lorder}, $B_1^+=\frac{1}{2}\mathcal{B}_+Q_1^\top$, $C_1^+=\mathcal{C}_+Q_1^\top$, Lemmas \ref{ker} and \ref{Kerb} that
\begin{equation}\label{1term}
-\partial_tQ_0^++\mathcal{L}Q_0^+=\HQ Q_2^\bot+\f12\mathcal{B}_+Q_1^\top+\mathcal{C}_+Q_1^\top\in \Kerb.
\end{equation}
It is a compatibility condition on the solvability of \eqref{lorder} for $k=2$.
We can derive the evolution  equation for the alignment  vector field $n$.
\begin{Lemma} The alignment vector field $n$ obeys \begin{equation}\label{nev}
(\pa_t n-h)\times n=0,
\end{equation}
where $h=-\frac{\delta E(n,\nabla n)}{\delta n}$ with the Oseen--Frank energy $E(n,\nabla n)$  defined by
\begin{equation*}
E(n,\nabla n)=\frac12\Big\{|\nabla n|^2+\frac{L}{2}\Big((\nabla\cdot n)^2+|n\times (\nabla \times n)|^2\Big)\Big\}.
\end{equation*}
\end{Lemma}
\begin{proof}
 Indeed, thanks to
Lemma \ref{ker}, we have
\begin{align*}\label{contract:Q0}
(\partial_tQ_0^+-\mathcal{L} Q_0^+):(nn^\bot+n^\bot n)=0,
\end{align*}
which together with $Q_0^+=(nn-\frac{1}{3}I)$ gives
\begin{equation*}
\pa_t n\cdot n^\bot-h\cdot n^\bot=0.
\end{equation*}
Hence, we conclude that $n$ satisfies \eqref{nev}.
\end{proof}
Now we are in a position to solve $Q_1^+$ and $Q_2^{\bot}$. For $Q_2^\bot$, we claim
\begin{equation}\label{Q2solve}
Q_2^\bot=\mathcal{H}_{+}^{-1}A_1,
\end{equation}
where $A_1$ only depends on $Q_0^+$ and $Q_1^+$.
Indeed, \eqref{Q2solve} follows by a reformulation of  \eqref{1term}:
\begin{equation*}
A_1:=-\partial_tQ_0^++\mathcal{L}Q_0^+-\f12\mathcal{B}_+Q_1^\top-\mathcal{C}_+Q_1^\top
=\HQ Q_2^\bot\in \Kerb.
\end{equation*}
For $Q_1^+$, thanks to \eqref{Q1}, it suffices to solve $q_{1,1}$ and $q_{1,2}$.

Once  $Q_{1}^+$ is solved, $Q_{2}^\bot$ is determined by \eqref{Q2solve}. In fact, {the outer expansion only gives parabolic equations for $q_{1,1}$ and $q_{1,2}$ on $U_+$, while the boundary condition on $\Gamma$ is obtained from  the inner expansion in  \eqref{bd:q111} and \eqref{bd:q121}. }

The equations are derived as follows:

 Thanks to the definitions \eqref{B} and \eqref{def:calB-calC}, it follows  $B_2^+=\mathcal{B}_+Q_2^+=\mathcal{B}_+Q_2^\top$ and $C_2^+=2\mathcal{C}_+Q_2^++3|Q_1^+|^2Q_1^+$. Consequently,  we write the equation \eqref{lorder} ($k=3$) as
$$\mathcal{H}_{+}Q_{3}^+=-\partial_tQ_{1}^\top+\mathcal{L}Q_{1}^\top-B_{2}^+-C_{2}^+.$$
Together with \eqref{Q1} and Lemma \ref{Kerb}, we arrive at  
\begin{equation*}\label{1terms}
-\partial_tQ_1^\top+\mathcal{L}Q_1^\top-2\mathcal{C}_+Q_2^\bot-3|Q_1^\top|^2Q_1^\top=\HQ Q_{3}^\bot+\mathcal{B}_+Q_2^\top+2\mathcal{C}_+Q_2^\top\in\Kerb.
\end{equation*}
Using span $\{E^1,E^2\}\subseteq \text{Ker}\mathcal{H}_+$, then for $i=1,2$, we have
\begin{equation}\label{comQ2}
\big(-\partial_tQ_1^\top+\mathcal{L}Q_1^\top-2\mathcal{C}_+Q_2^\bot-3|Q_1^\top|^2Q_1^\top\big):E^i=0.
\end{equation}
To  simplify \eqref{comQ2}, we need the following lemma.
\begin{Lemma}
For any $Q_1^\top=q_{1,1}E^1+q_{1,2}E^2$, we have
\begin{align*}
2\mathcal{C}_+&\mathcal{H}_+^{-1}(1/2\mathcal{B}_+Q_1^\top+\mathcal{C}_+Q_1^\top)-3|Q_1^\top|^2Q_1^\top=0.
\end{align*}
\end{Lemma}
\begin{proof}
For convenience, we give some basic calculations:
\begin{align*}
  &\mathcal{H}_+^{-1}E^0=E^0,\quad \mathcal{H}_+^{-1}E^3=\f19E^3,\quad\mathcal{H}_+^{-1}E^4=\f19E^4.
\end{align*}
Then we have
\begin{align*}
    \mathcal{C}_+Q_1^\top
    &=6(q_{1,1}^2+q_{1,2}^2)E^0,\\
    \mathcal{B}_+Q_1^\top
    &=-9(q_{1,1}^2+q_{1,2}^2)E^0-18q_{1,1}q_{1,2}E^3+9(q_{1,2}^2-q_{1,1}^2)E^4,\\
    1/2\mathcal{B}_+Q_1^\top+\mathcal{C}_+Q_1^\top&=\f32(q_{1,1}^2+q_{1,2}^2)E^0-9q_{1,1}q_{1,2}E^3+\f92(q_{1,2}^2-q_{1,1}^2)E^4,\\
    \mathcal{H}_+^{-1}(1/2\mathcal{B}_+Q_1^\top+\mathcal{C}_+Q_1^\top)&=\f32(q_{1,1}^2+q_{1,2}^2)E^0-q_{1,1}q_{1,2}E^3+\f12(q_{1,2}^2-q_{1,1}^2)E^4.
\end{align*}
Using the above calculations, we deduce that
\begin{align*}
2\mathcal{C}_+\mathcal{H}_+^{-1}(1/2\mathcal{B}_+Q_1^\top+\mathcal{C}_+Q_1^\top)-3|Q_1^\top|^2Q_1^\top
=2\mathcal{C}_+(3/2 (q_{1,1}^2+q_{1,2}^2)E^0)-3|Q_1^\top|^2Q_1^\top=0.
\end{align*}
\end{proof}
Using the lemma, we have
\begin{align*}
    -2\mathcal{C}_+&Q_2^\bot- 3|Q_1^\top|^2Q_1^\top=-2\mathcal{C}_+\mathcal{H}_+^{-1}A_1- 3|Q_1^\top|^2Q_1^\top\\
    &=-2\mathcal{C}_+\mathcal{H}_+^{-1}(-\pa_t Q_0^++\mathcal{L}Q_0^+)+2\mathcal{C}_+\mathcal{H}_+^{-1}(1/2\mathcal{B}_+Q_1^\top+\mathcal{C}_+Q_1^\top)- 3|Q_1^\top|^2Q_1^\top\\
    &=-2\mathcal{C}_+\mathcal{H}_+^{-1}(-\pa_t Q_0^++\mathcal{L}Q_0^+).
\end{align*}
By the definition of $\mathcal{C}_+$, we conclude that $-2\mathcal{C}_+Q_2^\bot- 3|Q_1^\top|^2Q_1^\top$ is linear on $q_{1,i}$.

Substituting \eqref{Q1} and \eqref{Q2solve} into \eqref{comQ2}, we obtain the desired  linear parabolic equations in $U_+$:
\begin{align}
\pa_tq_{1,1}&-\Delta q_{1,1}-\frac{L}{2}(\pa_{ik}q_{1,1}n_kn_i+\pa_{ik}q_{1,1}l_kl_i+\pa_{ik}q_{1,2}m_kl_i)=r_{1,1},\label{q11}
\end{align}
and
\begin{align}
\pa_tq_{1,2}&-\Delta q_{1,2}-\frac{L}{2}(\pa_{ik}q_{1,1}l_km_i+\pa_{ik}q_{1,2}n_kn_i+\pa_{ik}q_{1,2}m_km_i)=r_{1,2},\label{q12}
\end{align}
where $r_{1,i}$ is a linear function for $q_{1,1},\ q_{1,2},\ \nabla q_{1,1}$ and $\nabla q_{1,2}$.

Note that we are only concerned with  $Q_2^\bot$ and $Q_1$ in this step.
The treatment of $Q_2^\top$ will be addressed in the subsequent step.
\subsection{Solving $Q_{k+1}^\top$ and $Q_{k+2}^\bot$ $(k\geq 1)$}
The construction follows by induction.
Assuming that $Q_i^\top$ and $Q_{i+1}^\bot$ $(1\leq i\leq k)$  have been  determined, we will proceed to construct  $Q_{k+1}^\top$ and $Q_{k+2}^\bot$.

To obtain $Q_{k+2}^\bot$, we write the $O(\ve^{k})$ order equation as
$$\mathcal{H}_{+}Q_{k+2}^+=-\partial_tQ_{k}^++\mathcal{L}Q_{k}^+-B_{k+1}^+-C_{k+1}^+.$$
Together with the definitions \eqref{B},  \eqref{def:calB-calC} and the decomposition  \eqref{decom},  one has
\begin{align}
\mathcal{H}_{+}Q_{k+2}^\bot&=-(\pa_t-\mathcal{L})Q_{k}^+-\mathcal{B}_+Q_{k+1}^+-\tilde{B}_{k}-\mathcal{C}_+Q_{k+1}^+-\tilde{C}_{k}\nonumber\\
&=-\mathcal{B}_+Q_{k+1}^\top-2\mathcal{C}_+Q_{k+1}^\top+D_{k+1}.\label{Qm+1bot}
 \end{align}
Indeed, for $k\geq 1$, one has
 \begin{align}
 \notag B_{k+1}^+&=B(Q_{k+1}^+,Q_{1}^+)+\frac{1}{2}\sum_{i=2}^{k}B(Q_i^+,Q_{k+2-i}^+):=\mathcal{B}_+Q_{k+1}^\top+\tilde{B}_{k},\\
\label{Cm}
C_{k+1}^+&=2C(Q_1^+,Q_{k+1}^+,Q_0^+)+\sum_{i=2}^{k}C(Q_i^\pm,Q_{k+2-i}^\pm,Q_0^\pm)+3\sum_{\substack{i+j\leq k+1\\1\leq i,j\leq k+1}}Q_{k+2-i-j}^\pm(Q_i^\pm:Q^\pm_j)\\
&:=2\mathcal{C}_+Q_{k+1}^++\tilde{C}_{k},\notag\\
\notag D_{k+1}&=-\tilde{B}_{k}-\tilde{C}_{k}-(\pa_t-\mathcal{L})Q_{k}^+-2\mathcal{C}_+Q_{k+1}^\bot.
\end{align}
Note that $\tilde{B}_{k}$ and $\tilde{C}_{k}$ only depend on $Q_0^+,\dots,Q_{k}^+$, while $D_{k+1}$ depends on $Q_0^+,\dots,Q_{k}^+$ and $Q_{k+1}^\bot.$

Consequently, applying \eqref{Qm+1bot} and Lemma \ref{ker}, we have
 \begin{equation}\label{Qm+1=H-1}
Q_{k+2}^\bot=\mathcal{H}_{+}^{-1}\Big(-\mathcal{B}_+Q_{k+1}^\top-2\mathcal{C}_+Q_{k+1}^\top+D_{k+1}\Big)
:=\mathcal{H}_{+}^{-1}A_{k+1},
 \end{equation}
  where $A_{k+1}$ depends on $Q_0^+,\dots,Q_{k}^+$ and $Q_{k+1}^\bot$, which follows by the assumption.

 To solve $$Q_{k+1}^\top=q_{k+1,1}E^1+q_{k+1,2}E^2,$$ for the $O(\ve^{k+1})$ order, we deduce from \eqref{lorder}, \eqref{Cm} and \eqref{Qm+1=H-1}  that
\begin{align*}
-\pa_tQ_{k+1}^+&+\mathcal{L}Q_{k+1}^+\\
=&\mathcal{H}_{+}Q_{k+3}^++B_{k+2}^++C_{k+2}^+\notag\\
=&\mathcal{H}_{+}Q_{k+3}^\bot+(\mathcal{B}_++2\mathcal{C}_+)Q_{k+2}^\bot+(\mathcal{B}_+
+2\mathcal{C}_+)Q_{k+2}^\top
+\tilde{B}_{k+1}+\tilde{C}_{k+1}\nonumber\\
=&\mathcal{H}_{+}Q_{k+3}^\bot-(\mathcal{B}_++2\mathcal{C}_+)\mathcal{H}_{+}^{-1}\Big(\mathcal{B}_+Q_{k+1}^\top+2\mathcal{C}_+Q_{k+1}^\top-D_{k+1}\Big)\\
&+(\mathcal{B}_+
+2\mathcal{C}_+)Q_{k+2}^\top
+\tilde{B}_{k+1}+\tilde{C}_{k+1}. \nonumber
\end{align*}
Together with Lemma \ref{Kerb}, the above equation gives  following compatibility condition for $Q_{k+1}^\top$:
 \begin{align}
 -\pa_tQ_{k+1}^\top&+\mathcal{L}Q_{k+1}^\top+(\mathcal{B}_++2\mathcal{C}_+)\mathcal{H}_{+}^{-1}(\mathcal{B}_++2\mathcal{C}_+)Q^{\top}_{k+1}\nonumber\\
 &-(\pa_t Q_{k+1}^\bot-\mathcal{L}Q_{k+1}^\bot)-\tilde{B}_{k+1}-\tilde{C}_{k+1}
 -(\mathcal{B}_++2\mathcal{C}_+)\HQ^{-1}D_{k+1}\nonumber\\
 =&\HQ Q_{k+3}^\bot+(\mathcal{B}_++2\mathcal{C}_+)Q_{k+2}^\top\in \Kerb.\label{qm+1}
 \end{align}
Analogously to \eqref{comQ2}, we can deduce the linear equations for $q_{k+1,1}$ and $q_{k+1,2}$ in $ U_+$:
\begin{align}
\pa_tq_{k+1,1}&-\Delta q_{k+1,1}-\frac{L}{2}(\pa_{ik}q_{k+1,1}n_kn_i+\pa_{ik}q_{k+1,1}l_kl_i+\pa_{ik}q_{k+1,2}m_kl_i)=r_{k+1,1},\label{qk+11}
\end{align}
and
\begin{align}
\pa_tq_{k+1,2}&-\Delta q_{k+1,2}-\frac{L}{2}(\pa_{ik}q_{k+1,1}l_km_i+\pa_{ik}q_{k+1,2}n_kn_i+\pa_{ik}q_{k+1,2}m_km_i)=r_{k+1,2},\label{qk+12}
\end{align}
where $r_{k+1,i}$ is a linear function for $q_{k+1,1},\ q_{k+1,2},\ \nabla q_{k+1,1}$ and $\nabla q_{k+1,2}$.
The boundary conditions on $\Gamma$  will be given in the inner expansion. Once  $Q_{k+1}^\top$ is solved, $Q_{k+2}^\bot$ can be deduced by \eqref{Qm+1=H-1}.
\section{Inner expansion}\label{section3}
 In this section, we perform an inner expansion in $\Gamma(\delta)$. In this region, the true solution exhibits rapid variations. To capture these changes, we introduce a new variable, $z=\frac{d^\ve(x,t)}{\ve}$. We will regard $z$ as an independent variable and $(x, t)\in\Gamma(\delta)$ as parameters. Then we will solve a series of ODE systems with respect to $z$. Formally, we write the inner expansion as
 \begin{align}\label{innerexpansion}
\tilde{Q}^\ve(x,t)=Q^\ve(z,x,t)=\sum_{k=0}^{+\infty}\ve^k Q_k(\frac{d^\ve}{\ve},x,t),
\end{align}
with
\begin{align}\label{dve}
d^\ve(x,t)=\sum_{k=0}^{+\infty}\ve^kd_k(x,t),
\end{align}
 which is a signed distance function with respect to a surface $\G_t^\ve$. Then we have $|\nabla d^\ve|^2=1$, which gives
\begin{equation}\label{didj}
\nabla d_0\cdot\nabla d_k=\left\{
\begin{array}{ll}
1,&k=0,\\
0,&k=1,\\
-\frac{1}{2}\sum_{j=1}^{k-1}\nabla d_j\cdot\nabla d_{k-j},&k\geq 2.
\end{array}
\right.
\end{equation}

 Since we would like to approximate the sharp interface system \eqref{limit}, we take $d_0(x,t)$ as the signed distance function to $\G_t$. We deduce that $\nabla d_0\cdot\nabla=\pa_\nu$ for $(x,t)\in\G$, where $\pa_\nu$ is the normal derivative. In addition, if $\ve$ is sufficiently small, $\G_t^\ve$ should be a good approximation of $\G_t$.
 We shall keep in mind that $Q_k$ has to satisfy the matching conditions:
 \begin{equation}\label{mathingcondition}
     \lim_{z\rightarrow  \pm\infty}|\pa_t^m\pa_x^l\pa_z^n(Q_k(z,x,t)-Q_k^\pm(x,t))|= O(e^{-\gamma |z|})\ \text{for}\ (x,t)\in\Gamma(\delta),\ m,l,n\geq 0.
 \end{equation}
Otherwise we cannot glue the solution in the overlapped region $\Gamma(\delta)\backslash \Gamma(\frac{\delta}{2})$.

 Substituting \eqref{innerexpansion} and \eqref{dve} into  \eqref{equation:main}, we have
\begin{equation}\label{innereq}
D_tQ^\ve=\Delta_x Q^\ve+\frac{L}{2}\Big(D_{ik}Q^\ve_{kj}+D_{jk}Q^\ve_{ki}-\frac{2}{3}\delta_{ij}D_{kl}Q_{kl}^\ve\Big)-\frac{1}{\ve^2}f(Q^\ve),
\end{equation}
where
\begin{align}
&D_t Q^\ve:=\partial_t Q^\ve+\frac{\partial_td^\ve}{\ve}\partial_z Q^\ve,\nonumber\\
&\Delta_x Q^\ve:=\Delta Q^\ve+\frac{1}{\ve}(2\nabla d^\ve\cdot \nabla \partial_zQ^\ve+\Delta d^\ve\partial_z Q^\ve)+\frac{1}{\ve^2}\partial_z^2Q^\ve,\notag\\
&D_{ik}Q^\ve_{kj}:=\pa_{ik}Q^\ve_{kj}+\frac{1}{\ve}(\pa_{iz}Q^\ve_{kj}\pa_{k}d^\ve +\pa_{kz}Q^\ve_{kj}\pa_id^\ve+\pa_{z}Q^\ve_{kj}\pa_{ik}d^\ve)+\frac{1}{\ve^2}\pa_z^2Q^\ve_{kj} \pa_{k}d^\ve\pa_id^\ve,\nonumber\\
&D_{jk}Q^\ve_{ki}:=\pa_{jk}Q^\ve_{ki}+\frac{1}{\ve}(\pa_{jz}Q^\ve_{ki}\pa_kd^\ve+\pa_{kz}Q^\ve_{ki}\pa_jd^\ve+\pa_zQ^\ve_{ki}\pa_{jk}d^\ve)+\frac{1}{\ve^2}\pa_z^2Q^\ve_{ki}\pa_kd^\ve\pa_jd^\ve,\nonumber\\
&D_{kl}Q^\ve_{kl}:=\pa_{kl}Q^\ve_{kl}+\frac{1}{\ve}(2\pa_{kz}Q^\ve_{kl}\pa_ld^\ve+\pa_zQ^\ve_{kl}\pa_{kl}d^\ve)+\frac{1}{\ve^2}\pa_z^2Q^\ve_{kl}\pa_ld^\ve\pa_kd^\ve\nonumber,\\
&f(Q^\ve):=f(Q_0)+\ve\mathcal{H}_{Q_0}Q_1+\sum_{k\geq 2}\ve^k\Big(\mathcal{H}_{Q_0}Q_k+B_{k-1}+C_{k-1}\Big),\nonumber
\end{align}
with
\begin{align}
&\mathcal{H}_{Q_0}Q=Q+9\Big(\f23I(Q:Q_0)-QQ_0-Q_0Q\Big)
+3\Big(Q|Q_0|^2+2Q_0(Q:Q_0)\Big),\label{op:calH}\\
&B_{k-1}=9\Big(\f13I\sum_{1\leq i\leq k-1}\big(Q_i:Q_{k-i}\big)-\sum_{1\leq i\leq k-1}\big( Q_iQ_{k-i}\big)\Big),\label{bm-1}\\
&C_{k-1}=3\Big(2\sum_{1\leq i\leq k-1}Q_i(Q_0:Q_{k-i})+\sum_{\substack{i+j\leq k\\1\leq i,j\leq k-1}}Q_{k-i-j}(Q_i:Q_j)\Big).\label{cm-1}
\end{align}
As in \cite{ABC1994}, we modify the system \eqref{innereq} as follows:
\begin{equation}\label{modieq}
D_tQ^\ve=\Delta_x Q^\ve+\frac{L}{2}\Big(D_{ik}Q^\ve_{kj}+D_{jk}Q^\ve_{ki}-\frac{2}{3}\delta_{ij}D_{kl}Q_{kl}^\ve\Big)-\frac{1}{\ve^2}f(Q^\ve)-G^\ve(z,x,t)(d^\ve-\ve z),
\end{equation}
where $$G^\ve(z,x,t)=\sum_{k=0}^{+\infty}\ve^{k-2}G_k(z,x,t),\qquad \text{ for }(z,x,t)\in\mathbb{R}\times \Gamma(\delta).$$
Here $G_k(z,x,t)$ are $3\times 3$ trace-free matrices to be determined later.
The modified term does not change the system \eqref{innereq} on $\{d^\ve=\ve z\}$. For convenience, we define
\begin{equation}\label{op:L_nu}
 \mathcal{L}_{\nu}Q=\partial_z^2Q+\frac{L}{2}\Big(\partial_z^2Q(\nabla d_0\nabla d_0)+(\nabla d_0\nabla d_0)\partial_z^2 Q -\frac{2I}{3}(\nabla d_0\nabla d_0):\partial_z^2Q\Big).
 \end{equation}
Then we have the following systems of order $\ve^{k-2}$ for $k\geq 0$ respectively.
\begin{itemize}
\item The $O(\ve^{-2})$ system takes the form
\begin{align}\label{inner-2}
\mathcal{L}_{\nu}Q_0-f(Q_0)=G_0d_0.
\end{align}
\item For $k\geq 1$, the $O(\ve^{k-2})$ system takes the form
\begin{align}
\mathcal{L}_{\nu}&Q_{k}-\mathcal{H}_{Q_0}Q_{k}
=F_{k}+G_{k}d_0.\label{Qm}
\end{align}
\begin{itemize}[label=$\blacktriangle$]
\item For $k=1,$ we have
\begin{align}
F_1=&\partial_z Q_0(\partial_t d_0-\Delta d_0)-2\nabla d_0\cdot\nabla \partial_zQ_0-\frac{L}{2}\mathcal{N}_1(\pa_zQ_0,\nabla d_0)\nonumber\\
&-\frac{L}{2}\sum_{\substack{q+r=1\\0\leq q,r\leq 1}}\mathcal{N}_2(\pa_z^2Q_0,\nabla d_q,\nabla d_r)
+G_0(d_1-z),\label{f1}
\end{align}
with
\begin{align*}
\mathcal{N}_1(Q, \nabla d)_{ij}
&=\partial_{ik}d Q_{kj}+\partial_{jk}d Q_{ki}
-\frac{2}{3}\delta_{ij}\partial_{kl}d Q_{kl}\\
&\quad+\partial_i d\partial_{k} Q_{kj}+\partial_jd\partial_{k}Q_{ki}
+\partial_kd\partial_{i}Q_{kj}
+\partial_kd\partial_{j}Q_{ki}
-\frac{4}{3}\delta_{ij}\partial_kd\partial_{l}Q_{kl},
\end{align*}
and
\begin{align}
\mathcal{N}_2&(Q, \nabla d_A, \nabla d_B)_{ij}\notag\\
=& Q_{kj}\pa_kd_A\pa_id_B+Q_{ki}\pa_kd_A\pa_j d_B-\frac{2}{3}\delta_{ij}Q_{kl}\pa_kd_A\pa_ld_B. \label{N2}
\end{align}
\item For $k\geq 2$, one has
\begin{align}
F_k=&\sum_{\substack{p+q=k-1\\0\leq p,q\leq k-1\\}}(\pa_t-\Delta )d_p\pa_zQ_q-2\nabla d_p\cdot \nabla\pa_zQ_q+(\pa_t -\mathcal{L})Q_{k-2}\nonumber\\
&\quad-\frac{L}{2}\sum_{\substack{p+q=k-1\\0\leq p,q\leq k-1}}\mathcal{N}_1(\pa_zQ_p,\nabla d_q)
-\frac{L}{2}\sum_{\substack{p+q+r=k\\0\leq p\leq k-1}}\mathcal{N}_2(\pa_z^2Q_p,\nabla d_q,\nabla d_r)\nonumber\\
&\quad-B_{k-1}-C_{k-1}+\sum_{\substack{p+q=k\\1\leq q\leq k}}G_pd_q-G_{k-1}z,\label{Fm}
\end{align}
with
\begin{align*}
(\mathcal{L}Q)_{ij}&:=\Delta Q_{ij}+\frac{L}{2}\Big(\partial_{ik}Q_{kj}+\partial_{jk}Q_{ki}-\frac{2}{3}\delta_{ij}\partial_{kl}Q_{kl}\Big).
\end{align*}
\end{itemize}
\end{itemize}
We seek $Q_k$ for $k\geq 0$  which takes the form
\begin{equation}\label{Qmesim}
Q_k(z,x,t)=Q_k^\top(z,x,t) +Q_k^\bot(z,x,t),
\end{equation}
where
\begin{align*}
Q_k^\bot(z,x,t)&=s_{k,0}(z,x,t)E^0(x,t)+s_{k,3}(z,x,t)E^3(x,t)+s_{k,4}(z,x,t)E^4(x,t),\\
Q_k^\top(z,x,t)&=s_{k,1}(z,x,t)E^1(x,t)+s_{k,2}(z,x,t)E^2(x,t).
\end{align*}
Direct calculation by the definition gives the following.
\begin{lemma}
Let $Q,P$ be symmetric and traceless matrix.
It holds that
\begin{align}
\mathcal{N}_1(Q,\nabla d):P&=2(\pa_{ik}dQ_{kj}+\pa_id\pa_{k}Q_{kj}+\pa_kd\pa_{i}Q_{kj})P_{ij},\label{N1:def}\\
\mathcal{N}_2(Q,\nabla d_A,\nabla d_B):P&=2Q_{kj}\pa_k d_A \pa_id_BP_{ij}.\label{N2:def}
\end{align}
Moreover, one has
\begin{align}\label{N1:PQ}\nonumber
\mathcal{N}_1(Q,\nabla d):P+ \mathcal{N}_1(P,\nabla d):Q
&=2\Big(\pa_i\big(\pa_kd (Q_{ij}P_{kj}+Q_{kj}P_{ij})\big)\Big)\\
&=2\nabla\cdot \Big(\big(PQ+QP\big)\cdot \nabla d\Big),\\ \label{N2:PQ}
 \mathcal{N}_2(Q,\nabla d_A,\nabla d_B):P&=\mathcal{N}_2(P,\nabla d_B,\nabla d_A):Q.
\end{align}
In particular,
\begin{align}\label{N1:QQ}
\mathcal{N}_1(Q,\nabla d):Q
&=2\nabla\cdot \Big( Q^2\cdot\nabla d\Big).
\end{align}
\end{lemma}

\subsection{Solving $Q_0$ }
Thanks to the matching condition \eqref{mathingcondition}, we solve the $O(\ve^{-2})$ order equation on $(z,x,t)\in \mathbb{R}\times \G(\delta)$ as
\begin{equation}
\begin{cases}
\begin{array}{ll}
\mathcal{L}_{\nu}Q_0-f(Q_0)=G_0d_0,\\
Q_0(-\infty,x,t)=0,\ Q_0(+\infty,x,t)=Q_0^+(x,t).\label{in2}
\end{array}
\end{cases}
\end{equation}
We introduce an auxiliary operator:
 \begin{align}\label{op:Lstar}
 \mathcal{L}_{*}Q= &\partial_z^2Q+\frac{L}{2}\Big(\partial_z^2Q(nn)+(nn)\partial_z^2 Q -\frac{2I}{3}(nn):\partial_z^2Q\Big).
 \end{align}
Then the equation
\begin{align*}
 \mathcal{L}_{*}Q_0-f(Q_0)=0
\end{align*}
 has a solution
 \begin{equation}\label{q0}
    Q_0(x,t,z)=s(z)\Big(n(x,t)n(x,t)-\frac{1}{3}I\Big),
 \end{equation}
where
\begin{align}\label{def:s}
s(z)=\frac{\exp(\gamma z)}{1+\exp(\gamma z)} \qquad\text{with } \gamma= \Big(1+\frac{2L}{3}\Big)^{-1/2},
\end{align}
is the solution of
 \begin{equation}\label{sz}
\begin{cases}
(1+\frac{2L}{3})s''-(s-3s^2+2s^3)=0,\\
 s(-\infty)=0,\ s(+\infty)=1.
\end{cases}
\end{equation}
Thus, from the strong anchoring condition
 \begin{align}\label{anchoringcondition}
 n(x,t)=\nabla d_0 \qquad \text{for }(x,t)\in\Gamma,
 \end{align}
 we know that, for $(x,t)\in\Gamma$, \eqref{q0} is a solution of \eqref{in2}.

For $(x,t) \in \Gamma(\delta)$, the vectors $n(x,t)$ and $\nabla d_0$ are not necessarily equal.
In this case, we define a smooth vector field $\bar{h}$ in $\Gamma(\delta)$ such that
\begin{equation}\label{d0de}
n=\nabla d_0+d_0\bar{h}(x,t),\quad \text{for }  (x,t)\in\Gamma(\delta).
\end{equation}
In particular $\bar{h}(x,t)=\nabla d_0\cdot\nabla (n-\nabla d_0)=\partial_\nu n$ which implies
\begin{equation}\label{h0cdotn=0}
 \bar{h}\cdot n=0, \quad \text{for }(x,t)\in\Gamma.
 \end{equation}
Then one has
 \begin{align*}
 \mathcal{L}_{\nu}Q-\mathcal{L}_{*}Q=&~\frac{Ld_0}{2}\Big(\partial_z^2Q P_*+ P_*\partial_z^2 Q -\frac{2I}{3} P_*:\partial_z^2Q\Big),\\
 \text{with }P_*=&~d_0^{-1}\big((\nabla d_0\nabla d_0)-nn\big)=-(n\bar{h}+\bar{h}n)+d_0\bar{h}\bar{h}.
 \end{align*}

Consequently,  if we choose
 \begin{align}\nonumber
 G_0(z,x,t)=&\frac{1}{d_0}\Big(\mathcal{L}_{\nu}Q_0-\mathcal{L}_{*}Q_0\Big)\\
 =&\frac{Ls''}{2}\Big\{-\frac{1}{3}(n\bar{h}+\bar{h}n)-2nn(n\cdot \bar{h})+\frac{8}{9}(\bar{h}\cdot n) I\nonumber\\ \nonumber
 &+d_0\Big((n\cdot \bar{h})(n\bar{h}+\bar{h}n) -\frac{2}{3}(\bar{h}\cdot n)^2I-\frac{2}{3}\bar{h}\bar{h}+\frac{2}{9}|\bar{h}|^2I\Big)\Big\}\\
 \triangleq&\frac{Ls''}{2}g_0(x,t),  \label{G0}
 \end{align}
then $Q_0(z)$ defined in \eqref{q0} is a solution to \eqref{in2} in $\G(\delta)$.
In particular, for $(x,t)\in\G$, we have
\begin{align}
&g_0(x,t)=-\frac{1}{3}(n\bar{h}+\bar{h}n).\label{g0e}
\end{align}
\subsection{Decomposition of $Q_k$ ($k\geq 1$)}
The goal of this subsection is to decompose \eqref{Qm} into several
scalar ODEs. Similar to Section \ref{decom2.2}, by employing \eqref{def:E-al}, we write  $Q_k\ (k\geq 1)$ in terms of the orthogonal basis
\begin{align}\label{Qme}
Q_k(z,x,t)=\sum_{i=0}^{4}s_{k,i}(z,x,t)E^i(x,t),
\end{align}
where
\begin{align*}
\begin{split}
&E^0(x,t)=\Big(nn-\frac{1}{3}I\Big),\ \ E^1(x,t)=\Big(nl+ln\Big),\ \ E^2(x,t)=\Big(nm+mn\Big),
\\
&E^3(x,t)=\Big(ml+lm\Big),\ \ E^4(x,t)=\Big(ll-mm\Big),
\end{split}
\end{align*}
with $n$ determined in \eqref{limit}.

To begin with, utilizing \eqref{q0},  a tedious calculation of \eqref{op:calH}  gives
\begin{align}
\mathcal{H}_{Q_0}Q_k=&(1+6s+2s^2)Q_k-9s(nn Q_k+Q_k nn)+2(3s^2nn+s(3-s)I)(nn:Q_k)\nonumber\\
=&\left(1-6s+6s^2\right)s_{k,0}(z,x,t)E^0
+\left(1-3s+2s^2\right)\sum_{i=1,2}
s_{k,i}(z,x,t)E^{i}\notag\\
&+\left(1+6s
+2s^2\right)\sum_{{i}=3,4}s_{k,i}(z,x,t)E^{i}.\label{hq0m}
\end{align}
Calculation of \eqref{op:L_nu} gives
 \begin{align}
 \begin{split}\label{op:L_nuQk}
 &\mathcal{L}_{\nu}Q_k
 =\sum_{i=0}^4\pa_z^2s_{k,i}(z,x,t)P^{i}(x,t),\\
 &\text{where}\quad P^{i}:=E^{i}+\f{L}{2}\Big(E^{i}(\nabla d_0\nabla d_0)+\nabla d_0(\nabla d_0 E^{i})-\f23\left((\nabla d_0\nabla d_0):E^{i}\right)I\Big).
  \end{split}
\end{align}
Recall the auxiliary operator $\mathcal{L}_{*}Q$, we have
 \begin{align}\label{op:LstarQk} \nonumber
 \mathcal{L}_{*}Q_k= &\partial_z^2Q_k+\frac{L}{2}\Big(\partial_z^2Q_k (nn)+(nn)\partial_z^2 Q_k -\frac{2I}{3}(nn):\partial_z^2Q_k\Big)\\
=& \Big(1+\f{2L}{3}\Big)\pa_z^2s_{k,0}E^0
 +\Big(1+\f{L}{2}\Big)\sum_{i=1,2}\pa_z^2s_{k,i}E^{i}
+\sum_{i=3,4}\pa_z^2s_{k,i}E^{i}.
 \end{align}

Using the above expressions, we define
 \begin{align}\label{def:hm}
\tilde{g}_k(z,x,t):=
\begin{cases}
\begin{array}{ll}
\f{\left(\mathcal{L}_{\nu}Q_k-\mathcal{L}_{*}Q_k\right)
(z,x,t)}{d_0(x,t)},&(x,t)\in\Gamma(\delta)/\Gamma,\\
\pa_\nu{(\mathcal{L}_{\nu}Q_k-\mathcal{L}_{*}Q_k
)(z,x,t)},&(x,t)\in\Gamma,
\end{array}
\end{cases}
\end{align}
which implies
\begin{align}\label{fm:hm}
\tilde{g}_kd_0+\mathcal{L}_{*}Q_k-\mathcal{L}_{\nu}Q_k=0,\;\;\text{for}\;(x,t)\in\Gamma(\delta).
\end{align}
Since $s(z)$ is monotonic,  we define a smooth  cut-off function $\eta(z)$ as follows
 \begin{align*}\eta(z)=0\ \text{if} \ z\leq -1,\ \eta(z)=1\ \text{if}\ z\geq 1.\end{align*}
 In particular, we can take $\eta(-z)=1-\eta(z)$.
Then, we construct $G_k\ (k\geq 1)$ as follows
\begin{align}\label{Gm}
G_k(z,x,t):=\sum_{i=0}^{2}\eta_i'(z)g_{k,i}(x,t)E^i(x,t)+\tilde{g}_{k}(z,x,t),\;\;\text{for}\; (x,t)\in \Gamma(\delta),
\end{align}
where $g_{k,i}\ (i=0, 1,2)$ will be determined later, and
\begin{align*}
  \eta_0(z)=\eta(z),\quad \eta_1(z)=\eta_2(z)= -\frac{L}{6}s'(z).
\end{align*}
Thanks to \eqref{op:L_nuQk} and \eqref{def:hm}, there exists a constant $C>0$ only depending on $d_0,n,l,m$ such that
\begin{align}\label{est:hm}
\big|\tilde{g}_k|_{\G}\big|\leq C\sum_{i=0}^4\left|\pa_z^2s_{k,i}|_{\G}\right|,\;k\geq 1.
\end{align}
To fulfill the matching condition \eqref{mathingcondition}, we require
 \begin{align*}
  \lim_{z\to-\infty} Q_k(z,x,t)=Q_{k}^{-}(x,t)\;\;\text{and}\;\;\lim_{z\to+\infty} Q_k(z,x,t)=Q_{k}^{+}(x,t) .
 \end{align*}
Thus, it  suffices to take
  \begin{align}\label{boundary-sm}
  \lim_{z\to-\infty} s_{k,i}(z,x,t)=0\;\;\text{and}\;\; \lim_{z\to+\infty} s_{k,i}(z,x,t)=q_{k,i}(x,t),\quad i=0,\cdots,4.
 \end{align}

The system \eqref{Qm} is reduced to
\begin{align}\label{Qm-simple}
\mathcal{L}_{*}Q_{k}-\mathcal{H}_{Q_0}Q_{k}=F_{k}+d_0\sum_{i=0}^2g_{k,i}(x,t)\eta_i'(z)E^i(x,t).
\end{align}

Combining \eqref{hq0m}, \eqref{op:LstarQk} and \eqref{boundary-sm}, we obtain the following equivalent scalar equations of $s_{k,i}(0\le i\le 4)$ for $z\in\mathbb{R}$:
\begin{align}
\label{sm0-gamma}&
\begin{cases}
\begin{array}{ll}
-(1+\frac{2L}{3})\pa_z^2s_{k,0}+\left(1-6s+6s^2\right)s_{k,0}
=-\frac{3}{2}F_k:E^0-\eta'(z)d_0g_{k,0},\\
s_{k,0}(-\infty,x,t)=0,\ s_{k,0}(+\infty,x,t)=q_{k,0}(x,t);
\end{array}
\end{cases}\\
\label{sm1-gamma}
\quad &\begin{cases}
\begin{array}{ll}
-(1+\frac{L}{2}) \pa_z^2s_{k,i}
+\left(1-3s+2s^2\right)s_{k,i}
=-\frac{1}{2}F_k:E^{i}+\frac{L}{6}s''(z)d_0g_{k,i},\\
s_{k,i}(-\infty,x,t)=0,\ s_{k,i}(+\infty,x,t)=q_{k,i}(x,t),
\end{array}
\end{cases}&&\text{for $i=1,2$;}\\
\label{sm3-gamma}
\quad &\begin{cases}
\begin{array}{ll} -\pa_z^2s_{k,i}+\left(1+6s
+2s^2\right)s_{k,i}=-\frac{1}{2}F_k:E^{i},\\
s_{k,i}(-\infty,x,t)=0,\ s_{k,i}(+\infty,x,t)=q_{k,i}(x,t).
\end{array}
\end{cases}&&\text{for $i=3,4$.}
\end{align}

\subsection{Compatibility conditions} We consider the solvability of \eqref{sm0-gamma}-\eqref{sm3-gamma}.
Let
\begin{equation}\label{the-kap-l}
\theta(s)=1-6s+6s^2,\ \kappa(s)=1-3s+2s^2\ \text{and}\ \iota (s)=1+6s+2s^2.
\end{equation}
To begin with, we list some basic properties for $s(z),\ \theta(s),\ \kappa(s)$ and $\iota(s)$ as follows.
\begin{lemma}\label{lem:pro-s} For $s(z),\ \theta(s),\ \kappa(s)$ and $\iota(s)$, we have the following properties.
\begin{itemize}
\item[(1)]
 $s'(z)$ satisfies the homogeneous equation of \eqref{sm0-gamma}, i.e.,
\begin{align*}
-(1+\frac{2L}{3})s'''(z)+\theta(s(z))s'(z)=0.
\end{align*}
Moreover, it holds that
\begin{equation}\label{s'}
    s'=\gamma s(1-s),
\end{equation}
and for all $k\ge 0$ that
 \begin{align}\label{est:s(z)}
|s(z)|+\left|\partial_z^k s(z)\right|\leq Ce^{-\gamma |z|},\; z<0,\;\;\text{and}\; \;|s(z)-1|+\left|\partial_z^{k+1} s(z)\right|\leq Ce^{-\gamma |z|},\; z>0.
\end{align}
\item[(2)]
We have
\begin{align}\label{kappaf}
\kappa(s(z))=-\frac{f(s(z))}{s(z)}=\Big(1+\frac{2L}{3}\Big)\frac{s''(z)}{s(z)},\ \theta(s(z))=\Big(1+\frac{2L}{3}\Big)\frac{s'''(z)}{s'(z)},
\end{align}
and $\iota(s(z))>1$.
\item [(3)]
It holds that
\begin{align}
&s(-\infty)=0,\;\; \ s(+\infty)= 1,\label{spminfty}\\
&\theta(s(-\infty))=1,\ \ \theta(s(+\infty))=1,\label{thetas}\\
&\kappa(s(-\infty))=1,\ \ \kappa(s(+\infty))= 0,\label{kappas}\\
&\iota (s(-\infty))=1,\ \ \iota (s(+\infty))=9.\label{iotas}
\end{align}
\end{itemize}
\end{lemma}

Define the spaces:
\begin{align*}
&\mathcal{S}(\gamma)=\big\{f\in C^{\infty}(\mathbb{R}): f^\pm=\lim_{z\to\pm\infty}f(z) \text{ exist,}\\
&\qquad\text{and for } \forall j\ge 0,~\exists k\ge 0, \text{s.t.}~~
|\partial_z^j(f(z)-f^{\pm})|\lesssim |z|^{k}e^{-\gamma |z|}, \text{as }z\to \pm\infty\big\}.\\
&\mathcal{S}_{L,M}(\gamma)=\big\{f(\cdot,x,t)\in \mathcal{S}(\gamma) : \text{for }\forall ~(j,l',m')\in\mathbb{N}\!\times\![0,L]\!\times\![0,M],~\exists k\ge 0, \text{s.t.}~~\\
&\qquad\big|\partial_z^j\partial_x^{l'}\partial_t^{m'}\big(f(z,x,t)-f^{\pm}(x,t)\big)\big|\lesssim |z|^{k}e^{-\gamma |z|}, \text{as }z\to \pm\infty\big\}.
\end{align*}

The following lemma concerns the solvability of \eqref{sm0-gamma}. We refer to \cite[Lemma 2.6]{FWZZ2018} for the proof.
\begin{lemma}\label{lem1}
Assume that $f_0(z,x,t)\in \mathcal{S}_{L,M}(\gamma)$ for all $(x,t)\in \Gamma(\delta)$ with $f_0(+\infty,x,t)=f_{0}^{+}(x,t)$ and $f_0(-\infty,x,t)=0$.
The equation
 \begin{equation}\label{ode:s0}
\begin{cases}
-(1+\frac{2L}{3})\pa_z^2s_{0}(z,x,t)+\theta(s)s_0(z,x,t)=f_0(z,x,t),\\
s_0(0,x,t)=0,
\end{cases}
\end{equation}
 has a unique bounded solution $s_0(z,x,t)$ if and only if
\begin{equation}\label{solve1}
\int_{\mathbb{R}}f_0(z,x,t)s'(z)dz=0.
\end{equation}
In this case, the solution $s_0(z,x,t)$ is given by
\begin{align}
\begin{split}
s_0(z,x,t)&=\Big(1+\f{2L}{3}\Big)^{-1}
\int_{0}^{z}\f{s'(z)}{(s'(y))^2}\int_{y}^{+\infty}
s'(w)f_0(w,x,t)dw dy,
\end{split}
\end{align}
and satisfies
$s_0(z,x,t)\in \mathcal{S}_{L,M}(\gamma)$ for all $(x,t)\in \Gamma(\delta)$ with $s_0(+\infty,x,t)=f_{0}^{+}(x,t)$ and $s_0(-\infty,x,t)=0$.
\end{lemma}

For the solvability of \eqref{sm1-gamma}, we let $u_\pm(z)$ be smooth solutions of the homogeneous equations
\[
\begin{cases}
-\bigl(1+\frac{L}{2}\bigr) u_+'' + \kappa(s) u_+ = 0,\\
u_+(z) \sim 1 \quad \text{as } z \to +\infty,
\end{cases}
\quad\text{and}\quad
\begin{cases}
-\bigl(1+\frac{L}{2}\bigr) u_-'' + \kappa(s) u_- = 0,\\
u_-(z) \sim \exp\bigl( \frac{z}{\sqrt{1+L/2}} \bigr) \quad \text{as } z \to -\infty.
\end{cases}
\]
Denote by $W(u_-,u_+) = u_-'(z)u_+(z) - u_+'(z)u_-(z)$ the Wronskian, and let
$$W_A=\Big(1+\frac{L}{2}\Big)W(u_-,u_+).$$

\begin{lemma}\label{lem2}
Let $f_1(z,x,t)\in\mathcal{S}_{L,M}(\gamma)$  with $f_1(\pm\infty, x,t)=0$.
The boundary value problem
\begin{align}\label{ode:s1}
\begin{cases}
-\bigl(1 + \frac{L}{2}\bigr) \partial_z^2 s_1(z,x,t) + \kappa(s) s_1(z,x,t) = f_1(z,x,t),\\[4pt]
s_1(-\infty,x,t) = 0,\quad s_1(+\infty,x,t) = s_1^+(x,t),
\end{cases}
\end{align}
 has a unique bounded solution if and only if
\begin{align}\label{solve2in}
s_1^+(x,t)=W_A^{-1} \int_{\mathbb{R}} u_-(z) f_1(z,x,t) \, dz.
\end{align}
In this case, the solution is given by
\begin{align}\nonumber
s_1(z,x,t) = &W_A^{-1}\left( u_+(z) \int_{-\infty}^z u_-(w) f_1(w,x,t) \, dw + u_-(z) \int_z^{+\infty} u_+(w) f_1(w,x,t) \, dw \right),
\end{align}
satisfying
$s_1(z,x,t)\in \mathcal{S}_{L,M}(\gamma)$ for all $(x,t)\in \Gamma(\delta)$.
\end{lemma}
\begin{proof}
The lemma can be reduced to Lemma \ref{le:mmodelsolve}, whose proof is provided in Appendix \ref{odesolveappendix}. Indeed, it suffices to introduce the change of variables $z=y/\gamma$ (with  $\gamma=1/\sqrt{1+\frac{2L}{3}}$), and define $A=\frac{-L/6}{1+\frac{2}{3}L}>0$ (which gives  $\frac{\gamma^2}{1+A}=\frac{1}{1+L/2}$). Then, compared with Lemma \ref{le:mmodelsolve}, we can verify that
\begin{align*}
&s(z)=\frac{1}{1+e^{-\gamma z}}=\frac{1}{1+e^{-y}}=\hat{s}(y),\quad\kappa(s)=(1-\hat{s})(1-2\hat{s}),\\
& u_+(z)=\hat{u}_3(y),\quad u_-(z)=\hat{u}_1(y),\quad s_1(z)=u(y)\quad \text{and} \ f_1(z)=f(y).
\end{align*}
\end{proof}

Let $v_\pm(z)$ be smooth solutions of the homogeneous equations
 \begin{align*}
\begin{cases}
-v_+''
+\iota(s)v_+=0,\\
v_+(z)\sim e^{-3z},\;\text{as}\;z\to+\infty,
\end{cases} \qquad \text{ and }\qquad \begin{cases}
-v_-''
+\iota(s)v_-=0,\\
v_-(z)\sim e^{z},\;\text{as}\;z\to-\infty.
\end{cases}
\end{align*}
Denote by $W(v_-,v_+)=v_-'(z)v_+(z)-v_+'(z)v_-(z)$ the Wronskian, and let $W_B=W(v_-,v_+).$

As \cite[Lemma 2.8]{FWZZ2018}, we have
\begin{lemma}\label{lem3}
Assume that $f_2(z,x,t)\in \mathcal{S}_{L,M}(\gamma)$ for all $(x,t)\in \Gamma(\delta)$ with $f_2(+\infty,x,t)=f_{2}^{+}(x,t)$, $f_2(-\infty,x,t)=0$.
The boundary value problem
\begin{align}\label{ode:s2}
\begin{cases}
 -\pa_z^2s_2(z,x,t)+\iota (s)s_2(z,x,t)=f_2(z,x,t),\\
s_2(-\infty,x,t)=0,\ s_2(+\infty,x,t)=s^+_{2}(x,t),
\end{cases}
\end{align}
 has a unique bounded solution if and only if $s^+_{2}(x,t)=\frac19 f^+_{2}(x,t)$. In this case,
the solution is given by:
\begin{align} s_2(z,x,t)=&W_B^{-1}\Big(v_+(z)
\int_{-\infty}^{z}v_-(w)f_2(w,x,t)dw+v_-(z)
\int_{z}^{+\infty}v_+(w)f_2(w,x,t)dw\Big),\label{def:s_beta}
\end{align}
and $s_2(z,x,t)\in \mathcal{S}_{L,M}(\gamma)$ for all $(x,t)\in \Gamma(\delta)$.
\end{lemma}
We refer to \eqref{solve1} and \eqref{solve2in} as the compatibility conditions.
\subsection{Compatibility conditions for $O(\ve^{-1})$ system} 
In this subsection, we analyse equations \eqref{sm0-gamma}-\eqref{sm3-gamma} for $k=1$. The first step is to study the source term $F_1$, cf. \eqref{f1}.
According to the compatibility condition \eqref{solve1}, we derive that
\begin{align}\label{d0evol}
\int_{\mathbb{R}}s'F_1:E^0dz=0\qquad \text{ on }\Gamma.
\end{align}
which is equivalent to
\begin{equation}\label{d0evolution}
    \pa_t d_0-\Big(1+\frac{2L}{3}\Big)\Delta d_0=0,
\end{equation}
by the following lemma.
\begin{lemma}\label{lem:F1E0-1}
It holds on $\Gamma$ that
  \begin{align}
    F_1:E^0&=\frac{2}{3}s'\big(\pa_td_0-(1+\frac{2L}{3})\Delta d_0\big).\label{f1e0}
  \end{align}
\end{lemma}
\begin{Remark}
Since $d_0$ is defined as the signed distance function to $\Gamma_t$, the equation \eqref{d0evolution} is equivalent to the fact that $\Gamma_t$ evolves according to the mean curvature flow.
\end{Remark}
\begin{lemma}\label{lem:F1Ei}
For $(x,t)\in\G,$ it holds that
\begin{align}
F_1:E^1&=-(\bar{h}\cdot l)\Big((4+2L)s'-\frac{L}{3}zs''\Big)-\frac{L}{3}\Big(l\cdot \nabla d_1+d_1(\bar{h}\cdot l)\Big)s'', \label{f1e1}\\
F_1:E^2&=-(\bar{h}\cdot m)\Big((4+2L)s'-\frac{L}{3}zs''\Big)-\frac{L}{3}\Big(m\cdot \nabla d_1+d_1(\bar{h}\cdot m)\Big)s'',\label{f1e2}\\
F_1:E^0&=F_1:E^3=F_1:E^4=0.\label{f1e34}
\end{align}
Moreover, we have for $(x,t)\in\G(\delta)$ that
\begin{align}\label{f1e0:estimate}
|F_1:E^i|\leq &~Cd_0O(e^{-\gamma|z|}), \quad \text{for }i=0,3,4;\\  \label{f1e12:estimate}
|F_1:E^i|\leq &~CO(e^{-\gamma|z|}), \qquad \text{for }i=1,2,
\end{align}
where the constant $C$ is independent of $z$.

\begin{proof}
    The proof is provided in Appendix \ref{f1123proof}.
\end{proof}
\end{lemma}

By compatibility conditions \eqref{solve1} to equation \eqref{sm0-gamma} on $\Gamma(\delta)$ yields
\begin{align*}
&\int_{\mathbb{R}}\Big(\frac{3}{2}F_1:E^0+d_0g_{10}\eta'(z)\Big)s'(z)dz=0.
\end{align*}
Then it suffices to construct $g_{1,0}(x,t)$ as
\begin{equation}\label{g10}
g_{1,0}(x,t)=
\begin{cases}
\begin{array}{ll}
-\f{3}{2d_0}\frac{\int_\BR F_1:E^0s'dz}{\int_{\mathbb{R}}\eta'(z)s'(z)dz},&(x,t)\in \Gamma(\delta)\backslash\Gamma,\\
-\frac32\frac{\pa_{\nu}\big( \int_\BR F_1:E^0s'dz\big)}{\int_{\mathbb{R}}\eta'(z)s'(z)dz},&(x,t)\in \Gamma.
\end{array}
\end{cases}
\end{equation}
Then by Lemma \ref{lem1}, it follows
\begin{align}\label{s10}
s_{1,0}(z,x,t)&=\left(1+\f{2L}{3}\right)^{-1}
\int_{0}^{z}\f{s'(z)}{(s'(y))^2}\int_{y}^{+\infty}
\Big(\frac{3}{2}(F_1:E^0)+\eta'(z)d_0g_{1,0}\Big)s'(w)dw dy.
\end{align}

Since $s$ satisfies
 \begin{align*}
   \bigl(1 + \frac{L}{2}\bigr) s''- \kappa(s) s= -\f{L}{6}s'', \text{ with }s(-\infty)=0,~s(+\infty)=1,
 \end{align*}
we get from \eqref{f1e1}  and \eqref{f1e2} that
\begin{align}
  &s_{1,1}(z,x,t)= (\bar{h}\cdot l)s^*_{1}(z,x,t)+\Big(l\cdot \nabla d_1+d_1(\bar{h}\cdot l)\Big)s,\label{s11trans}\\
   &s_{1,2}(z,x,t)= (\bar{h}\cdot m)s^*_{1}(z,x,t)+\Big(m\cdot \nabla d_1+d_1(\bar{h}\cdot m)\Big)s,\label{s12trans}
\end{align}
where $s^*_{1}(z,x,t)$ is the solution of
\begin{align}\label{s11*}
(1+\frac{L}{2}) \pa_z^2s^*_{1}
-\kappa(s)s^*_{1}
=-\frac{1}{2}\big((4+2L)s'-\frac{L}{3}zs''\big),
\end{align}
with boundary conditions:
\begin{align*}
s_{1}^*(-\infty,x,t)=0,\ s_{1}^*(+\infty,x,t)=W_A^{-1} \int_{\mathbb{R}} u_-(z) \frac{1}{2}\big((4+2L)s'-\frac{L}{3}zs''\big) dz\triangleq s_1^+,
\end{align*}
Thus,
\begin{align}\label{bd:q111}
  q_{1,1}(x,t)=  s_{1,1}(+\infty,x,t)= (\bar{h}\cdot l)s^+_{1}+\big(l\cdot \nabla d_1+d_1(\bar{h}\cdot l)\big),
\end{align}
and
\begin{align}\label{bd:q121}
  q_{1,2}(x,t)=  s_{1,2}(+\infty,x,t)= (\bar{h}\cdot m)s^+_{1}+\big(m\cdot \nabla d_1+d_1(\bar{h}\cdot m)\big).
\end{align}
Equations \eqref{bd:q111}-\eqref{bd:q121} provide the boundary conditions for $q_{1,i}(x,t)\ (i=1,2)$ on $\G$. 
Together with \eqref{q11} and \eqref{q12},
 we can solve $q_{1,i}(x,t)\ (i=1,2)$ in $\Omega_+$ and extend them to $\Gamma(\delta)$ smoothly.

Applying the compatibility condition \eqref{solve2in} with $f_1=-\f12(F_1:E^i)$ for $(x,t)\in\Gamma(\delta)$, we get
\begin{align}\label{Compat-g1al}
q_{1,i}(x,t)-d_0g_{1,i}=&s_{1,i}(+\infty, x,t)-s(+\infty)d_0g_{1,i}\notag\\
=&-\f12W_A^{-1}\int_{\mathbb{R}}u_-(z)
\left(F_1:E^{i}\right)dz,\;\;i=1,2.
\end{align}
Therefore, it suffices to take for $(x,t)\in\Gamma(\delta)$:
\begin{equation}\label{g1i}
g_{1,i}(x,t)=\begin{cases}
\begin{array}{ll}
\frac{1}{d_0}
\Big(q_{1,i}+\f12W_A^{-1}\int_{\mathbb{R}}u_-(z)
\left(F_1:E^{i}\right)dz\Big ),&(x,t)\in \Gamma(\delta)\backslash\Gamma,\\
\pa_{\nu}\Big(q_{1,i}+\f12W_A^{-1}\int_{\mathbb{R}}u_-(z)
\left(F_1:E^{i}\right)dz\Big ),&(x,t)\in \Gamma.
\end{array}
\end{cases}
\end{equation}

For $i=3,4$, Lemma \ref{lem3} gives us
\begin{align}\label{s1j}
s_{1,i}&(z,x,t)=-\f12W_B^{-1}\nonumber\\
&\Big(v_+(z)
\int_{-\infty}^{z}v_-(w)(F_1:E^{i})(w,x,t)dw+v_-(z)
\int_{z}^{+\infty}v_+(w)(F_1:E^{i})(w,x,t)dw\Big).
\end{align}
\begin{Remark}\label{lem:g10-independent}
It follows from the construction that $g_{1,i}(x,t)$ $(i=0,1,2)$ and $s_{1,i}(z,x,t)$ $(i=0,\cdots,4)$ depends only on $d_0$ and $d_1$. Furthermore, from  \eqref{f1e34}, \eqref{s10} and \eqref{s1j}, we conclude that the boundary values $s_{1,i}|_\G$ $(i=0,3,4)$ depend only on $d_0$. Substituting the identity $\int_{\mathbb{R}}s''s' dz=0$ and equation \eqref{F1:fm'} into \eqref{g10}, we find that $g_{1,0}(x,t)$ only depends on $d_0$.

\end{Remark}
Substituting the estimates \eqref{f1e0:estimate}-\eqref{f1e12:estimate} into \eqref{s10}-\eqref{s12trans}  and \eqref{s1j}, we have the following lemma.
\begin{lemma}\label{s1iestimate}
 For $i=0,3,4$, one has $s_{1,i}(z,x,t)\equiv 0$ for $(x,t)\in\Gamma$. Moreover, for  $(x,t)\in\Gamma(\delta)$ and $l=0, 1, 2$ one has
\begin{align}\label{s1-0:estimate}
\left|\pa_z^l\big(s_{1,i}(z,x,t)-s_{1,i}(\pm\infty,x,t)\big)\right|&\leq C d_0O(e^{-\gamma|z|}), \quad \text{ for }i=0,3,4;\\
\left|\pa_z^l\big(s_{1,i}(z,x,t)-s_{1,i}(\pm\infty,x,t)\big)\right|&\leq C O(e^{-\gamma|z|}),\qquad \text{ for }i=1,2, \label{s1-12:estimate}
\end{align}
where the constant $C$ is independent of  $z$.
\end{lemma}
\begin{Remark}  \label{remark:analyitical}
For $i=1,2$, since $s'(z),\ s''(z)$ are real analytical in $z$, we note by the formulation  \eqref{F1:fm'} that $F_1:E^i$ are real analytical in $z$, which is the source term in the ODE \eqref{sm1-gamma}. Therefore, the solutions $s_{1,i}(z,x,t)$ are real analytical in $z$.
\end{Remark}
\subsection{Compatibility conditions for $O(1)$ system}
For $k=2$ in \eqref{Qm},  we have
\begin{equation}\label{q2eq}
\mathcal{L}_{\nu}Q_2-\mathcal{H}_{Q_0}Q_2=F_2+G_2d_0.
\end{equation}
By compatibility condition \eqref{solve1}, we have on $\Gamma$ that
 \begin{align}\label{comd1}
0=\int_{\mathbb{R}}F_2:\partial_zQ_0dz,
 \end{align}
 which gives the equation of $d_1$ on $\G$.
\begin{lemma}\label{lem:d1}
For $(x,t)\in\Gamma,$
$d_1(x,t)$ satisfies
\begin{equation}\label{d1q11ex}
\begin{cases}
\begin{array}{ll}
\pa_td_1-\big(1+\frac{2L}{3}\big)\Delta d_1=\mu_1(x,t)\cdot \nabla d_1(x,t)+\mu_{2}(x,t)d_1(x,t)+ \mu_0(x,t),\ \ \ \ (x, t)\in\Gamma, \\
d_1(x, 0)=0, \ \qquad\qquad\qquad\qquad\qquad\qquad\qquad\qquad\qquad \qquad \qquad \quad\ \ \ \ \quad t=0,
\end{array}
\end{cases}
\end{equation}
where $\mu_i(x,t)$\ ($i=0,1,2$) are smooth functions only depending on $d_0$.
 \end{lemma}
\begin{proof}
See Appendix \ref{lem:d1proof} for the proof.
\end{proof}

Then we can extend $d_1$ from $\G$ to $\G(\delta)$ by the ODE:
\begin{align}\label{d1delta}\nabla d_0\cdot \nabla d_1=0.
\end{align}
\subsection{Compatibility conditions for $O(\ve^{k-1})\ (k\geq 2)$ system}
Recall
\begin{align}
F_{k+1}=&\sum_{\substack{p+q=k\\0\leq p,q\leq k}}(\pa_t-\Delta )d_p\pa_zQ_q-2\nabla d_p\cdot \nabla\pa_zQ_q+(\pa_t -\mathcal{L})Q_{k-1}\nonumber\\
&\quad-\frac{L}{2}\sum_{\substack{p+q=k\\0\leq p,q\leq k}}\mathcal{N}_1(\pa_zQ_p,d_q)
-\frac{L}{2}\sum_{\substack{p+q+r=k+1\\0\leq p\leq k}}\mathcal{N}_2(\pa_z^2Q_p,\nabla d_q,\nabla d_r)\nonumber\\
&\quad-B_{k}-C_{k}+\sum_{\substack{p+q=k+1\\1\leq q\leq k}}G_p(z,x,t)d_q+G_0d_{k+1}-G_{k}z.\label{Fm+1}
\end{align}
It follows that
\begin{lemma}\label{lem:boudaryFk+1Ei}
Let $k\geq 1$. It holds on $\Gamma$ that
\begin{align}
\begin{split}\label{fm:boudaryFk+1Ei}
&F_{k+1}:E^1=-\frac{L}{3}\Big(l\cdot \nabla d_{k+1}+(\bar{h}\cdot l)d_{k+1}\Big)s''+H_{k+1,1},\\
 &F_{k+1}:E^2=-\frac{L}{3}\Big(m\cdot \nabla d_{k+1}+(\bar{h}\cdot m)d_{k+1}\Big)s''+H_{k+1,2}, \\
&F_{k+1}:E^i\quad \text{only}\;\text{depend} \;\text{on}\; d_p,s_{p,j}\ (0\leq p\leq k,\;0\leq j\leq 4) \text{ for }i=0,3,4,
\end{split}
\end{align}
where $H_{k+1,i}$ for $i=1,2$ depends only on $d_p$ and $s_{p,j}\ (0\leq p\leq k,\;0\leq j\leq 4)$.
Assuming $d_p\ (p\leq k+1)$ are known,  for $i=0,\cdots,4$ it holds on $\Gamma $ that
\begin{align}\label{fke12estimate}
|F_{k+1}:E^i|\leq CO(e^{-\gamma|z|}),
\end{align}
where $C$ is independent of $z$.
\end{lemma}
\begin{proof}
See the proof of \eqref{fm:boudaryFk+1Ei} in the Appendix \ref{f1123proof}. \eqref{fke12estimate} follows by the fact in the definition \eqref{Fm+1} that $F_{k+1}|_{\G}$ only depends on $s_{l,i}'|_{\G}(z,x,t)$ for $l\leq k$, $i=0,1,2,3,4$ and the estimates in \eqref{skiestimate}.
\end{proof}

According to the  compatibility condition \eqref{solve1} in $\Gamma(\delta)$, we have
\begin{align*}
&\int_{\mathbb{R}}\Big(\frac{3}{2} F_{k+1}:E^0+g_{k+1,0}d_0\eta'(z)\Big)s'(z)dz=0,
\end{align*}
which gives
\begin{equation}\label{gm+10}
g_{k+1,0}(x,t)=\begin{cases}
\begin{array}{ll}
-\f{3}{2d_0}\frac{\int_\BR F_{k+1}:E^0s' dz}{\int_{\mathbb{R}}\eta'(z)s'(z)dz},&(x,t)\in \Gamma(\delta)\backslash\Gamma,\\
-\frac32\frac{\pa_{\nu}\big( \int_\BR F_{k+1}:E^0s' dz\big)}{\int_{\mathbb{R}}\eta'(z)s'(z)dz},&(x,t)\in \Gamma.
\end{array}
\end{cases}
\end{equation}

It follows by Lemma \ref{lem1} that
\begin{align}\label{sm+10}
s_{k+1, 0}(z, x, t)=&\left(1+\f{2L}{3}\right)^{-1}\notag\\
&\int_{0}^{z}\f{s'(z)}{(s'(y))^2}\int_{y}^{+\infty}
s'(w)
\Big(\frac{3}{2}(F_{k+1}:E^0)+\eta'(z)d_0g_{{k+1}, 0}\Big)\ud w \ud y.
\end{align}

Recall $s$ satisfies
 \begin{align*}
   \bigl(1 + \frac{L}{2}\bigr) s''- \kappa(s) s= -\f{L}{6}s'', \text{ with }s(-\infty)=0,~s(+\infty)=1,
 \end{align*}
we get from \eqref{fm:boudaryFk+1Ei}  that
\begin{align}
  &s_{k+1,1}(z,x,t)= s^*_{k+1,1}(z,x,t)+\Big(l\cdot \nabla d_{k+1}+d_{k+1}(\bar{h}\cdot l)\Big)s,\label{sk+11trans}\\
   &s_{k+1,2}(z,x,t)= s^*_{k+1,2}(z,x,t)+\Big(m\cdot \nabla d_{k+1}+d_{k+1}(\bar{h}\cdot m)\Big)s,\label{sk+12trans}
\end{align}
where $s^*_{k+1,1}(z,x,t)$ and   $s^*_{k+1,2}(z,x,t)$ being the solution of
\begin{align}
(1+\frac{L}{2}) \pa_z^2s^*_{k+1,1}
-\kappa(s)s^*_{k+1,1}
=\frac{1}{2}H_{k+1,1},\label{sk+11*}\\
(1+\frac{L}{2}) \pa_z^2s^*_{k+1,2}
-\kappa(s)s^*_{k+1,2}
=\frac{1}{2}H_{k+1,2},\label{sk+12*}
\end{align}
with boundary conditions:
\begin{align*}
&s_{k+1,1}^*(-\infty,x,t)=0,\ s_{k+1,1}^*(+\infty,x,t)=W_A^{-1} \int_{\mathbb{R}} u_-(z) \frac{1}{2} H_{k+1,1} dz\triangleq s_{k+1,1}^+,\\
&s_{k+1,2}^*(-\infty,x,t)=0,\ s_{k+1,2}^*(+\infty,x,t)=W_A^{-1} \int_{\mathbb{R}} u_-(z) \frac{1}{2} H_{k+1,2} dz\triangleq s_{k+1,2}^+.
\end{align*}
Thus,
\begin{align}\label{bd:qk+11}
  q_{k+1,1}(x,t)=  s_{k+1,1}(+\infty,x,t)= s^+_{k+1,1}+\big(l\cdot \nabla d_{k+1}+d_{k+1}(\bar{h}\cdot l)\big).
\end{align}
Similarly, we have
\begin{align}\label{bd:qk+12}
  q_{k+1,2}(x,t)=  s_{k+1,2}(+\infty,x,t)= s^+_{k+1,2}+\big(m\cdot \nabla d_{k+1}+d_{k+1}(\bar{h}\cdot m)\big).
\end{align}
Equations \eqref{bd:qk+11}-\eqref{bd:qk+12} provide the boundary conditions for $q_{k+1,i}(x,t)\ (i=1,2)$ on $\G$. Together with \eqref{qm+1}, we can solve $q_{k+1,i}(x,t)\ (i=1,2)$ in $\Omega_+$ and extend them to $\Gamma(\delta)$ smoothly.
Similar as the compatibility condition \eqref{Compat-g1al}, we deduce by \eqref{solve2in} that
\begin{equation}\label{gk+1i}
g_{k+1,i}(x,t)=\begin{cases}
\frac{1}{d_0}
\big\{q_{k+1,i}
+\f12W_A^{-1}
\int_{\mathbb{R}}u_-(z)
\left(F_{k+1}:E^{i}\right)\ud z\big\},&(x,t)\in \Gamma(\delta)\backslash\Gamma,\\
\pa_{\nu}\big\{q_{k+1,i}+
\f12W_A^{-1}\int_{\mathbb{R}}u_-(z)
\left(F_{k+1}:E^{i}\right)\ud z\big\},&(x,t)\in \Gamma.
\end{cases}
\end{equation}

For $i=3,4$, we deduce from  Lemma \ref{lem3} that
\begin{align}
s_{k+1, i}(z, x, t)=&-\frac{1}{2}W_B^{-1}
\Big\{v_+(z)
\int_{-\infty}^{z}v_-(w)(F_{k+1}:E^{i})(w, x, t)\ud w\notag\\
&+v_-(z)
\int_{z}^{+\infty}v_+(w)(F_{k+1}:E^{i})(w, x, t)\ud w\Big\}. \label{sm+13}
\end{align}
In summary, we can obtain $$Q_{k+1}(z,x,t)=\sum_{i=0}^4s_{k+1,i}(z,x,t)E^i(x,t).$$
It follows by   the above solution and Lemma \ref{lem1}, Lemma \ref{lem2}, Lemma \ref{lem3}, Lemma \ref{lem:boudaryFk+1Ei} and \eqref{fke12estimate} that
\begin{lemma}\label{skiestimate}
For $k\geq1$, $s_{k+1,i}|_\G$ ($i=0,3,4$) only depend on  $d_l(l\leq k)$, and it holds on $\Gamma$ that
\begin{align}
\label{est:s_{k,i}^(k)}
&
\left|\pa_z^ms_{k+1,i}(z,x,t)\right|\leq CO(e^{-\gamma|z|}), \text{ for } m=0,1,2.
\end{align}
For $i=1,2$, it holds on $\Gamma$ that
\begin{align}
\left|\pa_z^m\Big(s_{k+1,i}(z,x,t)-s_{k+1,i}(\pm\infty,x,t)\Big)\right|&\leq CO(e^{-\gamma|z|}), \text{ for } m=0,1,2.
\end{align}
Here the constants $C$ are independent of  $z$.
\end{lemma}
\begin{Remark}\label{lem:gm+1-independent}
It follows by the construction that $g_{k+1,i}(x,t)\ (i=0,1,2)$ and $Q_{k+1}$ only depend on $d_l$ for $l\leq k+1$.
In particular, using $\int_{\mathbb{R}}s''s'dz=0$ and \eqref{Fm+1}, one gets that $g_{k+1,0}(x,t)$ is independent of $d_{k+1}$, cf. \eqref{gm+10}.
\end{Remark}
\subsection{Compatibility conditions for $O(\ve^{k})$ system}
For \eqref{Qm}, the $O(\ve^k)$ system can be written as
\begin{align}
\mathcal{L}_{\nu}Q_{k+2}-\mathcal{H}_{Q_0}Q_{k+2}=F_{k+2}+G_{k+2}d_0.\label{inm}
\end{align}
On $\Gamma,$ by  the compatibility condition
$$\int_{\mathbb{R}}F_{k+2}:\pa_zQ_0 \ud z=0,$$
we obtain the evolution equation of $d_{k+1}$:
\begin{lemma}\label{lem:dk+1}
Let $k+1\geq 2$. For $(x,t)\in\Gamma,$
$d_{k+1}(x,t)$ satisfies
\begin{equation}\label{dk+1q11ex}
\begin{cases}
\begin{array}{ll}
\pa_td_{k+1}-\big(1+\frac{2L}{3}\big)\Delta d_{k+1}\\
\qquad \qquad =\mu_1^{k+1}(x,t)\cdot \nabla d_{k+1}(x,t)+\mu_{2}^{k+1}(x,t)d_{k+1}(x,t)+ \mu_0^{k+1}(x,t),\ \ (x, t)\in\Gamma, \\
d_{k+1}(x, 0)=0\ \qquad\qquad\qquad\qquad\qquad\qquad\qquad\qquad\quad\quad\qquad\qquad \qquad \ \ x\in\G_0,
\end{array}
\end{cases}
\end{equation}
where $\mu_i^{k+1}(x,t)$($i=0,1,2$) are smooth functions depending only on $d_l,s_{l,0}|_{\G},s_{l,1}|_{\G},s_{l,2}|_{\G}(l\leq k)$.
\end{lemma}
\begin{Remark}
Combining Lemma \ref{skiestimate}, we know $\mu_i^{k+1}(x,t)$ ($i=0,1,2$)
only depends on $d_l\ (l\leq k)$.
\end{Remark}
The function $d_{k+1}$ can be extended from $\G$ to $\G(\delta)$ via the following ODE:
 \begin{equation}\label{dm+1delta}
     \nabla d_0\cdot \nabla d_{k+1}=-\f12\sum_{\substack{1\leq j\leq k}}\nabla d_j\cdot\nabla d_{k+1-j}.
 \end{equation}
\section{Solving outer/inner expansion}\label{section4}
Let $k\geq 0$. The goal of this section is to show the procedure determining $Q_k$, $G_k$, $d_k$ in $\Gamma(\delta)$ and $Q_k^{\pm}$ in $U_\pm$. For convenience, we introduce
$$\mathcal{V}^k=\{Q_k,Q_k^\pm,d_k,G_k\}\text{ for }k\geq 0.$$
\subsection{Solving $\mathcal{V}^0$}
\subsubsection{Determining $Q_0^{\pm}$ in $U_\pm$.}
$Q_0^{\pm}$ are determined in \eqref{Q0}.
\subsubsection{Determining  $Q_0$ in $\Gamma(\delta)$}
\eqref{q0} determines $Q_0$ in $\Gamma(\delta)$.
\subsubsection{Determining  $d_0$ in $\Gamma(\delta)$}
  Equation \eqref{d0evolution} implies that $\Gamma_t$ evolves by mean curvature flow. In the neighborhood $\Gamma(\delta)$, $d_0$ is determined by the following system:
  \begin{equation}\label{d0near}
\left\{
\begin{array}{ll}
\pa_t d_0-\left(1+\frac{2L}{3}\right)\Delta d_0=0,&(x,t)\in\Gamma,\\
d_0=0,&(x,t)\in\Gamma,\\
\nabla d_0\cdot \nabla d_0=1,&(x,t)\in\Gamma(\delta).
\end{array}
\right.
\end{equation}
\subsubsection{Determining $G_0$ in $\Gamma(\delta)$}
Substituting $d_0$ into \eqref{G0}, one can get $G_0$.
\subsection{Solving $\mathcal{V}^1$}
\subsubsection{Determining  $Q_{1}^{\pm}$ in $U_\pm$ and $d_1$ on $\G$.}
$Q_1^-$ is determined in \eqref{def:Q_m^-}.

For $Q_1^+$, thanks to \eqref{Q1}, it is sufficient to determine $q_{1,1}$ and $q_{1,2}$. We can find the parabolic system  \eqref{q11} and \eqref{q12} in $ U_+$. The boundary condition $q_{1,i}(x,t)|_{\Gamma}\ (i=1,2)$ of the system is given in \eqref{bd:q111}, \eqref{bd:q121} which only depends on $d_0|_{\Gamma}$ and $d_1|_{\Gamma}$. Let us determine $d_1|_\G$ and $q_{1,1}$, $q_{1,2}$.

 By Lemma \ref{lem:d1},  we have
\begin{align*}\label{d1q11ex'}
\begin{cases}
\begin{array}{ll}
\pa_td_1-\big(1+\frac{2L}{3}\big)\Delta d_1=\mu_1(x,t)\cdot \nabla d_1(x,t)+\mu_{2}(x,t)d_1(x,t)+ \mu_0(x,t)\ \ \ \ (x, t)\in\Gamma, \\
d_1(x, 0)=0, \ \qquad\qquad\qquad\qquad\qquad\qquad\quad\;\qquad\qquad\qquad\qquad \qquad \ \ \ \ \quad t=0,
\end{array}
\end{cases}
\end{align*}
where $\mu_i(x,t)$ ($i=0,1,2$) are smooth functions that only depend on $d_0$. We can find  $d_1|_\G$ by the above system.

By \eqref{q11}, \eqref{q12}, \eqref{bd:q111} and \eqref{bd:q121}, we have the following system:
\begin{equation}\label{q1q12d1}
\begin{cases}
\begin{array}{ll}
\pa_tq_{1, 1}-\Delta q_{1, 1}-L/2(\pa_{ik}q_{1, 1}n_kn_i+\pa_{ik}q_{1, 1}l_kl_i+\pa_{ik}q_{1, 2}m_kl_i)=r_{1,1}, &(x,t)\in  U_+, \\
\pa_tq_{1, 2}-\Delta q_{1, 2}-L/2(\pa_{ik}q_{1, 1}l_km_i+\pa_{ik}q_{1, 2}n_kn_i+\pa_{ik}q_{1, 2}m_km_i)=r_{1,2}, &(x,t)\in U_+, \\
q_{1, 1}(x, t)|_{t=0}=q_{1, 1}(x, 0), \ q_{1, 2}(x, t)|_{t=0}=q_{1, 2}(x, 0), &x\in \Omega^+_0,
\end{array}
\end{cases}
\end{equation}
with the following Dirichlet boundary conditions on $\G$:
\begin{align}\label{bd:q11}
  q_{1,1}(x,t)=  s_{1,1}(+\infty,x,t)= (\bar{h}\cdot l)s^+_{1}+\big(l\cdot \nabla d_1+d_1(\bar{h}\cdot l)\big),
\end{align}
and
\begin{align}\label{bd:q12}
  q_{1,2}(x,t)= (\bar{h}\cdot m)s^+_{1}+\big(m\cdot \nabla d_1+d_1(\bar{h}\cdot m)\big),
\end{align}
where $r_i$ is a linear function for $q_{1,1},\ q_{1,2},\ \nabla q_{1,1}$ and $\nabla q_{1,2}$.
Then we obtain $q_{1,1}$ and $q_{1,2}$ in $ U_+$

\subsubsection{Determining  $d_1$, $Q_1$ and $G_1$ in $\G(\delta)$}
We can easily get $d_1$ in $\G(\delta)$ by the ODE:
\begin{align*}
\nabla d_0\cdot \nabla d_1=0.
\end{align*}
Once $d_0$, $Q_0$, $G_0$, $d_1$ and $q_{1,i}\ (i=1,2)$ are determined, $G_1$ is constructed via \eqref{Gm} using $g_{1,0}$  from \eqref{g10} and $g_{1,i}\ (i=1,2)$ from \eqref{g1i}. We can then solve for $s_{1,i}$ ($i=0,\cdots,4$) using equations \eqref{sm0-gamma}-\eqref{sm3-gamma} with the explicit formulas provided in \eqref{s10}, \eqref{s11trans}, \eqref{s12trans} and \eqref{s1j}.  Remark \ref{lem:g10-independent} ensures that this construction is independent of any unknown functions.
\subsection{Solving $\mathcal{V}^{k+1}\ (k\geq 1)$}
Given that $\mathcal{V}^p$ for $p \leq k$ are known, we now aim to determine $\mathcal{V}^{k+1}$. This allows the expansion system to be solved via an induction argument.
\subsubsection{Determining $d_{k+1}$ on $\G$ and $Q_{k+1}^{\pm}$ in $U_\pm$ .}
 We first solve the  following linear equation to obtain $d_{k+1}|_{\G}$:
  \begin{align*}
\begin{cases}
\begin{array}{ll}
\pa_td_{k+1}-\big(1+\frac{2L}{3}\big)\Delta d_{k+1}\\
\qquad\qquad =\mu_1^{k+1}(x,t)\cdot \nabla d_{k+1}(x,t)+\mu_{2}^{k+1}(x,t)d_{k+1}(x,t)+ \mu_0^{k+1}(x,t),\ \ (x, t)\in\Gamma, \\
d_{k+1}(x, 0)=0,\ \qquad\qquad\qquad\qquad\qquad\qquad\qquad\qquad\quad\quad\qquad\qquad\qquad \ \ x\in\G_0,
\end{array}
\end{cases}
\end{align*}
where $\mu_i^{k+1}(x,t)$ ($i=0,1,2$) are smooth functions only depend on $d_l,s_{l,0}|_{\G},s_{l,1}|_{\G},s_{l,2}|_{\G}(l\leq k)$. $Q_{k+1}^-$ is given by \eqref{def:Q_m^-}. For $Q_{k+1}^\bot$, we have
 \begin{equation}\label{Qm+1out}
Q_{k+1}^\bot=\mathcal{H}_{+}^{-1}\Big(-\mathcal{B}_+Q_{k}^\top-2\mathcal{C}_+Q_{k}^\top+D_{k}\Big)
:=\mathcal{H}_{+}^{-1}A_{k},
 \end{equation}
  where $A_{k}$ depends on $Q_0^+,\dots,Q_{k-1}^+$ and $Q_{k}^\bot$, which is known by the induction. For $Q_{k+1}^\top$, combining parabolic systems \eqref{qk+11}, \eqref{qk+12} and the boundary
 conditions \eqref{bd:qk+11}, \eqref{bd:qk+12},
 we can obtain $Q_{k+1}^\top$ in $ U_+$.
\subsubsection{Determining  $d_{k+1}$, $Q_{k+1}$ and $G_{k+1}$ in $\G(\delta)$}
After $d_{k+1}|_\G$ is determined from the equation \eqref{dk+1q11ex}, we extend $d_{k+1}$ from $\G$ to $\G(\delta)$ via the ODE:
\begin{align}\label{dk+1extend}
\nabla d_0\cdot \nabla d_{k+1}=-\f12\sum_{\substack{1\leq j\leq k}}\nabla d_j\cdot\nabla d_{k+1-j}\ \ \ \ \ \ &\text{in }\Gamma(\delta).
\end{align}
Since $g_{k+1,0}$ is independent of $d_{k+1}$, it can be obtained directly from \eqref{gm+10}. Once $d_{k+1}$ in $\G(\delta)$ and $Q_{k+1}^\top$ in $ U_+$ are determined, we can write $F_{k+1}$. The functions $g_{k+1,i}$ ($i=1,2$) are then given by \eqref{gk+1i}. Subsequently, we solve for $s_{k+1,i}$ ($i=0,\cdots, 4$) using \eqref{sm0-gamma}-\eqref{sm3-gamma} with the explicit expressions provided in \eqref{sm+10}, \eqref{sk+11trans}, \eqref{sk+12trans} and \eqref{sm+13}. Finally, $\tilde{g}_{k+1}$ is defined via \eqref{def:hm}.  Remark \ref{lem:gm+1-independent} ensures that this entire construction depends only on known functions.
\begin{figure}\small
\begin{tikzpicture}[man/.style={rectangle,draw,fill=blue!20},
  woman/.style={rectangle,draw,fill=red!20,rounded corners=.8ex},
 VIP/.style={rectangle,very thick, draw}, fill=white]
\draw [dotted](-3,0.25) to (-3,-4.8);
\draw[dotted](-3,-5.2) to (-3,-7.4);
\draw[dotted] (-3,-7.8) to (-3,-12.8);
\draw[dotted] (-3,-13.2) to (-3,-15.6);
\draw[dotted] (-3,-16) to (-3,-1*17.3);

\node (Out) at (0,0){Inner expansion};
\node (Inn) at (-6,0){Outer expansion};

\node[man] (A0) at (-0.7,-1){$Q_0,G_0, d_0 $ in ${\Gamma(\delta)}$};
\node[man] (A0-out) at (-7.2,-1) {$Q^+_{0}$ in ${\Omega}^+$, $Q^-_k$ in $ U_-$};

\node[man, text width=5.2cm, align=left] (F1) at (0.2,-2.5)  {$s_{1,i}|_\G$ for $i=0,3,4$: \eqref{s10}, \eqref{s1j}\\
$s_{1}^*$:  \eqref{s11*}\\and $g_{1,0}|_{\G(\delta)}$: \eqref{g10}};

\node[man] (d11) at (-0.7,-4) {Equations for $d_1|_\G$};
\node[woman] at (2,-4) {\eqref{d1q11ex}};
\node[man] (d12gamma) at (-0.7,-5)  {$d_1|_\G$};
\node[man] (d12) at (-0.7,-6)  {$d_1$ in $\Gamma(\delta)$};
\node[man] (P1star) at (-0.7,-7) {$F_1$ in $\G(\delta)$};

\node[man] (P1i) at (-0.7,-1*8) {$Q_1, G_1$ in $\Gamma(\delta)$};
\node[woman] at (2.7,-8) {\eqref{g10}-\eqref{s12trans},\eqref{g1i},\eqref{s1j}};
\node[woman] at (2.7,-1*1) {\eqref{q0}, \eqref{G0}, \eqref{d0near}};

\node[woman] at (1,-6) {\eqref{d1delta}};

\node[woman] at (1,-7) {\eqref{f1}};

\node[man] (M1) at (-7.2,-1*2) {$Q_1^\bot=0$ in $ U_+$};
\node[man] (V1) at (-6,-7)  {$Q_1^\top$ in $ U_+$};
\node[woman] at (-4.35,-1*1) {\eqref{Q0}, \eqref{def:Q_m^-}};
\node[woman] at (-9.2,-1*2) {(\ref{Q1})};
\node[man] (Eq1) at (-7.2,-1*3)  {Equations for $Q_1^{\top}|_{ U_+}$: \eqref{q1q12d1}};

\node[man] (Eqbd1) at (-6,-5){Boundary conditions for $Q_1^\top$};
\node[woman] at (-4.5,-5.7){\eqref{bd:q111},\eqref{bd:q121}};

\node (stepk) at (-0.4,-9)  {\underline{$\{d_p, ~Q_p, ~G_p| p\leq k \}$ are solved}};

\node (stepk) at (-5.5,-1*8-1)  {\underline{$\{Q_p^\pm|p\le k \}$ are solved}};

\node[man,text width=5.6cm, align=left] (F2) at (0.2,-10.5) {$s_{k+1,i}|_\G$ for $i=0,3,4$: \eqref{sm+10}, \eqref{sm+13}\\
$s_{k+1,1}^*,\ s_{k+1,2}^*$: \eqref{sk+11*}, \eqref{sk+12*} \\
and $g_{k+1,0}|_{\G(\delta)}$: \eqref{gm+10}};
\node[man] (d21) at (-0.7,-12)  {Equations for $d_{k+1}|_\G$};
\node[man] (d22gamma) at (-0.7,-13)  {$d_{k+1}|_\G$};
\node[man] (d22) at (-0.7,-14)  {$d_{k+1}$ in $\Gamma(\delta)$};
\node[man] (F22) at (-0.7,-15)  {$F_{k+1}$ in $\G(\delta)$};
\node[man] (P2i) at (-0.7,-16)  {$Q_{k+1},G_{k+1}$ in $\Gamma(\delta)$};
\node[woman] at (3.2,-16)  {\eqref{gm+10}-\eqref{sk+12trans},\eqref{gk+1i},\eqref{sm+13}};
\node[woman] at (2,-12)  {\eqref{dk+1q11ex}};
\node[woman] at (1.2,-15)  {\eqref{Fm+1}};
\node[woman] at (1.2,-14)  {\eqref{dm+1delta} };

\node[man] (Mk) at (-7.3,-1*10)  {$Q_{k+1}^\bot$ in $ U_+$};
\node[man] (Vk) at (-5.7,-15.6)  {$Q_{k+1}^\top$ in $ U_+$};

\node[man] (Eqk) at (-7.3,-11){Equations for $Q_{k+1}^\top|_{ U_+}$: \eqref{qk+11},\eqref{qk+12}};
\node[man] (Eqbdk) at (-5.7,-13)  {Boundary conditions for
 $Q_{k+1}^\top$};
\node[woman] at (-4.5,-13.8){\eqref{bd:qk+11},\eqref{bd:qk+12}};
\node[woman] at (-9,-1*10)  {(\ref{Qm+1=H-1})};

\node (stepk) at (-0.3,-17)  {\underline{$\{d_{k+1},\ G_{k+1},\ Q_{k+1}\}$ are solved}};
\node (stepk1) at (-5.5,-17)    {\underline{$\{Q_{k+1}^\pm\}$ are solved}};

\draw [->](A0) to (F1);
\draw [->](d12) to (P1star) ;
\draw [->](F1) to (d11);
\draw [->](d11) to (d12gamma);
\draw [->](d12gamma) to (d12);
\draw [->](V1) to (P1i);
\draw [->](A0-out) to  (M1);
\draw [->](d12gamma) to (Eqbd1);
\draw [->](M1) to (Eq1);
\draw [->](P1star) to (P1i);
\draw [->](V1) to (P1i);
\draw [->](Eqbd1) to (V1);
\draw [->](-9,-3.3) |- (V1);

\draw [->](F2) to (d21);
\draw [->](d21) to (d22gamma);
\draw [->](d22gamma) to (d22);
\draw [->](d22) to (F22);
\draw [->](F22) to (P2i);
\draw [->](Vk) to (P2i);
\draw [->](Eqbdk) to (Vk);
\draw [->](-9,-11.3) |-(Vk);
\draw [->](Mk) to (Eqk);
\draw [->](d22gamma) to (Eqbdk);

\draw [very thick, dotted](-0.7,-8.4) to (-0.7,-8.82);
\draw [very thick, dotted](-5.8, -7.4) to (-5.8,-8.82);
\draw [very thick, dotted](-0.7,-17.2) to (-0.7,-17.7);
\draw [very thick, dotted](-5.8, -17.2) to (-5.8, -17.7);
\end{tikzpicture}
\caption{The whole procedure to solve the outer and inner expansion systems}
\label{fig:expansion}
\end{figure}
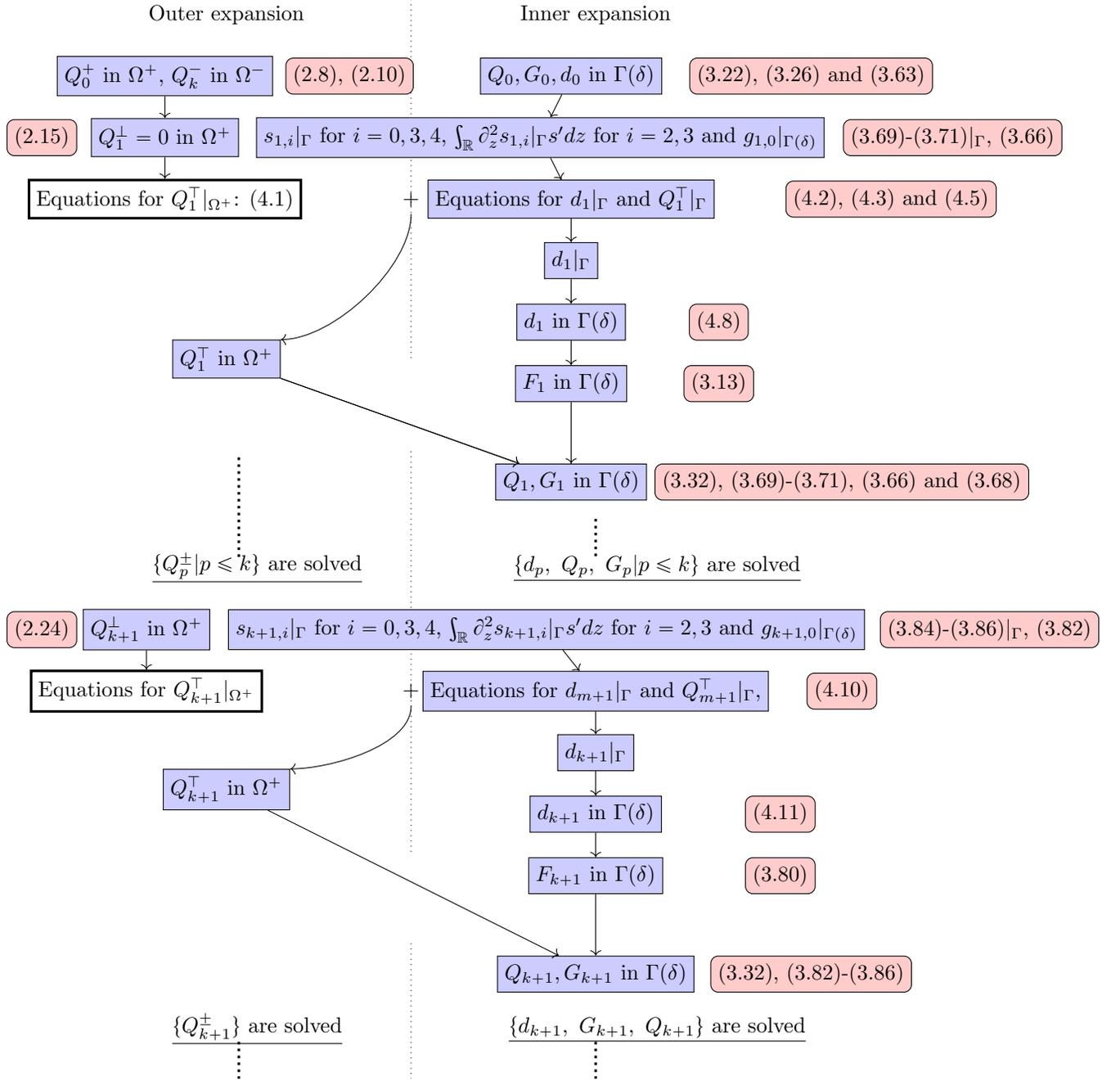
\section{Construction of approximate solutions}\label{section5}
In this section, we glue together the inner expansion and outer expansion.
The outer expansion has the form:
$$Q^K_\pm(x,t)=Q_0^\pm(x,t)+\ve Q_1^\pm(x,t)+\ve^2Q_2^\pm(x,t)+\cdots+\ve^KQ_K^\pm(x,t).$$
We define
$$Q_o^K=Q_+^K\chi _{ U_+}+Q_-^K\chi_{ U_-},$$
which gives
\begin{align}
\partial_t&Q^K_{o}-\mathcal{L}(Q^K_{o})+\frac{1}{\ve^2}f(Q^K_{o})\nonumber\\
&=\sum_{k=0}^{K-2}\ve^{k}\big(\partial_tQ_k^\pm-\mathcal{L}Q_k^\pm+\mathcal{H}_{Q_0^\pm}Q_{k+2}^\pm+C_{k+1}^\pm+B_{k+1}^\pm\big)+O(\ve^{K-1})\nonumber\\
&=O(\ve^{K-1}).\label{outer}
\end{align}
In $\Gamma(\delta)$, let
\begin{align}
Q^K_{in}(x,t)=Q_0(\frac{d^K}{\ve},x,t)+\ve Q_1(\frac{d^K}{\ve},x,t)+\cdots+\ve^KQ_K(\frac{d^K}{\ve},x,t),
\end{align}
with
\begin{align}\label{nabladK-1}
 d^K(x,t)= d_0+\sum_{i=1}^K\ve^{i}d_i= d_0+O(\ve).
\end{align}
Then  we have by \eqref{didj} that
\begin{align}\label{nabladK}
|\nabla d^K|^2=1+\sum_{\substack{1\leq i,j\leq K\\i+j\geq K+1}}\ve^{i+j}\nabla d_j\nabla d_i=1+O(\ve^{K+1}).
\end{align}
Upon comparing the coefficients of $\ve^k$, it is apparent that
\begin{align}
\partial_tQ^K_{in}&-\mathcal{L}(Q^K_{in})+\frac{1}{\ve^2}f(Q^K_{in})\nonumber\\
=&\partial_tQ^K_{in}-\mathcal{L}(Q^K_{in})+\frac{1}{\ve^2}f(Q^K_{in})+\Big(\sum_{k=0}^{K}\ve^kd_k-\ve z\Big)\sum_{k=0}^{K}\ve^{k-2}G_k\big|_{z=\frac{d^K}{\ve}}\nonumber\\
=&O(\ve^{K-1}).\label{inner}
\end{align}
Due to the matching condition \eqref{mathingcondition}, we have
$$ \lim_{z\rightarrow  \pm\infty}|\pa_t^m\pa_x^n\pa_z^l(Q_I^K-Q_O^K)|=O(e^{-\gamma|z|}),\ \text{for}\ (x,t)\in\Gamma(\delta),\ m,n,l\geq 0.$$
Define
\begin{align}\label{def:QK}
Q^K=\eta(\frac{d_0}{\delta})Q_{in}^K+(1-\eta(\frac{d_0}{\delta}))Q^K_o,
\end{align}
where $\eta$ is a smooth cut-off function with
\begin{align}\label{eta-function}
 \text{$\eta(y)=1$ for $|y|\le 1/2$ and $\eta(y)=0$ for $|y|\ge 1$}.
\end{align}
Then it holds that
\begin{align}
\partial_tQ^K-\mathcal{L} Q^K+\ve^{-2}f(Q^K)=\mathcal{R}^K,
\end{align}
where $\mathcal{R}^K=O(\ve^{K-1})$. The proof of Theorem \ref{th1} is completed.
\section{Spectral lower bound estimate for the linearized operator}\label{section6}
This section is devoted to the proof of Theorem \ref{th:uA}, specifically inequality \eqref{mainineq}. It suffices to consider sufficiently small $\ve > 0$. The proof proceeds in six steps, detailed in Sections \ref{6.1}-\ref{6.6}:
\begin{itemize}
 \item [Step 1.] Applying the divergence-curl decomposition to the gradient of the $Q$-tensor field yields:
 \begin{align*}
|\nabla Q|^2 = \frac{3}{2}|\nabla\cdot Q|^2 + \frac{1}{4}|\mathcal{T}(Q)|^2 + (\partial_kQ_{li}\partial_lQ_{ki} - \partial_kQ_{ki}\partial_lQ_{li}),
\end{align*}
where
\begin{equation*}
\mathcal{T}(Q) = \big(T_{ij}(Q)\big)_{1\leq i, j\leq 3}, \quad T_{ij}(Q) = \varepsilon^{ikl}\partial_kQ_{lj} + \varepsilon^{jkl}\partial_kQ_{li}.
\end{equation*}
This demonstrates that inequality \eqref{mainineq} is equivalent to \eqref{omegaestimate}.
\item [Step 2.] By applying the coordinate transformation \eqref{transformation} and the basis decomposition \eqref{Qdecom}, the problem is reduced to establishing scalar inequalities on a one-dimensional interval. Specifically, Lemma \ref{neartheinterface} further reduces \eqref{omegaestimate} to \eqref{ineqmain} in the domain $\Gamma_t^K(\delta/4)$. Within this framework, Lemma \ref{tqestimate} provides an estimate for $\frac{L}{6}|\mathcal{T}(Q)|^2$. It therefore remains to prove \eqref{prove}. Finally, the term $J$ is removed from \eqref{prove}, thereby reducing the problem to verifying \eqref{prove-Linear}--\eqref{prove-Correction}.
 \item [Step 3.] We establish coercive estimates for the scalar linearized operators $\mathcal{G}_j(p)$ ($j=0,1$) in Lemmas \ref{le:0} and \ref{le:1}, and obtain $L^\infty$ estimates at the endpoints in Lemmas \ref{leminfty} and \ref{le:endpoints}. Furthermore, we derive second spectral estimates for $\mathcal{G}_j(p)$ ($j=0,1$) in Lemmas \ref{secondeigenvalue} and \ref{g1second}.
 \item [Step 4.] We establish estimates for the convolution terms in \eqref{prove-Linear} using Lemmas \ref{lem:tq1-1}--\ref{lem:tq1-3}.
 \item [Step 5.] We bound the cross terms in \eqref{prove-Cross} using strong coercive estimates, as established in Lemmas \ref{cross1}--\ref{cross3}.
 \item [Step 6.] We establish estimates for the correction terms in \eqref{prove-Correction}, with the detailed bounds provided in Lemma \ref{correction}.
 \end{itemize}
 \subsection{Divergence-curl decomposition for gradient of $Q$-tensor field}\label{6.1}
\begin{lemma}\label{divcurl}For any symmetric traceless tensor $Q$, one has
\begin{align}\label{nablaqdecom}
|\nabla Q|^2=\frac{3}{2}|\nabla\cdot Q|^2+\frac{1}{4}|\mathcal{T}(Q)|^2+(\pa_kQ_{li}\pa_lQ_{ki}-\pa_kQ_{ki}\pa_lQ_{li}),
\end{align}
where
\begin{equation}\label{tq2}
\mathcal{T}(Q)=\big(T_{ij}(Q)\big)_{1\leq i,j\leq 3}\;\;\;\text{ with }\;\;T_{ij}(Q)=\ve^{ikl}\pa_kQ_{lj}+\ve^{jkl}\pa_kQ_{li},
\end{equation}
where $\ve^{ijk}$  is the Levi-Civita symbol denoted in \eqref{levisivita}.
\begin{proof}
For a symmetric traceless tensor $Q$, direct calculation gives that
\begin{align*}
|\mathcal{T}(Q)|^2&=|\ve^{ikl}\pa_kQ_{lj}+\ve^{jkl}\pa_kQ_{li}|^2\\
&=\ve^{ikl}\ve^{imn}\pa_kQ_{lj}\pa_mQ_{nj}+\ve^{jkl}\ve^{jmn}\pa_kQ_{li}\pa_mQ_{ni}+2\ve^{ikl}\ve^{jmn}\pa_kQ_{lj}\pa_mQ_{ni}\\
&=2\ve^{ikl}\ve^{imn}\pa_kQ_{lj}\pa_mQ_{nj}+2\ve^{ikl}\ve^{jmn}\pa_kQ_{lj}\pa_mQ_{ni}.
\end{align*}
As $\ve^{ikl}\ve^{imn}=\delta_{km}\delta_{ln}-\delta_{kn}\delta_{lm}$,  one has
\begin{align*}
\ve^{ikl}\ve^{imn}\pa_kQ_{lj}\pa_mQ_{nj}&=(\delta_{km}\delta_{ln}-\delta_{kn}\delta_{lm})\pa_kQ_{lj}\pa_mQ_{nj}\\
&=|\nabla Q|^2-\pa_kQ_{li}\pa_lQ_{ki}.
\end{align*}
In addition, from the properties of the Levi-Civita symbol, we obtain
\begin{align*}
\ve^{ikl}\ve^{jmn}\pa_kQ_{lj}\pa_mQ_{ni}=&|\pa_mQ_{nj}|^2-\pa_nQ_{mj}\pa_mQ_{nj}-\pa_jQ_{nj}\pa_mQ_{nm}\\
&+\pa_nQ_{jj}\pa_mQ_{nm}+\pa_jQ_{mj}\pa_mQ_{nn}-\pa_mQ_{jj}\pa_mQ_{nn}\\
=&|\nabla Q|^2-\pa_nQ_{mj}\pa_mQ_{nj}-\pa_jQ_{nj}\pa_mQ_{nm}\\
=&|\nabla Q|^2-\pa_kQ_{li}\pa_lQ_{ki}-|\nabla\cdot Q|^2,
\end{align*}
where in the second identity we used $\tr(Q)=0$. Then \eqref{nablaqdecom} follows.
\end{proof}
\end{lemma}

Integration by parts gives us that
\begin{align*}
\int_{\Omega}\pa_kQ_{li}\pa_lQ_{ki}-\pa_kQ_{ki}\pa_lQ_{li}dx=\int_{\pa\Omega}Q_{li}(\nu_k\pa_lQ_{ki}-\nu_l\pa_kQ_{ki})dS=0,
\end{align*}
where $\nu$ represents the outward pointing unit normal vector field along $\pa\Omega$. Then  it follows by Lemma \ref{divcurl} that
$$\int_{\Omega}|\nabla \cdot Q|^2dx=\frac{2}{3}\int_{\Omega}|\nabla Q|^2dx-\frac{1}{6}\int_{\Omega}|\mathcal{T}(Q)|^2dx.$$
Thus, the inequality \eqref{mainineq} is equivalent to
\begin{align}\label{omegaestimate}
\int_{\Omega}&\Big(1+\frac{2L}{3}\Big)|\nabla Q|^2+\frac{1}{\ve^2}\mathcal{H}_{Q^K}Q:Qdx-\frac{L}{6}\int_{\Omega}|\mathcal{T}(Q)|^2dx\geq-C\int_{\Omega}|Q|^2dx.
\end{align}
\subsection{Reduction to inequalities for scalar functions on an interval}
 Using the distance function $d^K$  constructed in the inner expansion (cf. \eqref{nabladK-1}), we define
\begin{align}\label{def:Gamma-t-K}
\Gamma_t^K(\delta):=\{x:|d^K(x,t)|<\delta\}.
\end{align}
We also recall
\begin{align*}
\Gamma_t(\delta):=\{x:|d_0(x,t)|<\delta\}.
\end{align*}
\subsubsection{Reduction of the inequality to $\Gamma_t^K(\delta/4)$}
To begin with, due to \eqref{nabladK-1},  we choose $\ve$ sufficiently small  such that
\begin{equation}
\Gamma_t(\delta/8)\subset \Gamma_t^K(\delta/4)\subset \Gamma_t(\delta/2),\ \ \text{for}\ t\in [0,T].\label{subsect:Gamma-delta}
\end{equation}
Let
\begin{align*}
z=\f{d^{K}(x,t)}{\ve}.
\end{align*}
We then  define
\begin{align}
\begin{split}
\tilde{s}(z, x, t):=&\eta\big(\frac{d_0(x, t)} {\delta}\big)s(z)+\Big(1-\eta\big(\frac{d_0(x, t)}{\delta}\big)\Big)
\chi_{  U_+}(x, t) , \label{tildes}\\
\tilde{s}_{1, i}(z, x, t):=&\eta\big(\frac{d_0(x, t)}{\delta}\big)s_{1, i}
(z, x, t), \;  \;  i=0, 3, 4, \\
\tilde{s}_{1, i}(z, x, t):=&\eta\big(\frac{d_0(x, t)}{\delta}\big)
s_{1, i}(z, x, t)+\Big(1-\eta\big(\frac{d_0(x, t)}{\delta}\big)\Big)
\chi_{ U_+}(x, t)q_{1, i}(x, t), \;  \; i=1, 2, \\
\end{split}
\end{align}
where $\eta$ is the cut-off function given in \eqref{eta-function}.
Let $\{E^{i}\}_{i=0}^4$ be the basis defined in \eqref{def:E-al}, we define
\begin{align}
\begin{split}
\tilde{Q}_0(z, x, t):&=\tilde{s}(z, x, t)E^0(x, t), \\
\tilde{Q}_1(z, x, t):&=\sum_{i=0}^{4}\tilde{s}_{1, i}(z, x, t)E^i(x, t), \quad \text{with}\quad z=
\f{d^K(x, t)}{\ve}.
\label{q1near}
\end{split}
\end{align}
Therefore, according to the definition of the approximate solution $Q^K$, (cf. \eqref{def:QK}), we have
\begin{align}
Q^K(x,t)&=\tilde{Q}_0\Big(\frac{d^{K}(x,t)}{\ve},x,t\Big)
+\ve \tilde{Q}_1\Big(\frac{d^{K}(x,t)}{\ve},x,t\Big)+O(\ve^2),
\label{Q^k=Q0+O(ve)}
\end{align}
The symmetric traceless $3 \times 3$ matrix $Q$ is decomposed with respect to the basis $\{E^i\}_{i=0}^{4}$ as
\begin{equation}\label{Qdecom}
 Q(x,t)=\sum_{i=0}^{4}q_{i}(x,t)E^{i}(x,t),
 \end{equation}
 which implies
 \begin{equation}\label{Qest}
 \sum_{i=0}^{4}q_{i}^2\leq C|Q|^2.
 \end{equation}
Inserting \eqref{Q^k=Q0+O(ve)} into the linear operator \eqref{hq0estimta-K}, we obtain
\begin{align}\label{hqaqq}
\frac{1}{\ve^2}\mathcal{H}_{Q^K}Q:Q=&\frac{1}{\ve^2}\Big(Q+9\Big(\f23I(Q:Q^K)
    -QQ^K-Q^KQ\Big)
    \notag\\
    &+3\Big(Q|Q^K|^2+2Q^K(Q:Q^K)\Big)\Big):Q\notag\\
=& \frac{1}{\ve^2}\mathcal{H}_{\tilde{Q}_0}Q:Q
+\frac{1}{\ve}f''_{\tilde{Q}_0}(Q, Q):\tilde{Q}_1+O(1),
\end{align}
where
\begin{align}
    f''_{\tilde{Q}_0}(Q,Q)&=B(Q,Q)+2C(Q,Q,\tilde{Q}_0)\notag\\
    &=9\Big( \f23I|Q|^2-2Q^2\Big)+6\Big(2Q(Q:\tilde{Q}_0)
    +\tilde{Q}_0|Q|^2\Big)\label{f''q0}.
\end{align}

Recall \eqref{the-kap-l} that
\begin{equation*}
\theta(s)=1-6s+6s^2,\ \kappa(s)=1-3s+2s^2\ \text{and}\ \iota (s)=1+6s+2s^2.
\end{equation*}
\begin{lemma}\label{neartheinterface}
  For $t\in[0,T]$, $x\in \Omega \backslash \G_t^K(\delta/4)$ and  $\ve$ sufficiently small, there exists a constant $C>0$ independent of $\ve$, $x$ and $t$ such that
  \begin{align}\label{linear}
  \frac{1}{\ve^2}\mathcal{H}_{\tilde{Q}_0}Q:Q
  +\frac{1}{\ve}f''_{\tilde{Q}_0}&(Q,Q):\tilde{Q}_1\geq -C|Q|^2.
  \end{align}
\begin{proof}
Substituting   \eqref{q1near} and \eqref{Qdecom} into \eqref{hqaqq} and \eqref{f''q0}, we have
\begin{align}\label{hq0qq}
\mathcal{H}_{\tilde{Q}_0}Q:Q&=\frac{2}{3}\theta (\tilde{s}(z)) q_0^2+2\kappa (\tilde{s}(z))( q_1^2+ q_2^2)+2\iota(\tilde{s}(z))(q_3^2+q_4^2)
\end{align}
and
\begin{align*}
f''_{\tilde{Q}_0}&(Q,Q):\tilde{Q}_1\nonumber\\
=&-18(\tilde{Q}_1:Q^2)+6(\tilde{Q}_0:\tilde{Q}_1)|Q|^2+12(\tilde{Q}_0:Q)(
\tilde{Q}_1:Q)\nonumber\\
=&\bigg(4(-1+2\tilde{s})\tilde{s}_{1,0}q_0^2+2(-3+4\tilde{s})
\tilde{s}_{1,0}(q_1^2+q_2^2)+4(3+2\tilde{s})\tilde{s}_{1,0}(q_3^2+q_4^2)
\bigg)\nonumber\\
&+\bigg(4(-3+4\tilde{s})q_0(\tilde{s}_{1,1}q_1+\tilde{s}_{1,2}q_2)
-36\tilde{s}_{1,1}(q_1q_4+q_3q_2)-36\tilde{s}_{1,2}(q_1q_3-q_4q_2)\bigg)\nonumber\\
&+\bigg(8(3+2\tilde{s})q_0(\tilde{s}_{1,3}q_3+\tilde{s}_{1,4}q_4)
-18\tilde{s}_{1,4}(q_1^2-q_4^2)-36\tilde{s}_{1,3}q_1q_2\bigg)\\
\triangleq&I_1+I_2+I_3.\nonumber
\end{align*}

First, we estimate $\mathcal{H}_{\tilde{Q}_0}Q:Q$. Using  \eqref{est:s(z)}, \eqref{tildes} and $\theta(s(\pm\infty))=1$,  we deduce that there exists $M>0$ such that
\begin{align*}
\theta(\tilde{s}(z))=1-6\tilde{s}(z)+6\tilde{s}^2(z)\geq \frac{1}{2}\ \text{for }|z|\geq M.
\end{align*}
 Using  \eqref{est:s(z)}, \eqref{tildes}, $\kappa(s(+\infty))=0$ and  $\kappa(s(-\infty))=1$, we deduce for large $M>0$ that,
\begin{align*}
\kappa(\tilde{s}(z))=1-3\tilde{s}(z)+2\tilde{s}^2(z)=(1-\tilde{s}(z))(1-2\tilde{s}(z))\geq
\left \{
\begin{array}{ll}
\frac{1}{2},\ \ \ \ \ \ \ \ \ \ z\leq -M,\\
-Ce^{-\gamma z},\ z\geq M,
\end{array}
\right.
\end{align*}
where $\gamma>0$ without further explanation.
Thus, thanks to the definition \eqref{def:Gamma-t-K},  for $z=\f{d^{K}(x,t)}{\ve}$ and small $\ve$, one has
\begin{align}\label{theta}
\theta\left(\tilde{s}\Big(\f{d^{K}(x,t)}{\ve}\Big)\right)&\geq \frac{1}{2},\ \ \ \ \ \ \ \ \quad\quad\; x\in\Omega \backslash \G_t^K(\delta/4),\\
\kappa\left(\tilde{s}\Big(\f{d^{K}(x,t)}{\ve}\Big)\right)&\geq -Ce^{-\frac{\gamma}{\ve}},\ \ \ \ \ \ \ x\in\Omega \backslash \G_t^K(\delta/4).\label{kappaest}
\end{align}
Notice that $\iota(s)> 1$ always holds. Therefore, we have by  \eqref{Qest}, \eqref{hq0qq}, \eqref{theta} and \eqref{kappaest} that
\begin{align}\label{hq0}
\mathcal{H}_{\tilde{Q}_0}Q:Q&\geq \frac{2}{3}\theta(\tilde{s})q_0^2-Ce^{-\frac{\gamma}{\ve}}
(q_1^2+q_2^2)+2\iota(\tilde{s})(q_3^2+q_4^2)\notag\\
&\geq \frac{1}{3}q_0^2-C\ve^2|Q|^2+2(q_3^2+q_4^2),\;\;\;\;\quad\quad\quad\ x\in \Omega \backslash \G_t^K(\delta/4),
\end{align}
where we have used $e^{-\frac{\gamma}{\ve}}\leq C\ve^2$.

Next, we estimate $f''_{\tilde{Q}_0}(Q,Q):\tilde{Q}_1$, i.e., $I_1-I_3$. For $\ve$ sufficient small, by definitions \eqref{tildes}, \eqref{q1near} and the construction of $s_{1,i}$ (see Lemmas \ref{lem1}--\ref{lem3}), we have for $ x\in  \Omega \backslash \G_t^K(\delta/4)$ that
\begin{align}
\label{s1034}
\left|\tilde{s}_{1,i}\Big(\f{d^K(x,t)}{\ve},x,t\Big)\right|&\leq Ce^{-\frac{\gamma}{\ve}}\leq C,\;i=0,3,4;\\
\left|\tilde{s}_{1,i}\Big(\f{d^K(x,t)}{\ve},x,t\Big)\right|&\leq C,\;i=1,2;\;\;
\left|\tilde{s}\Big(\f{d^K(x,t)}{\ve},x,t\Big)\right|\leq C.
\label{s112-s}
\end{align}
Then for $x\in  \Omega \backslash \G_t^K(\delta/4),$ we have
\begin{align}
|I_1|&\leq Ce^{-\frac{\gamma }{\ve}}|Q|^2\leq C\ve |Q|^2 .\label{first}
\end{align}
Moreover,  applying \eqref{theta}, \eqref{s112-s} and  Young's inequality, 
one has
\begin{align}
|I_2|&\leq \frac{1}{12\ve}q_0^2+\frac{1}{2\ve}(q_3^2+q_4^2)+C\ve(q_1^2+q_2^2)\nonumber\\
&\leq \frac{1}{12\ve}q_0^2+\frac{1}{2\ve}(q_3^2+q_4^2)+C\ve|Q|^2.\label{second}
\end{align}
Similarly, it follows that
\begin{align}
|I_3|&\leq \frac{1}{12\ve}q_0^2+\frac{1}{2\ve}(q_3^2+q_4^2)+Ce^{-\frac{\gamma}{\ve}}|Q|^2\nonumber\\
&\leq \frac{1}{12\ve}q_0^2
+\frac{1}{2\ve}(q_3^2+q_4^2)+C\ve |Q|^2.\label{third}
\end{align}
Combining \eqref{hq0}, \eqref{first}--\eqref{third}, we obtain \eqref{linear}.
\end{proof}
\end{lemma}
\begin{Remark}
For $(x,t)\in \G_t^K(\delta/4)$, it follows by \eqref{subsect:Gamma-delta} that  $\tilde{s}(z,x,t)=s(z)$, $\tilde{s}_{1,i}(z,x,t)=s_{1,i}(z,x,t)$ for $i=0,\cdots,4$ and  $\tilde{Q}_0=Q_0$, $\tilde{Q}_1=Q_1$, which are all constructed in the inner expansion.
 Thus, thanks to Lemma \ref{neartheinterface}, the inequality  \eqref{omegaestimate} is reduced to
\begin{align}\label{ineqmain}
\int_{\G_t^K(\delta /4)}&(1+\frac{2L}{3})|\nabla Q|^2+\frac{1}{\ve^2}\mathcal{H}_{Q_0}Q:Qdx-\frac{L}{6}\int_{\G_t^K(\delta /4)}|\mathcal{T}(Q)|^2dx\nonumber\\
&+\frac{1}{\ve}\int_{\G_t^K(\delta /4)}f''_{Q_0}(Q,Q):Q_1dx\geq -C\int_{\G_t^K(\delta /4)}|Q|^2dx.
\end{align}
Here,
\begin{align}\label{hq0qqin}
\mathcal{H}_{Q_0}Q:Q&=\frac{2}{3}\theta (s(z)) q_0^2+2\kappa (s(z))( q_1^2+ q_2^2)+2\iota(s(z))(q_3^2+q_4^2),
\end{align}
and
\begin{align}\label{f''in}
f''_{Q_0}&(Q, Q):Q_1\notag\\
=&4(-1+2s)s_{1, 0}q_0^2+2(-3+4s)s_{1, 0}(q_1^2+q_2^2)+4(3+2s)s_{1, 0}(q_3^2+q_4^2)\nonumber\\
&+4(-3+4s)q_0(s_{1, 1}q_1+s_{1, 2}q_2)-36s_{1, 1}(q_1q_4+q_3q_2)-36s_{1, 2}(q_1q_3-q_4q_2)\nonumber\\
&+8(3+2s)q_0(s_{1, 3}q_3+s_{1, 4}q_4)-18s_{1, 4}(q_1^2-q_2^2)-36s_{1, 3}q_1q_2.
\end{align}
\end{Remark}
\subsubsection{Transformation of coordinates and basis decomposition in $\G_t^K(\delta/4)$}
For each $t\in [0,T]$, we give a parameter transformation defined as follows
\begin{align}\label{transformation}
(\sigma,r)\in \Gamma_t^K\times(-\delta /4,\delta /4)\rightarrow x(\sigma,r)\in\Gamma_t^K(\delta/4),
\end{align}
with
\begin{align}\label{parx}
x(\sigma,0)=\sigma\in\Gamma_t^K\text{ and }\pa_r x(\sigma,r)=\frac{\nabla d^K}{|\nabla d^K|^2}\circ(x(\sigma,r)).
\end{align}
Then it follows
$$\frac{d}{dr}(d^K(x(\sigma,r))-r)=\pa_rx\cdot\nabla d^K-1=0.$$
Together with $d^K(x(\sigma,0))=0$, we have
\begin{align}\label{dK=r}
d^K(x(\sigma,r),t)=r.
\end{align}
Applying $d^K=d_0+\sum_{k=1}^K\ve^kd_k$, one has
\begin{align}
\label{d0=r+O(ve)}
d_0=r+O(\ve).
\end{align}
Indeed, $r$ is the distance to $\G_t^K$ and $\sigma$ is the orthogonal projection onto $\G_t^K$.
 Let $$J(\sigma,r)=\det\Big(\frac{\pa x(\sigma,r)}{\pa(\sigma,r)}\Big)$$ be the Jacobian of the transformation. Therefore it follows
\begin{align}\label{est:J}
dx=Jd\sigma dr,\ J|_{\Gamma_t^K}=1,\ \text{ and } J(\sigma,r)=1+O(r).
\end{align}
Moreover, we define the tangent derivative as follows
 \begin{align}\label{def:tan-deri}
\nabla_\Gamma :=\nabla -\nabla d^K\pa_r,
\end{align}
which together with \eqref{parx} gives the orthogonality
\begin{align}\label{bot}
\nabla_\G Q\cdot \nabla d^K&=(\nabla Q-\pa_rQ\nabla d^K)\cdot \nabla d^K\notag\\
&=\nabla Q\cdot \nabla d^K-|\nabla d^K|^2\frac{\nabla d^K}{|\nabla d^K|^2}\cdot \nabla Q\notag\\
&=0.
\end{align}
Recall \eqref{nabladK} that $|\nabla d^K|^2=1+O(\ve^{K+1}).$ Then we have by  definition \eqref{def:tan-deri} and \eqref{bot} that
\begin{align}\label{decom:r-Gamma}
|\nabla Q|^2&=|\nabla d^K|^2|\pa_r Q|^2+|\nabla_\G Q|^2\notag\\
&=\left(1+O(\ve^{K+1})\right)|\pa_r Q|^2+|\nabla_\G Q|^2\notag\\
&=|\pa_r Q|^2+|\nabla_\G Q|^2+O(\ve^{K+1}).
\end{align}
It follows from \eqref{Qdecom} that
\begin{align}
|\pa_r Q|^2&=|\sum_{i=0}^4\pa_rq_i E^i+q_i\pa_rE^i |^2\nonumber\\
&\geq \sum_{i=0}^{4}a_i|\pa_rq_i|^2+2\sum_{\al\neq \beta}\pa_r q_\al q_\beta E^\al:\pa_rE^\beta\nonumber\\
&=\sum_{i=0}^{4}a_i|\pa_rq_i|^2+\sum_{\al\neq \beta}(\pa_r q_\al q_\beta-q_\al\pa_rq_\beta)E^\al:\pa_rE^\beta,\label{parq}
\end{align}
with
\begin{equation*}
a_{i}=\left \{
\begin{array}{ll}
\frac{2}{3},&i=0,\\
2,&i=1,\cdots ,4.
\end{array}
\right.
\end{equation*}

The following lemma shows that $-\frac{L}{6}|\mathcal{T}(Q)|^2$ gives some positive terms on the normal direction of $\Gamma_t^K$ which is the heart of the spectral estimate.
\begin{lemma} \label{tqestimate}
Let $0<\delta_0\leq -\f{L}{12}$, $\pa_r$ and $\nabla_{\Gamma}$ represent the normal and tangent derivative along $\Gamma^K_t$ as in \eqref{transformation} and \eqref{def:tan-deri}. Then there exists  a positive constant $C$ independent of $\delta_0$, $\ve$, $x$, $t$ such that for small $\ve$,
\begin{align}\label{maint}
-\frac{L}{6}\Big|\mathcal{T}(Q)\Big|^2\geq & ~\delta_0
\Big((\pa_rq_1)^2+(\pa_rq_2)^2\Big)\notag\\
&-{C\delta_0}
\sum_{i=0}^4\left((r^2+\ve^2)|\pa_rq_i|^2+|\nabla _\Gamma q_i|^2\right)-C|Q|^2,
\end{align}
\begin{proof}
It follows by the definition \eqref{tq2} and the decomposition \eqref{Qdecom} that
\begin{align*}
|\mathcal{T}(Q)|^2&=|\ve^{ijk}(\pa_j q_\al E_{kl}^\al+q_\al\pa_jE_{kl}^\al)+\ve^{ljk}(\pa_j q_\al E_{ki}^\al+q_\al\pa_jE_{ki}^\al)|^2\\
&\geq \frac{1}{2}|\ve^{ijk}\pa_j q_\al E_{kl}^\al+\ve^{ljk}\pa_j q_\al E_{ki}^\al|^2-|\ve^{ijk}q_\al\pa_jE_{kl}^\al
+\ve^{ljk}q_\al\pa_jE_{ki}^\al|^2\\
&\geq \frac{1}{2}|\ve^{ijk}\pa_j q_\al E_{kl}^\al+\ve^{ljk}\pa_j q_\al E_{ki}^\al|^2-C|Q|^2.
\end{align*}
A tedious calculation gives that
\begin{align}\label{keyineqt}
\frac{1}{2}|\ve^{ijk}\pa_j q_\al E_{kl}^\al+\ve^{ljk}\pa_j q_\al E_{ki}^\al|^2\geq &~ (m\cdot\nabla q_0-n\cdot\nabla q_2+l\cdot\nabla q_3-m\cdot\nabla q_4)^2\nonumber\\
&+(l\cdot\nabla q_0-n\cdot\nabla q_1+m\cdot\nabla q_3 +l\cdot\nabla q_4)^2.
\end{align}
Therefore, it follows for $\delta_0\in(0,-\f{L}{12})$ that
\begin{align}\label{lower:TQ}
-\frac{L}{6}|\mathcal{T}(Q)|^2
\geq&~ 2\delta_0\Big((n\cdot\nabla q_1)^2+(n\cdot\nabla q_2)^2\Big)\notag\\
&-C\delta_0\Big((l\cdot\nabla q_0)^2+(m\cdot\nabla q_0)^2
+(l\cdot\nabla q_3)^2\notag\\
&+(m\cdot\nabla q_3)^2+(m\cdot\nabla q_4)^2+(l\cdot\nabla q_4)^2\Big)-C|Q|^2.
\end{align}
In $\Gamma_t^K(\delta/4),$ we deduce by \eqref{nabladK-1}, $\nabla d_0=n-d_0\bar{h}$(\eqref{d0de}) and \eqref{d0=r+O(ve)} that
\begin{align*}
n\cdot\nabla d^K&=n\cdot\Big(\nabla d_0+\sum_{i=1}^K\ve^i\nabla d_i\Big)\notag\\
&=n\cdot\Big(n-d_0\bar{h}+O(\ve)\Big)\\
&\notag=
n\cdot\Big(n-r\bar{h}+O(\ve)\Big)\\
&=1+O(1)r+O(\ve).
\end{align*}
 Then we have by 
\eqref{def:tan-deri} that for $i=1,2$,
\begin{align}
n\cdot \nabla q_i&=(n\cdot \nabla \notag d^K)\pa_rq_i+n\cdot\nabla_{\Gamma}q_i\\
&=\pa_rq_i+\big(O(\ve)+O(1)r\big)\pa_rq_i+n\cdot \nabla _{\Gamma}q_i,\label{deduce:ncdotnablaq1}
\end{align}
which gives for $i=1,2$ and $\ve$ small enough that
\begin{align}\label{nl}
(n\cdot \nabla q_i)^2&\geq \f34(\pa_rq_i)^2-Cr^2(\pa_r q_i)^2-C|\nabla _\Gamma q_i|^2.
\end{align}
 Similarly as \eqref{deduce:ncdotnablaq1}, for $i=0,3,4$  and $\tau\in\{l,m\}$, we get  by \eqref{d0de}, \eqref{nabladK-1},    \eqref{d0=r+O(ve)}, \eqref{def:tan-deri} and $\tau\cdot n=0$ that
\begin{align*}
\tau\cdot\nabla q_i&=\tau\cdot\Big (n-d_0\bar{h}+\sum_{i=1}^K\ve^i\nabla d_i\Big)\pa_rq_i+\tau\cdot\nabla_{\Gamma}q_i\\
&=\Big(O(\ve)+O(1)d_0\Big)\pa_rq_i+\tau\cdot\nabla_\Gamma q_i,\\
&=\Big(O(\ve)+O(1)r\Big)\pa_rq_i+\tau\cdot\nabla_\Gamma q_i,
\end{align*}
which gives for $i=0,3,4$ and $\tau\in\{l,m\}$ that
\begin{align}\label{nl034}
(\tau\cdot\nabla q_i)^2
&\leq C\ve^2 (\pa_rq_i)^2+Cr^2(\pa_rq_i)^2+C|\nabla _\Gamma q_i|^2.
\end{align}
Summing up \eqref{lower:TQ},\eqref{nl} and \eqref{nl034}, we obtain the desired inequality \eqref{maint}.

To complete the proof, we give details for the estimate \eqref{keyineqt}. It suffices to calculate for $\al=0,\cdots,4$:
\begin{align*}
(\ve^{ijk}\pa_j q_\al E_{kl}^\al)_{\text{sym}}:=\ve^{ijk}\pa_j q_\al E_{kl}^\al+\ve^{ljk}\pa_j q_\al E_{ki}^\al.
 \end{align*}
Recall the basis of symmetric traceless $3\times 3$ matrix, cf. \eqref{def:E-al}:
\begin{align*}
\begin{split}
&E^0(x,t)=\Big(nn-\frac{1}{3}I\Big),\ \ E^1(x,t)=\Big(nl+ln\Big),\ \ E^2(x,t)=\Big(nm+mn\Big),\\
&E^3(x,t)=\Big(ml+lm\Big),\ \ E^4(x,t)=\Big(ll-mm\Big),
\end{split}
\end{align*}
where $\{l,m,n\}$ forms a right-handed coordinate system in $\mathbb{R}^3$:
\begin{align}\label{lmntimes}
l\times m=n,\ m\times n=l,\ n\times l=m,\;\;ll+mm+nn=I.
\end{align}
Then, we have for any vector $u$ that
\begin{align*}
u \times n=l(u \cdot m)-m(u \cdot l),\quad
u \times l=m(u \cdot n)-n(u \cdot m),\quad
u\times m=n(u\cdot l)-l(u\cdot n).
\end{align*}

 Direct calculations give for $\al=0$ that
\begin{align*}
(\ve^{ijk}\pa_j q_0 E_{kl}^0)_{\text{sym}}=&\big(\ve^{ijk}\pa_j q_0 (n_kn_l-\f13 \delta_{kl})\big)_{\text{sym}}\\
=&(\nabla q_0\cdot m)E^1-(\nabla q_0\cdot l)E^2.
\end{align*}
Similarly, we have for $\al=1,2,3,4$ that
\begin{align*}
(\ve^{ijk}\pa_j q_1 E_{kl}^1)_{\text{sym}}=&\big(l(\nabla q_1\times n)+(\nabla q_1\times n)l\big)+\big((\nabla q_1\times l)n+n(\nabla q_1\times l)\big)\\
=&2(\nabla q_1\cdot m)(ll-nn)-(\nabla q_1\cdot l)E^3+(\nabla q_1\cdot n)E^2,\\
(\ve^{ijk}\pa_j q_2 E_{kl}^2)_{\text{sym}}=&\big(m(\nabla q_2\times n)+(\nabla q_2\times n)m\big)+\big((\nabla q_2\times m)n+n(\nabla q_2\times m)\big)\\
=&2(\nabla q_2\cdot l)(nn-mm)+(\nabla q_2\cdot m)E^3-(\nabla q_2\cdot n)E^1,\\
(\ve^{ijk}\pa_j q_3 E_{kl}^3)_{\text{sym}}=&\big(l(\nabla q_3\times m)+(\nabla q_3\times m)l\big)+\big((\nabla q_3\times l)m+m(\nabla q_3\times l)\big)\\
=&(\nabla q_3\cdot l)E^1-(\nabla q_3\cdot n)E^4-(\nabla q_3\cdot m)E^2,\\
(\ve^{ijk}\pa_j q_4 E_{kl}^4)_{\text{sym}}=&\big(l(\nabla q_4\times l)+(\nabla q_4\times l)l\big)-\big((\nabla q_4\times m)m+m(\nabla q_4\times m)\big)\\
=&2(\nabla q_4\cdot n)E^3-(\nabla q_4\cdot m)E^1-(\nabla q_4\cdot l)E^2.
\end{align*}
Thanks to the above calculations, $ll-nn=-\frac{3}{2}E^0+\frac{1}{2}E^4$ and $nn-mm=\frac{3}{2}E^0+\frac{1}{2}E^4$, we infer that
\begin{align*}
\frac{1}{2}|(\ve^{ijk}\pa_j q^\al E_{kl}^\al)_{\text{sym}}|^2=&3(-\nabla q_1\cdot m+\nabla q_2\cdot l)^2\\
&+(\nabla q_0\cdot m-\nabla q_2\cdot n+\nabla q_3\cdot l-\nabla q_4\cdot m)^2\\
&+(\nabla q_0\cdot l-\nabla q_1\cdot n+\nabla q_3\cdot m+\nabla q_4\cdot l)^2\\
&+(-\nabla q_1\cdot l+\nabla q_2\cdot m+2\nabla q_4\cdot n)^2\\
&+(\nabla q_1\cdot m+\nabla q_2\cdot l-\nabla q_3\cdot n)^2\\
\geq &(\nabla q_0\cdot m-\nabla q_2\cdot n+\nabla q_3\cdot l-\nabla q_4\cdot m)^2\nonumber\\
&+(\nabla q_0\cdot l-\nabla q_1\cdot n+\nabla q_3\cdot m+\nabla q_4\cdot l)^2.
\end{align*}
Thus  \eqref{keyineqt} follows. This finishes the proof of the lemma.
\end{proof}
\end{lemma}
Let  $I(\delta)=[-\delta/4,\delta/4].$
Thanks to  \eqref{hqaqq},  \eqref{hq0qqin}, \eqref{decom:r-Gamma}, \eqref{parq} and \eqref{maint}, to prove the inequality \eqref{ineqmain}, it suffices to show  for $\delta_0$ and $\ve$ small enough that
\begin{align*}
\f23&\int_{I(\delta)}\int_{\G_t^K}\Big((1+\frac{2L}{3})|\pa_r q_0|^2+\frac{1}{\ve^2}\theta\big(s(\frac{r}{\ve})\big)q_0^2\Big)Jd\sigma dr\\
&+2\int_{I(\delta)}\int_{\G_t^K}\Big((1+\frac{2L}{3})(|\pa_r q_1|^2+|\pa_rq_2|^2)+\frac{1}{\ve^2}\kappa\big(s(\frac{r}{\ve})\big)
(q_1^2+q_2^2)\Big)Jd\sigma dr\\
&+2\int_{I(\delta)}\int_{\G_t^K}\Big
((1+\frac{2L}{3})(|\pa_rq_3|^2+|\pa_rq_4|^2)+\frac{1}{\ve^2}
\iota\big(s(\frac{r}{\ve})\big)(q_3^2+q_4^2)\Big)Jd\sigma dr\\
&+\frac{1}{\ve}\int_{I(\delta)}\int_{\G_t^K}f''_{Q_0}(Q, Q):Q_1Jd\sigma dr\\
&+(1+\frac{2L}{3})\int_{I(\delta)}\int_{\G_t^K}\sum_{\al\neq \beta}(\pa_r q_\al q_\beta-q_\al \pa_rq_\beta)E^\al:\pa_rE^\beta Jd\sigma dr\\
&+(1+\frac{2L}{3})\int_{I(\delta)}\int_{\G_t^K}|\nabla_\Gamma Q|^2Jd\sigma dr+\delta_0\int_{I(\delta)}\int_{\G_t^K}
\big((\pa_rq_1)^2+(\pa_rq_2)^2\big)Jd\sigma dr\\
&-{C\delta_0}\sum_{i=0}^4
\int_{I(\delta)}\int_{\G_t^K}
\left(|\nabla _\Gamma q_{i}|^2+(r^2+\ve^2)|\pa_rq_{i}|^2\right)Jd\sigma dr\\
&+C\sum_{i=0}^4\int_{I(\delta)}\int_{\G_t^K}q_{i}^2Jd\sigma dr\geq 0.
\end{align*}
To eliminate the tangent derivative terms, we can choose  $\delta_0$ small enough such that
\begin{align}
-{C\delta_0}\sum_{i=0}^4|\nabla _\Gamma q_{i}|^2+(1+\frac{2L}{3})|\nabla_\G Q|^2+ C\sum_{i=0}^4|q_{i}|^2\geq 0.\label{nablagamma}
\end{align}
For convenience, we will assume $\delta/4=1$, i.e., $I=[-1,1]$. And it follows by the definition and \eqref{est:J} that
\begin{align}\label{est:J-I}
|\pa_rJ|+|\pa_r^2J|+|J^{-1}|\leq C\;\;\text{on}\;\;I.
\end{align}
Together with \eqref{nablagamma}, we only need to prove that for each $\sigma \in \G_t^K$ the following inequality holds:
\begin{align}\label{prove}
\f23\int_I&\Big((1+\frac{2L}{3})|\pa_r q_0|^2+\frac{1}{\ve^2}\theta\big(s(\frac{r}{\ve})\big)q_0^2\Big)Jdr\notag\\
&+2\int_I \Big((1+\frac{2L}{3})(|\pa_r q_1|^2+|\pa_rq_2|^2)+\frac{1}{\ve^2}\kappa
\big(s(\frac{r}{\ve})\big)(q_1^2+q_2^2)\Big)J   dr\nonumber\\
&+2\int_I \Big((1+\frac{2L}{3})(|\pa_rq_3|^2+|\pa_rq_4|^2)+\frac{1}{\ve^2}\iota\big(s(\frac{r}{\ve})\big)
(q_3^2+q_4^2)\Big)J dr\notag\\
&+\delta_0\int_I
\Big((\pa_rq_1)^2+(\pa_rq_2)^2\Big)J  dr-C\delta_0\sum_{i=0}^4\int_I (r^2+\ve^2)\big|\pa_rq_{i}\big|^2Jdr\notag\\
&+(1+\frac{2L}{3})\int_I \sum_{\al\neq \beta}(\pa_r q_\al q_\beta-q_\al \pa_rq_\beta)E^\al:\pa_rE^\beta J   dr\notag\\
&+\frac{1}{\ve}\int_I f''_{Q_0}(Q, Q):Q_1J   dr+C\sum_{i=0}^4\int_I q_{i}^2J   dr\geq 0.
\end{align}
\subsubsection{Removing $J$ from the inequality \eqref{prove}}
Let $P=J^{\f12}Q$, i.e., $p_{i}(r)=J^{\frac{1}{2}}(r)q_{i}(r)$, $i=0,\cdots,4$.
And we introduce the following quadratic forms for $p(\cdot)\in H^1(I)$:
\begin{equation}\label{linearo}
\mathcal{G}_{j}(p)=
\begin{cases}
\begin{array}{ll}
\int_I(1+\frac{2L}{3})|\pa_rp|^2(r)
+\frac{1}{\ve^2}\theta\big(s(\frac{r}{\ve})\big)p^2(r)dr,&\text{for}\ j=0,\\
\int_I(1+\frac{2L}{3})|\pa_rp|^2(r)
+\frac{1}{\ve^2}\kappa\big(s(\frac{r}{\ve})\big)p^2(r)dr, &\text{for}\ j=1,\\
\int_I(1+\frac{2L}{3})|\pa_rp|^2(r)+\frac{1}{\ve^2}\iota\big(s(\frac{r}{\ve})\big)
p^2(r)dr,&\text{for}\ j=2.
\end{array}
\end{cases}
\end{equation}
Define
\begin{equation}\label{def:j-al}
j(i)=\left \{
\begin{array}{ll}
0,&i=0,\\
1,&i=1,2,\\
2,&i=3,4.
\end{array}
\right.
\end{equation}
We want to eliminate the $J$ from \eqref{prove}.
It follows that, for $i=0,\cdots,4$ and $0<\nu_0\ll 1$,
\begin{align}\label{lowerbd:pa_rq-J}
\int_I |\pa_rq_{i}|^2Jdr\geq \int_I |\pa_rp_{i}|^2dr-\nu_0\mathcal{G}_{j(i)}(p_{i})-C(\nu_0)
\int_I |p_{i}|^2dr.
\end{align}
Indeed, we deduce by integration by parts and \eqref{est:J-I} that
\begin{align*}
\int_I (\pa_rq_{i})^2Jdr
&=\int_I \Big(J^{-\f12}\pa_rp_{i}+\pa_r(J^{-\f12})p_{i}\Big)^2Jdr\\
&=\int_I (\pa_rp_{i})^2dr+\int_I\big(\pa_r(J^{-\f12})\big)^2p_{i}^2dr
+\f12\int_I\pa_r\big(J^{-1}\big)\pa_r\big(p_{i}^2\big)dr\\
&=\int_I (\pa_rp_{i})^2dr+\int_I\Big(\big(\pa_r(J^{-\f12})\big)^2-\f12\pa_r^2\big(J^{-1}\big)
\Big)p_{i}^2dr
+\f12\pa_r\big(J^{-1}\big)p_{i}^2\Big|_{-1}^1\\
&\geq \int_I (\pa_rp_{i})^2dr-C(\nu_0)\int_I p_{i}^2dr-\nu_0\mathcal{G}_{j(i)}(p_{i}),
\end{align*}
where we have used Corollary \ref{cor:p(pm1)} to bound $p_{i}^2(\pm 1)$.

 In a similar manner, we deduce  that, for $i=0,\cdots,4$ and $0<\nu_0\ll 1$,
\begin{align}\label{lowerbd:pa_rq-J-ve2r2}
-\int_I (r^2+\ve^2)\big(\pa_rq_{i}\big)^2Jdr\geq -\int_I (r^2+\ve^2)\big(\pa_rp_{i}\big)^2
dr-\nu_0\mathcal{G}_{j(i)}(p_{i})-C(\nu_0)
\int_I |p_{i}|^2dr.
\end{align}
Thanks to \begin{align*}
(\partial_rq_\al q_\beta-q_\al\partial_rq_\beta)J=&\partial_r(J^{-\frac{1}{2}}p_\al)p_\beta J^{\frac{1}{2}}-\partial_r(J^{-\frac{1}{2}}p_\beta)p_\al J^{\frac{1}{2}}
=\partial_rp_\al p_\beta-p_\al\partial_rp_\beta,
\end{align*}
it follows that
 \begin{align}\label{removeJ2}
\int_I \sum_{\al\neq \beta}(\pa_r q_\al q_\beta-q_\al \pa_rq_\beta)E^\al:\pa_rE^\beta Jdr=\int_{I}\sum_{\al\neq \beta}(\pa_r p_\al p_\beta-p_\al\pa_rp_\beta)E^\al:\pa_rE^\beta dr.
\end{align}
Since $P=J^{\f12}Q$, one has by definition \eqref{f''in} that
\begin{align}\label{def:fPP}
\int_I&f''_{Q_0}(Q,Q):Q_1Jdr=
\int_If''_{Q_0}(P,P):Q_1dr.
\end{align}
Choosing $\nu_0$ sufficiently small such that
\begin{align*}
&\f23-\f23(1+\f{2L}{3})\nu_0-\delta_0\nu_0+C\delta_0\nu_0\geq \f13,\\
&2-2(1+\f{2L}{3})\nu_0-\delta_0\nu_0+C\delta_0\nu_0\geq 1,\\
&2-2(1+\f{2L}{3})\nu_0-\delta_0\nu_0+C\delta_0\nu_0\geq 1.
\end{align*}
Combining \eqref{linearo}--\eqref{def:fPP}, the inequality
\eqref{prove} can be reduced to
\begin{align}\label{prove2}
\f13\mathcal{G}_0(p_0)&+
\sum_{i=1, 2}\mathcal{G}_1(p_{i})+
\sum_{i=3, 4}\mathcal{G}_2(p_{i})+\delta_0\int_{I}
(\pa_rp_1)^2+(\pa_rp_2)^2dr\notag\\
&+C\sum_{i=0}^4\int_{I} p_{i}^2dr+(1+\frac{2L}{3})\underbrace{\sum_{\al\neq \beta}\int_{I}(\pa_r p_\al p_\beta-p_\al \pa_rp_\beta)E^\al:\pa_rE^\beta dr}_{\text{cross terms}}\notag\\
&\underbrace{-C\delta_0\sum_{i=0}^4\int_I (r^2+\ve^2)\big(\pa_rp_{i}\big)^2
dr}_{\text{convolution terms}}+\underbrace{\frac{1}{\ve}\int_{I}f''_{Q_0}(P, P):Q_1dr}_{
\text{correction terms}}\geq 0,
\end{align}
for some $0<\delta_0\ll 1$ and $\ve$ small enough. It suffices to prove the following inequalities.
\begin{itemize}[label=$\bullet$]
\item Convolution terms: For some $0<\delta_0\ll 1$ and  sufficiently small $\ve$, it holds that
\begin{align}\label{prove-Linear}
&\quad C\delta_0\sum_{i=0}^4\int_I (r^2+\ve^2)\big(\pa_rp_{i}\big)^2
dr\leq \f19\sum_{i=0}^4\mathcal{G}_{j(i)}(p_{i})+C\sum_{i=0}^4\int_{I} p_{i}^2dr.
\end{align}

\item Cross terms:
For some $0<\delta_0\ll 1$ and  sufficiently small $\ve$, it holds that

\begin{align}\label{prove-Cross}
(1+\frac{2L}{3})&\sum_{\substack{\al,\beta\in\{0,\cdots,4\}\\ \al\neq \beta}}\int_{I}(p_{\al}\pa_r p_{\beta} -p_{\beta} \pa_rp_{\al}) E^{\al}:\pa_rE^{\beta} dr\notag\\
&\leq \f19\sum_{i=0}^4\mathcal{G}_{j(i)}(p_{i})+\f{\delta_0}{2}
\int_{I}
(\pa_rp_1)^2+(\pa_rp_2)^2dr
+\sum_{i=0}^4\int_{I} p_{i}^2dr.
\end{align}

\item Correction terms:
For some $0<\delta_0\ll 1$ and  sufficiently small $\ve$, it holds that
\begin{align}\label{prove-Correction}
\frac{1}{\ve}&\Big|\int_{I}f''_{Q_0}(P, P):Q_1dr\Big|\notag\\
&\leq \f{\delta_0}{2}
\int_{I}
(\pa_rp_1)^2+(\pa_rp_2)^2dr+\f19\sum_{i=0}^4\mathcal{G}_{j(i)}(p_{i})
+\sum_{i=0}^4\int_{I} p_{i}^2dr.
\end{align}
\end{itemize}
\subsection{Estimate for quadratic forms $\mathcal{G}_j$ ($j=0,1,2$) }
About the estimate for $\mathcal{G}_{j}$, we can refer to \cite{C1994, FWZZ2018, FW}.
Obviously, thanks to $\iota (s)>1$, we have
\begin{align}\label{lowerG2}
\mathcal{G}_{2}(p)
&=\int_I(1+\frac{2L}{3})|\pa_rp|^2+\frac{1}{\ve^2}\iota\big(s(\frac{r}{\ve})\big)
p^2(r)dr\notag\\
&\geq \int_I(1+\frac{2L}{3})|\pa_rp|^2+\frac{1}{\ve^2}p^2(r)dr>0. 
\end{align}
It suffices to give the estimates for  $\mathcal{G}_j(p)$ for $j=0,1$.
In the sequel, we  always assume $p(r)\in H^1(I)$.
Recall \eqref{def:s(z)}:
\begin{align*}
s(z)=\frac{1}{1+e^{-\gamma z}}\;\;\text{with}\;\;\gamma=\frac{1}{\sqrt{1+\frac{2L}{3}}}>0,
\end{align*}
with properties listed in Lemma \ref{lem:pro-s}, i.e., \eqref{s'}--\eqref{kappaf}. Set
\begin{equation}\label{s01ve}
\xi_{1,\ve}(r):=s\big(\frac{r}{\ve}\big),\qquad \ \xi_{0,\ve}(r):=s'\big(\frac{r}{\ve}\big)=\ve\pa_r\xi_{1,\ve}.
\end{equation}
\subsubsection{Coercive estimates}
\begin{lemma}\label{le:0}
Let $p_0=\xi_{0,\ve}\bar{p}_0$. Then for any $\nu_0\in(0,\f12)$, there exists a positive constant $C(\nu_0)$ such that
\begin{align}\label{es:0}
\mathcal{G}_0(p_0)\geq &(1+\frac{2L}{3})
\big(\frac{1}{2}-\nu_0-O(e^{-\frac{\gamma}{\ve}})\big)\int_{I}(\xi_{0, \ve}\partial_r \bar{p}_0)^2dr\notag\\
&-C(\nu_0)\ve^{-2}e^{-\frac{2\gamma}{\ve}}\int_{I}\xi_{0, \ve}^2\bar{p}_0^2dr.
\end{align}
\proof
By the definition \eqref{s01ve} and \eqref{kappaf}, one has $$(1+\frac{2L}{3})\ve^2\partial_r^2\xi_{0,\ve}=\xi_{0,\ve}\theta(\xi_{1,\ve}).$$
Then we integrate by parts and arrive at
\begin{align}\label{est:G0-inte}
\mathcal{G}_0(p_0)&=\int_{I}(1+\frac{2L}{3})\big((\partial_r\xi_{0,\ve}\bar{p}_0+\xi_{0,\ve}\partial_r\bar{p}_0)^2\big)+\ve^{-2}\theta(\xi_{1,\ve})s^2_{0,\ve}\bar{p}_0^2dr\nonumber\\
&=(1+\frac{2L}{3})\xi_{0,\ve}\partial_r\xi_{0,\ve}\bar{p}_0^2\Big|_{-1}^1+(1+\frac{2L}{3})\int_{I}s^2_{0,\ve}(\partial_r\bar{p}_0)^2dr.
\end{align}
It suffices to estimate the boundary terms. Notice
\begin{align}
\begin{split}
\xi_{0,\ve}(r)&=
\f{\gamma}{4}\mathrm{sech}^2 (\f{\gamma r}{2\ve}),
\;\; \pa_r\xi_{0,\ve}(r)=\f{-\gamma^2}{4\ve}\mathrm{sech}^2(\f{\gamma r}{2\ve}) \tanh (\f{\gamma r}{2\ve}),\\
|\tanh z|&\leq 1,\; \mathrm{sech} z\leq 2e^{-|z|},\; \cosh z=\f12e^{|z|}\big(1+O(e^{-2|z|})\big).\label{s0ve:basic}
\end{split}
\end{align}
Thanks to \eqref{s0ve:basic}, we can directly get that
\begin{align}\label{est:sve01}
\Big|\big(\xi_{0, \ve}\pa_r\xi_{0, \ve}\big)(\pm 1)\Big|=&\frac{\gamma^3}{16\ve}\mathrm{sech}^4 (\f{\pm\gamma }{2\ve})|\tanh  (\f{\pm\gamma }{2\ve})|\leq \frac{\gamma^3 e^{-\frac{2\gamma}{\ve}}}{\ve}, \notag\\
\int_{0}^{1}\xi_{0, \ve}^{-2}dr=&\frac{16}{\gamma^2}\int_{0}^{1}
\mathrm{cosh}^4 (\f{\gamma r}{2\ve})dr\notag\\
=& \frac{1}{\gamma^2}\int_{0}^{1}e^{\f{2\gamma r}{\ve}}\big(1+O(e^{-\frac{\gamma r}{\ve}})\big)dr\notag\\
=&\frac{\ve e^{\frac{2\gamma}{\ve}}}{2\gamma^3} \Big(1+O(e^{-\frac{\gamma}{\ve}})\Big),
\end{align}
and
\begin{align}\label{s0ve:basic-2}
\xi_{0,\ve}(r)\geq \f{\gamma}{4}\mathrm{sech}^2 (\f{\gamma }{2})\gtrsim 1,\;|r|\leq \ve.
\end{align}
Notice there exists $a\in [0,\ve]$ such that
\begin{align*}
\bar{p}_0(a)^2&\leq \ve^{-1}
\int_{0}^{\ve}\bar{p}_0^2dr.
\end{align*}
Together with  Young's inequality, \eqref{est:sve01} and \eqref{s0ve:basic-2} imply
\begin{align}\label{est:p1-1}
\bar{p}_0(1)^2&=
\Big(\int_{a}^{1}\partial_r\bar{p}_0dr+\bar{p}_0(a)\Big)^2\notag\\
&\leq (1+\f{\nu_0}{2})
\Big(\int_{0}^{1}\partial_r\bar{p}_0dr\Big)^2+(1+\f{2}{\nu_0})\bar{p}_0(a)^2\notag\\
&\leq (1+\f{\nu_0}{2})\big(\int_{0}^{1}\xi_{0,\ve}^{-2}dr\big)
\Big(\int_{0}^{1}\xi_{0,\ve}^2(\partial_r\bar{p}_0)^2dr\Big)
+(1+\f{2}{\nu_0})C\ve^{-1}\int_{0}^{\ve}\xi_{0,\ve}^2\bar{p}_0^2dr\notag\\
&\leq 
(1+\f{\nu_0}{2})\f{\ve e^{\f{2\gamma}{\ve}}}{2\gamma^3}\big(1+O(e^{-\frac{\gamma}{\ve}})\big)
\int_{-1}^1\xi_{0,\ve}^2(\partial_r\bar{p}_0)^2dr
+\frac{C(\nu_0)}{\ve}\int_{-1}^{1}\xi_{0,\ve}^2\bar{p}_0^2dr.
\end{align}
Combining \eqref{est:sve01} and \eqref{est:p1-1},  we obtain by $\nu_0\in(0,\f12)$ that
\begin{align}\label{es:0f}
|\xi_{0, \ve}\partial_r\xi_{0, \ve}\bar{p}_0^2(1)|
\leq& \Big(\frac{1}{2}+\nu_0+O(e^{-\frac{\gamma}{\ve}})\Big)\int_{I} \xi_{0, \ve}^2(\pa_r\bar{p}_0)^2 dr\notag\\
&+C(\nu_0)\ve^{-2}e^{-\frac{2\gamma}{\ve}}\int_{I}\xi_{0, \ve}^2\bar{p}_0^2dr.
\end{align}
The boundary term $|\xi_{0,\ve}\partial_rs _{0,\ve}\bar{p}_0^2(-1)|$ can be estimated similarly. Thus, \eqref{es:0} follows by applying  these two estimates to \eqref{est:G0-inte}.
\qed
\end{lemma}
\begin{lemma}\label{le:1}
Let $p=\xi_{1,\ve}\bar{p}$. For any $\nu_0\in(0,\f12),$ there exists $C(\nu_0)>0$  such that
\begin{align*}
\mathcal{G}_1(p)\geq&\big(1+\frac{2L}{3})\big(\frac{1}{2}-\nu_0-O(
e^{-\frac{\gamma}{\ve}})\big)\int_{I}\xi_{1, \ve}^2(\pa_r\bar{p})^2dr\notag\\
&-C(\nu_0)\ve^{-1}e^{-\frac{2\gamma}{\ve}}\int_{I}\xi_{1, \ve}^2\bar{p}^2dr.
\end{align*}
\begin{proof}
By \eqref{kappaf}, a direct computation shows that
\begin{align}\label{trans}
\int_{I}&\big(1+\frac{2L}{3})|\pa_rp|^2+\frac{1}{\ve^2}\kappa(s)p^2dr\notag\\
&=\int_{I}\big(1+\frac{2L}{3})(\pa_r \xi_{1, \ve}\bar{p}+s_{1\ve}\pa_r\bar{p})^2+\frac{1}{\ve^2}\kappa(s)\bar{p}^2\xi_{1, \ve}^2dr\nonumber\\
&=\big(1+\frac{2L}{3})\xi_{1, \ve}
\pa_r\xi_{1, \ve}\bar{p}^2|_{-1}^1+\big(1+\frac{2L}{3})\int_{I}
\xi_{1, \ve}^2(\pa_r\bar{p})^2dr\notag\\
&\geq -\big(1+\frac{2L}{3})\xi_{1, \ve}\pa_r\xi_{1, \ve}\bar{p}^2(-1)+\big(1+\frac{2L}{3})\int_{I}\xi_{1, \ve}^2(\pa_r\bar{p})^2dr,
\end{align}
where the last line is due to $\xi_{1,\ve}(r)\pa_r\xi_{1,\ve}(r)\geq 0$.
Notice
\begin{align*}
\xi_{1,\ve}(r)&=\f{1}{1+e^{-\f{\gamma r}{\ve}}}\in(0,1),\; \pa_r\xi_{1,\ve}(r)=\f{\gamma}{4\ve}\mathrm{sech}^2(\f{\gamma r}{2\ve}) ,\;
\end{align*}
which  together with $\mathrm{sech} z\leq 2e^{-|z|}$ gives 
\begin{align}\label{est:sve12}
&(\xi_{1,\ve}\pa_r\xi_{1,\ve})(-1)\leq \frac{1}{\ve}\gamma e^{-\frac{2\gamma}{\ve}},\\
\int_{-1}^{1}\xi_{1,\ve}^{-2}dr&=\int_{-1}^{0}(1+e^{-\frac{\gamma r}{\ve}})^2dr+\int_{0}^{1}(1+e^{-\frac{\gamma r}{\ve}})^2dr\notag\\
&=\frac{\ve e^{\frac{2\gamma}{\ve}}}{2\gamma}\big(1+O(e^{-\frac{\gamma}{\ve}})\big)
+O(1)=\frac{\ve e^{\frac{2\gamma}{\ve}}}{2\gamma}\big(1+O(e^{-\frac{\gamma}{\ve}})\big).
\label{s1ve:basic-2}
\end{align}
We get by $\xi_{1,\ve}(r)\geq \f12,\;r\in[0,1]$ that there exists $a\in [0,1]$ such that
\begin{align}\label{est:p(a)-1}
\bar{p}(a)^2&\leq
\int_{0}^{1}\bar{p}^2dr\leq 4\int^{1}_0\xi_{1,\ve}^2\bar{p}^2dr.
\end{align}
This together with  Young's inequality, \eqref{est:sve12}-\eqref{est:p(a)-1} implies
\begin{align*}
(\xi_{1, \ve}\pa_r\xi_{1, \ve}\bar{p}^2)(-1)\leq &\xi_{1, \ve}\pa_r\xi_{1, \ve}(-1)\Big
(\int_{-1}^{a}\pa_r\bar{p}dr-\bar{p}(a)\Big)^2\notag\\
\leq &\frac{1}{\ve}\gamma e^{-\frac{2\gamma}{\ve}}
\Big((1+\nu_1)\int_{-1}^{1}\xi_{1, \ve}^{-2}dr\int_{-1}^{1}\xi_{1, \ve}^2(\pa_r\bar{p})^2dr
\notag\\
&+4(1+\nu_1^{-1})\int_{0}^1\xi_{1, \ve}^2\bar{p}^2dr\Big)\notag\\
\leq &\Big(\frac{1}{2}+\nu_0+O(e^{-\frac{\gamma}{\ve}})\Big)
\int_{I}\xi_{1, \ve}^2(\pa_r\bar{p})^2dr+C(\nu_0)\ve^{-1}
e^{-\frac{2\gamma}{\ve}}\int_{I}\xi_{1, \ve}^2\bar{p}^2dr.
\end{align*}
Applying the estimate to \eqref{trans}, we obtain the lemma.
\end{proof}
\end{lemma}
\begin{lemma}[Estimate of the second eigenvalue of $\mathcal{G}_0$]\label{secondeigenvalue}
Let $p_0=\mu \xi_{0,\ve}+\hat{p}_0$ with $\int_{I}\xi_{0,\ve}\hat{p}_0dr=0$ and $\mu\in\BR$. Then there exists $c_0>0$ such that for $\ve$ small,
\begin{align}
&\mathcal{G}_0(p_0)+o(\ve^2)\int_Ip_0^2dr\geq \frac{c_0}{\ve^2}\int_I\hat{p}_0^2dr,\label{firstfirst}\\
&\mathcal{G}_0(p_0)+o(\ve^2)\int_{I}p_0^2dr\geq c_0\int_{I}(\partial_r\hat{p}_0)^2dr.\label{firstdsecond}
\end{align}
\end{lemma}
\begin{proof}
Let $\bar{p}_0=p_0/\xi_{0,\ve}-\mu=\hat{p}_0/\xi_{0,\ve}$. Then it follows from Lemma \ref{le:0} that
\begin{align}\label{6.6main}
\mathcal{G}_0(p_0)+o(\ve^2)
\int_I{p}_0^2\ud r\ge \frac14\Big(1+\f23L\Big)\int_I\xi_{0,\ve}^2
\big(\partial_r\bar{p}_0\big)^2\ud r.
\end{align}
Then \eqref{firstfirst} can be reduced to
\begin{align*}
  \int_I\xi_{0,\ve}^2\big(\partial_r\bar{p}_0\big)^2\ud r
  \ge\frac{c_0}{\ve^2}\int_I\hat{p}_0^2dr.
\end{align*}
whose proof is analogous to that of \cite[Lemma 7.5 and Corollary 7.6]{FW}.
The remaining argument is identical, and it suffices to verify the following estimate for $\xi_{0,\ve}(r)=s'(\f{r}{\ve})$:
\begin{align*}
I_0(r)
&:=
\int_r^1 \xi_{0,\ve}^2(\tau)
\left(\int_0^\tau \frac{1}{\xi_{0,\ve}(y)}\,dy\right)d\tau \\
&=
\ve^2
\int_{r/\ve}^{1/\ve}
(s')^2(z)
\left(\int_0^z \frac{1}{s'(w)}\,dw\right)dz \\
&\lesssim
\ve^2\int_{r/\ve}^{1/\ve}s'(z)\,dz \le
\ve^2\int_{r/\ve}^{+\infty}s'(z)\,dz
\lesssim
\ve^2s'(r/\ve)
=
\ve^2\xi_{0,\ve}(r).
\end{align*}
\eqref{firstdsecond} is a direct consequence of \eqref{6.6main} and \eqref{firstfirst} using \begin{align*}|\pa_r\hat{p}_0|\lesssim |\xi_{0,\ve}\pa_r\bar{p}_0|+\ve^{-1}|s''(r/\ve)\bar{p}_0|
\lesssim |\xi_{0,\ve}\pa_r\bar{p}_0|
+\ve^{-1}|\xi_{0,\ve}\bar{p}_0|
\lesssim |\xi_{0,\ve}\pa_r\bar{p}_0|
+\ve^{-1}|\hat{p}_0|.
\end{align*}
\end{proof}
To estimate the second eigenvalue of $\mathcal{G}_1$, we introduce an exponentially weighted $L^2$ norm.
Let $\omega(r)$ be a non-negative and bounded function which decays exponentially to zero at infinity. For convenience, we write $\omega_\ve(r)=\omega (\frac{r}{\ve})$. We also require $\omega>0$ a.e., which ensures that $\|p\|_{\omega}=\big(\int_{I}p^2\omega_{\ve}\;dr\big)^{\f12}$ is a norm.
\begin{lemma}[Estimate of the second eigenvalue of $\mathcal{G}_1$]\label{g1second}
Let $p_1=\nu \xi_{1,\ve}+\hat{p}_1$ with $\int_{I}\xi_{1,\ve}\hat{p}_1\omega_\ve dr=0$, where  $\nu \in\BR$. Then there exists $c_0>0$ such that for   sufficiently small $\ve$,
\begin{align}
    &\mathcal{G}_1(p_1)+o(\ve^2)\int_Ip_1^2dr\geq \frac{c_0}{\ve^2}\int_I \hat{p}_1^2 \omega_\ve dr,\\
    &\mathcal{G}_1(p_1)+o(\ve^2)\int_{I}p_1^2dr\geq c_0\int_{I}(\partial_r\hat{p}_1)^2\omega_{\ve} dr.\label{secondsecond}
\end{align}
\end{lemma}
The following lemma is the key to prove Lemma \ref{g1second}.
\begin{lemma}
For
\begin{align*}
&K_+(r):=
\int_r^1
\xi_{1,\ve}^2\omega_\ve(\tau)
\left(
\int_0^\tau
\xi_{1,\ve}^{-1}(y)\,dy
\right)d\tau,\text{ for }r\in[0,1],\\
&K_-(r):=
\int_{-1}^r
\xi_{1,\ve}^2\omega_\ve(\tau)
\left(
\int_\tau^0
\xi_{1,\ve}^{-1}(y)\,dy
\right)d\tau,\text{ for }r\in[-1,0],
\end{align*}
one has
\begin{align}
&K_+(r)
\lesssim M_{\omega} \varepsilon^2\xi_{1,\ve}(r),
\qquad r\in[0,1],\label{eq:K-plus}\\
&K_-(r)\lesssim M_{\omega}\varepsilon^2\xi_{1,\ve}(r),
\qquad r\in[-1,0],\label{eq:K-minus}
\end{align}
where
\begin{align*}
M_\omega:=
\|\omega\|_{L^\infty(\mathbb R)}
+\int_0^{+\infty} z\omega(z)\,dz<+\infty.
\label{eq:M-omega}
\end{align*}
\end{lemma}
\begin{proof}
Since $\f12\le\xi_{1,\ve}\le 1$ on $[0,1]$, it follows that
\begin{align*}
K_+(r)
&\leq
2\int_r^1
\tau\omega\left(\frac{\tau}{\varepsilon}\right)d\tau
=
2\varepsilon^2
\int_{r/\varepsilon}^{1/\varepsilon}
z\omega(z)\,dz
\leq
2M_\omega\varepsilon^2.
\end{align*}
Using again $\xi_{1,\ve}(r)\geq\frac12$, we obtain \eqref{eq:K-plus}.

For $\tau\leq0$, using $\xi_{1,\ve}^{-1}(y)
=
1+e^{-\gamma y/\varepsilon},$
we find
\begin{equation}
\int_\tau^0\xi_{1,\ve}^{-1}(y)\,dy
\leq
C\varepsilon\xi_{1,\ve}^{-1}(\tau).
\label{eq:xi-inverse-left}
\end{equation}
Moreover,
$
\partial_r\xi_{1,\ve}
=
\frac{\gamma}{\varepsilon}
\xi_{1,\ve}(1-\xi_{1,\ve})
\geq
\frac{\gamma}{2\varepsilon}\xi_{1,\ve},
$ for $r \leq 0$,
and hence
\begin{equation}
\int_{-1}^r\xi_{1,\ve}(\tau)\,d\tau
\leq
\frac{2\varepsilon}{\gamma}\xi_{1,\ve}(r).
\label{eq:xi-left-integral}
\end{equation}
It follows from \eqref{eq:xi-inverse-left},
\eqref{eq:xi-left-integral}
that
\begin{align}\nonumber
K_-(r)
\leq
C\varepsilon\|\omega\|_{L^\infty}
\int_{-1}^r\xi_{1,\ve}(\tau)\,d\tau
\lesssim M_{\omega}
\varepsilon^2\xi_{1,\ve}(r),
\qquad -1\leq r\leq0,
\end{align}
which concludes the proof.
\end{proof}
\begin{proof}[Proof of Lemma \ref{g1second}]
Let $\bar{p}_1=p_1/\xi_{1,\ve}-\nu=\hat{p}_1/\xi_{1,\ve}$. Then it follows from Lemma \ref{le:1} that
\begin{align*}
\mathcal{G}_1(p_1)+o(\ve^2)
\int_I{p}_1^2\ud r\ge \frac14\Big(1+\f23L\Big)\int_I\xi_{1,\ve}^2
\big(\partial_r\bar{p}_1\big)^2\ud r.
\end{align*}
Now we prove that for $\int_I\xi_{1,\ve}^2
\bar{p}_1\omega_{\ve}\ud r=0$ and some $c_0(\omega)>0$,
\begin{align}\label{6.7main}
  \int_I\xi_{1,\ve}^2\big(\partial_r\bar{p}_1\big)^2\ud r
  \ge\frac{c_0(\omega)}{\ve^2} \int_I\xi_{1,\ve}^2 \bar{p}_1^2 \omega_\ve\ud r=\frac{c_0(\omega)}{\ve^2}\int_I\omega_\ve(r)\hat{p}_1^2dr.
\end{align}
We  follow the decomposition used in the proof of \cite[Lemma 7.5]{FW}. Set
\begin{align*}
&A_+:=\int_0^1\xi_{1,\ve}^2\omega_\ve \bar{p}_1^2\,dr,
\qquad
B_+:=\int_0^1\xi_{1,\ve}^2
|\partial_r \bar{p}_1|^2\,dr,
\qquad
D_+:=\int_0^1\xi_{1,\ve}^2\omega_\ve \bar{p}_1\,dr,\\
&A_-:=\int_{-1}^0\xi_{1,\ve}^2\omega_\ve \bar{p}_1^2\,dr,
\qquad
B_-:=\int_{-1}^0\xi_{1,\ve}^2
|\partial_r \bar{p}_1|^2\,dr,
\qquad
D_-:=\int_{-1}^0\xi_{1,\ve}^2\omega_\ve \bar{p}_1\,dr.
\end{align*}
For $\tau\in[0,1]$, write
\[
\bar{p}_1(\tau)=\bar{p}_1(0)+\int_0^\tau\partial_r \bar{p}_1(r)\,dr.
\]
It follows that
\begin{align}
A_+
&=
\bar{p}_1(0)D_+
+
\int_0^1
\xi_{1,\ve}^2\omega_\ve(\tau)\bar{p}_1(\tau)
\left(
\int_0^\tau\partial_r \bar{p}_1(r)\,dr
\right)d\tau
\nonumber\\
&\leq
\bar{p}_1(0)D_++A_+^{1/2}\Big[\int_0^1
\xi_{1,\ve}^2\omega_\ve(\tau)
\left(
\int_0^\tau\partial_r \bar{p}_1(r)\,dr
\right)^2d\tau\Big]^{1/2}.\label{eq:A-plus-first}
\end{align}
By the Cauchy--Schwarz inequality,
\[
\left(
\int_0^\tau\partial_r \bar{p}_1(r)\,dr
\right)^2
\leq
\left(
\int_0^\tau
\xi_{1,\ve}(r)|\partial_r \bar{p}_1(r)|^2\,dr
\right)
\left(
\int_0^\tau\xi_{1,\ve}^{-1}(r)\,dr
\right).
\]
Therefore, by Fubini's theorem and \eqref{eq:K-plus},
\begin{align}
\int_0^1
\xi_{1,\ve}^2\omega_\ve(\tau)
\left(
\int_0^\tau\partial_r \bar{p}_1(r)\,dr
\right)^2d\tau
&\leq
\int_0^1
\xi_{1,\ve}(r)|\partial_r \bar{p}_1(r)|^2K_+(r)\,dr
\nonumber\\
&\leq
C(\omega)\varepsilon^2
\int_0^1
\xi_{1,\ve}^2(r)|\partial_r \bar{p}_1(r)|^2\,dr
\nonumber\\
&=
C(\omega)\varepsilon^2B_+.
\label{eq:J-plus}
\end{align}
Combining \eqref{eq:A-plus-first} and \eqref{eq:J-plus},
we obtain
\begin{equation}
A_+
\leq
\bar{p}_1(0)D_+
+
C(\omega)\varepsilon A_+^{1/2}B_+^{1/2}.
\label{eq:A-plus}
\end{equation}
Similarly, for $\tau\in[-1,0]$, one has
\begin{equation}
A_-
\leq
\bar{p}_1(0)D_-
+
C(\omega)\varepsilon A_-^{1/2}B_-^{1/2}.
\label{eq:A-minus}
\end{equation}
Adding \eqref{eq:A-plus} and \eqref{eq:A-minus}, and  noticing $\int_I\xi_{1,\ve}^2
\bar{p}_1\omega_{\ve}\ud r=0$ gives
 $D_++D_-=0$, we have
\begin{align}
A_++A_-
&\leq
C(\omega)\varepsilon
\left(
A_+^{1/2}B_+^{1/2}
+
A_-^{1/2}B_-^{1/2}
\right)
\nonumber\\
&\leq
C(\omega)\varepsilon
(A_++A_-)^{1/2}(B_++B_-)^{1/2}.
\label{eq:A-B-total}
\end{align}
If $A_++A_-=0$, then \eqref{6.7main}
is immediate. Otherwise, dividing
\eqref{eq:A-B-total} by $(A_++A_-)^{1/2}$ gives
\[
A_++A_-
\leq
C(\omega)\varepsilon^2(B_++B_-).
\]
Consequently, for some $c_0(\omega)>0$ independent of
$\varepsilon$,
\begin{align}
\int_I\xi_{1,\ve}^2
|\partial_r \bar{p}_1|^2\,dr
&\geq
\frac{c_0}{\varepsilon^2}
\int_I\xi_{1,\ve}^2\omega_\ve \bar{p}_1^2\,dr=
\frac{c_0}{\varepsilon^2}
\int_I
\omega_\varepsilon(r)\widehat p_1^2(r)\,dr.
\label{eq:weighted-poincare-final}
\end{align}
This proves \eqref{6.7main}. At last, as $|\partial_r\xi_{1,\ve}|\le \frac{C}{\ve}\xi_{1,\ve}$, and $|\partial_r\hat{p}_1|\le|\partial_r\xi_{1,\ve} \bar{p}_1|+|\xi_{1,\ve}\partial_r\bar{p}_1| $, we obtain \eqref{secondsecond}.
The proof is completed.
\end{proof}
\subsubsection{Endpoints $L^\infty$ estimates}
Again, we assume that $\ve$ is sufficiently small.
\begin{lemma}\label{leminfty}
Let $j\in\{1,2\}$. Then for any $\nu_0\in(0,\f12)$, there exists a positive constant $C(\nu_0)$ such that
\begin{equation}\label{infty}
\|p\|_{L^\infty([-1,1])}^2\leq \nu_0\mathcal{G}_j(p)+C(\nu_0)\int_{I}p^2dr.
\end{equation}
\begin{proof}
For $j=1,$ let $p=\xi_{1,\ve}\bar{p}$. It follows from Gagliardo-Nirenberg inequality and $\xi_{1,\ve}\in(\f12,1)$ for $r\in[0,1]$ that
$$\|\bar{p}\|^2_{L^\infty([0,1])}\leq \nu_0\int_{0}^{1}\xi_{1,\ve}^2(\pa_r\bar{p})^2dr+C(\nu_0)\int_{0}^{1}\bar{p}^2dr.$$
For $r\in [-1,0]$, since $\xi_{1,\ve}\leq\f12$, we have by \eqref{est:p(a)-1} that
\begin{align*}
|p(r)|&=\xi_{1,\ve}(r)|\bar{p}(r)|\leq \xi_{1,\ve}(r)\Big(|\bar{p}(0)|+\int_{r}^{0}|\pa_r\bar
{p}|dr\Big)\nonumber\\
&\leq \frac{1}{2}|\bar{p}(0)|+\xi_{1,\ve}(r)\Big(\int_{r}^{0}\xi_{1,\ve}^2|\pa_r\bar{p}|^2dr\Big)^{\frac{1}{2}}\Big(\int_{r}^{0}\xi_{1,\ve}^{-2}dr\Big)^{\frac{1}{2}}\nonumber\\
&\leq \frac{1}{2}|\bar{p}(0)|+O(\sqrt{\ve})\Big(\int_{r}^{0}\xi_{1,\ve}^2|\pa_r\bar{p}|^2dr\Big)^{\frac{1}{2}}.
\end{align*}
It then follows from Lemma \ref{le:1} that \eqref{infty} holds.
For $j=2,$ one can infer that $\mathcal{G}_2(p)\geq \int_{I}|\pa_rp|^2dr$ as a fact of $\iota(s)> 1$. Similarly, we have
$$\|p\|^2_{L^\infty[-1,1]}\leq \nu_0\mathcal{G}_2(p)+C(\nu_0)\int_{I}p^2dr.$$
The proof is completed.
\end{proof}
\end{lemma}
\begin{lemma}\label{le:endpoints}
For $j=0,$ there exists $C>0$ such that
\begin{align}\label{endp}
|p(\pm 1)|^2\leq C\ve\Big(\mathcal{G}_0(p)+\int_{I}p^2dr\Big).
\end{align}
\proof
Let $p=\xi_{0,\ve}\bar{p}$. It is straightforward to check that $\xi_{0,\ve}^2(\pm 1)\leq C\ve |(\xi_{0,\ve}\partial_r\xi_{0,\ve})(\pm 1)|$. \eqref{es:0f} and Lemma \ref{le:0} yield \eqref{endp}.
\qed
\end{lemma}
Lemmas \ref{leminfty} and \ref{le:endpoints} give the following corollary for the endpoint estimates.
\begin{Corollary}\label{cor:p(pm1)}
Let $i=0,\cdots,4$, and let $j(i)\in\{0,1,2\}$ be defined as \eqref{def:j-al}, $0<\nu_0\ll 1$. It follows for $\ve$ small enough such that
\begin{align}\label{endp-al}
|p_{i}(\pm 1)|^2\leq \nu_0\mathcal{G}_{j(i)}(p_{i})+C(\nu_0)\int_{I}p_{i}^2dr.
\end{align}
\end{Corollary}
\subsection{Estimate for convolution terms}
In this section, to prove \eqref{prove-Linear} separately for $i=1,\cdots,4$, we only need to prove Lemmas \ref{lem:tq1-1}--\ref{lem:tq1-3}.
\begin{lemma}\label{lem:tq1-1}
Let $i=0$. For sufficiently small $\ve$, we have
\begin{equation}\label{lem:tq1mainineq}
\int_I (\ve^2+r^2)(\pa_rp_{0})^2dr\leq \nu_0\mathcal{G}_{0}(p_{0})+C\int_Ip_{0}^2dr.
\end{equation}
\end{lemma}
\begin{proof}
The proof follows steps similar to those in Lemma \ref{le:0}.
Set
$$p_0=\xi_{0,\ve}\bar{p}_0,\ \gamma =\big(1+\frac{2L}{3}\big)^{-\f12}.$$
Similar as \eqref{est:G0-inte}, we get by integration by parts that
\begin{align}
\notag\int_I (\ve^2+r^2)(\pa_rp_0)^2dr
&=\int_I(\ve^2+r^2)\xi_{0,\ve}^2
(\pa_r\bar{p}_0)^2dr+(\ve^2+r^2)\xi_{0,\ve}\pa_r\xi_{0,\ve}
\bar{p}_0^2\big|_{-1}^1\notag\\
&\quad-\int_I(\ve^2+r^2)\xi_{0,\ve}\pa_r^2\xi_{0,\ve}
\bar{p}_0^2dr-\int_I2r\xi_{0,\ve}
\pa_r\xi_{0,\ve}\bar{p}_0^2dr\notag\\
&=\int_I(\ve^2+r^2)\xi_{0,\ve}^2
(\pa_r\bar{p}_0)^2dr+(\ve^2+r^2)\xi_{0,\ve}\pa_r\xi_{0,\ve}
\bar{p}_0^2\big|_{-1}^1\notag\\
&\quad-\gamma\int_I\Big(\gamma\big(1+(\f{r}{\ve})^2\big)
\theta(\xi_{1,\ve})+\f{2r}{\ve}(1-2\xi_{1,\ve})\Big)
p_0^2dr,\label{fm:ve2r2parp0}
\end{align}
where we use $\partial_r^2\xi_{0,\ve}
=\frac{1}{\ve^2}\gamma^2 \xi_{0,\ve}\theta(\xi_{1,\ve})$ and $\pa_r\xi_{0,\ve}=\frac{1}{\ve}\gamma \xi_{0,\ve}(1-2\xi_{1,\ve})$.
Thanks to Lemma \ref{le:0}, 
we have
\begin{align}\label{est:p01}
\int_I(\ve^2+r^2)\xi_{0,\ve}^2
(\pa_r\bar{p}_0)^2dr\leq \frac{\nu_0}{3}\mathcal{G}_0(p_0)+C\int_Ip_0^2dr.
\end{align}
Therefore, we use \eqref{es:0f} (and the same bound at $r=1$) to obtain
\begin{align}\label{est:p02}
(\ve^2+1)|\xi_{0, \ve}\partial_r\xi_{0, \ve}\bar{p}_0^2(\pm1)|
&\leq \int_{I} \xi_{0, \ve}^2(\pa_r\bar{p}_0)^2 dr+C\int_{I}\xi_{0, \ve}^2\bar{p}_0^2dr\notag\\
&\leq \frac{\nu_0}{3}\mathcal{G}_0(p_0)+C\int_Ip_0^2dr.
\end{align}
We claim
\begin{align}\label{est:p03}
A:=\int_I\Big(\gamma\big(1+(\f{r}{\ve})^2\big)\theta(\xi_{1,\ve})
+\f{2r}{\ve}(1-2\xi_{1,\ve})\Big)p_0^2dr\geq -C\int_Ip_0^2dr,
\end{align}
which together with  \eqref{fm:ve2r2parp0}--\eqref{est:p02} gives \eqref{lem:tq1mainineq}. It remains to prove \eqref{est:p03}. Indeed, we observe from
$\theta(\xi_{1,\ve})=1-6\xi_{1,\ve}+6\xi_{1,\ve}^2$ and $\xi_{1,\ve}(r)=\f{1}{1+e^{-\f{\gamma r}{\ve}}}\in(0,1)$ that there exists $M>0$ such that
\begin{align*}
\theta(\xi_{1,\ve}(r))\geq \f12,\;\;|\f{r}{\ve}|\geq M.
\end{align*}
Together with $|1-2\xi_{1,\ve}|\leq 1$, it follows by Young's inequality that
\begin{align*}
\Big(\gamma\big(1+(\f{r}{\ve})^2\big)\theta(\xi_{1,\ve})
+\f{2r}{\ve}(1-2\xi_{1,\ve})\Big)
\geq \f{\gamma}{2}\big(1+|\f{r}{\ve}|^2\big)-2|\f{r}{\ve}|\geq \f{\gamma}{2}-\f{2}{\gamma},\;\;\;|\f{r}{\ve}|\geq M.
\end{align*}
Therefore, we can decompose $A$ in \eqref{est:p03} as
\begin{align*}
A=\int_{\{|\f{r}{\ve}|\geq M\}}+\int_{\{|\f{r}{\ve}|\leq M\}}:=A_1+A_2,
\end{align*}
and obtain
\begin{align*}
A_1\geq -C\int_{\{|\f{r}{\ve}\}|\geq M}p_0^2dr\geq -C\int_{I}p_0^2dr,\quad\text{and}\quad
|A_2|\leq C\int_{\{|\f{r}{\ve}|\leq M\}}p_0^2dr\leq C\int_{I}p_0^2dr.
\end{align*}
The estimate of $A_2$ uses $|\f{r}{\ve}|\leq M$. Then we have $A\geq -C\int_{I}p_0^2dr$.
\end{proof}
\begin{lemma}\label{lem:tq1-2}
Let $i=1,2$. For sufficiently small $\ve$, we have
\begin{equation}\label{lem:tq1mainineq2}
\int_I (\ve^2+r^2)(\pa_rp_{i})^2dr\leq \nu_0\mathcal{G}_{1}(p_{i})+C\int_Ip_{i}^2dr.
\end{equation}
\end{lemma}
\begin{proof}
Set $p_{i}=\xi_{1,\ve}\bar{p}_{i}$, $\gamma =\big(1+\frac{2L}{3}\big)^{-\frac{1}{2}}$. Similar as \eqref{trans}, we get by integration by parts that
\begin{align}
\notag\int_I (\ve^2+r^2)(\pa_rp_{i})^2dr
&=\int_I(\ve^2+r^2)\xi_{1,\ve}^2
(\pa_r\bar{p}_{i})^2dr+(\ve^2+r^2)\xi_{1,\ve}\pa_r\xi_{1,\ve}
\bar{p}_{i}^2\big|_{-1}^1\notag\\
&\quad-\int_I(\ve^2+r^2)\xi_{1,\ve}\pa_r^2\xi_{1,\ve}
\bar{p}_{i}^2dr-\int_I2r\xi_{1,\ve}
\pa_r\xi_{1,\ve}\bar{p}_{i}^2dr\notag\\
&=\int_I(\ve^2+r^2)\xi_{1,\ve}^2
(\pa_r\bar{p}_{i})^2dr+(\ve^2+r^2)\f{\gamma(1-\xi_{1,\ve})}{\ve}
p_{i}^2\big|_{-1}^1\notag\\
&\quad-\gamma\int_I\Big(\gamma \big(1+(\f{r}{\ve})^2\big)
\kappa(\xi_{1,\ve})+\f{2r}{\ve}(1-\xi_{1,\ve})\Big)
p_i^2dr,\label{fm:ve2r2parp1}
\end{align}
where we have used $\pa_r^2\xi_{1,\ve}=\frac{1}{\ve^2}\gamma^2 \xi_{1,\ve}\kappa(\xi_{1,\ve})
$ and $\pa_r\xi_{1,\ve}=\frac{\gamma }{\ve}(1-\xi_{1,\ve}) \xi_{1,\ve}$.
Thanks to Lemma \ref{le:1}, 
we have
\begin{align}\label{est:p11}
\int_I(\ve^2+r^2)\xi_{1,\ve}^2
(\pa_r\bar{p}_{i})^2dr\leq \int_I\xi_{1,\ve}^2
(\pa_r\bar{p}_{i})^2dr\leq \frac{\nu_0}{2}\mathcal{G}_1(p_{i})+C\int_Ip_{i}^2dr.
\end{align}
Using $1-\xi_{1,\ve}>0$, $\f{1-\xi_{1,\ve}(1)}{\ve}=\f{\ve^{-1}e^{-\f{\gamma }{\ve}}}{1+e^{-\f{\gamma }{\ve}}}\leq C$ and Corollary \ref{cor:p(pm1)} ($i=1,2$, $j(i)=1$), we have
\begin{align}\label{est:p12}
(\ve^2+r^2)\f{\gamma(1-\xi_{1, \ve})}{\ve}
p_{i}^2\big|_{-1}^1&\leq (\ve^2+1)\f{\gamma(1-\xi_{1, \ve}(1))}{\ve}
p_{i}^2(1)\notag\\
&\leq \frac{\nu_0}{2}\mathcal{G}_1(p_{i})+C\int_Ip_{i}^2dr.
\end{align}
We claim that
\begin{align}\label{est:p13}
A:=\int_I\Big(\gamma \big(1+(\f{r}{\ve})^2\big)
\kappa(\xi_{1,\ve})+\f{2r}{\ve}(1-\xi_{1,\ve})\Big)p_{i}^2dr\geq -C\int_Ip_{i}^2dr,
\end{align}
which together with  \eqref{fm:ve2r2parp1}--\eqref{est:p12} gives \eqref{lem:tq1mainineq2}. It remains to prove \eqref{est:p13}. Indeed, we observe from
$\kappa(\xi_{1,\ve})=1-3\xi_{1,\ve}+2\xi_{1,\ve}^2$ and $\xi_{1,\ve}(r)=\f{1}{1+e^{-\f{\gamma r}{\ve}}}\in(0,1)$ that there exists $M>0$ such that
\begin{align*}
&\kappa(\xi_{1,\ve}(r))\geq \f12,\;\;\quad\f{r}{\ve}\leq -M.
\end{align*}
and
\begin{align*}
&-\f12e^{-\f{\gamma r}{\ve}}\leq\kappa(\xi_{1,\ve}(r))\leq 0 ,\quad\quad 0<1-\xi_{1,\ve}(r)\leq e^{-\f{\gamma r}{\ve}},\;\;\quad r>0.
\end{align*}
Together with $|1-\xi_{1,\ve}|\leq 1$, it follows by Young's inequality that
\begin{align*}
\Big(\gamma \big(1+(\f{r}{\ve})^2\big)\kappa(\xi_{1,\ve})
+\f{2r}{\ve}(1-\xi_{1,\ve})\Big)
\geq \f{\gamma }{2}\big(1+|\f{r}{\ve}|^2\big)-2|\f{r}{\ve}|\geq
\f{\gamma }{2}-\f{2}{\gamma },\;\;\;\f{r}{\ve}\leq -M,
\end{align*}
and
\begin{align*}
\Big|\gamma \big(1+(\f{r}{\ve})^2\big)\kappa(\xi_{1,\ve})
+\f{2r}{\ve}(1-\xi_{1,\ve})\Big|
\leq \Big( \gamma \big(1+|\f{r}{\ve}|^2\big)+2|\f{r}{\ve}|\Big)e^{-\f{\gamma r}{\ve}}\leq C
,\;\;\;r>0.
\end{align*}
Therefore, we can decompose $A$ in \eqref{est:p13} as
\begin{align*}
A=\int_{\{\f{r}{\ve}\leq -M\}}+\int_{\{-M\leq \f{r}{\ve}\leq 0\}}+\int_{\{r>0\}}:=A_1+A_2+A_3,
\end{align*}
and obtain
\begin{align*}
A_1\geq -C\int_{\{\f{r}{\ve}\leq -M\}}p_{i}^2dr&\geq -C\int_{I}p_{i}^2dr,\quad
|A_2|\leq C\int_{\{-M\leq \f{r}{\ve}\leq 0\}}p_{i}^2dr\leq C\int_{I}p_{i}^2dr,\\
&\text{and}\quad\; |A_3|\leq C\int_{0}^1p_{i}^2dr\leq C\int_{I}p_{i}^2dr .
\end{align*}
Note that the estimate of $A_2$ uses $|\f{r}{\ve}|\leq M$. Thus we have $A\geq -C\int_{I}p_{i}^2dr$.
\end{proof}
\begin{lemma}\label{lem:tq1-3}
Let $i=3,4$. For sufficiently small $\ve$, we have
\begin{equation}\label{lem:tq1mainineq34}
\int_I (\ve^2+r^2)(\pa_rp_{i})^2dr\leq \nu_0\mathcal{G}_{2}(p_{i})+C\int_Ip_{i}^2dr.
\end{equation}
\end{lemma}
\begin{proof}
The proof follows directly by \eqref{lowerG2}.
\end{proof}

\subsection{Estimate for cross terms}
Denote $a(r)=E^\al:\pa_rE^\beta$. \eqref{prove-Cross} can be reduced to the following three lemmas.
\begin{lemma}\label{cross1}
For any $\nu_0\in(0,\f12)$, there exists $C=C(\nu_0,\|a(r)\|_{W^{1,\infty}})>0$ such that
\begin{align}
\int_{I}(p_\al\pa_rp_\beta-p_\beta\pa_r p_\al )E^\al:\pa_rE^\beta dr\leq \nu_0 \big(\mathcal{G}_{2}(p_\al)+\mathcal{G}_{2}(p_\beta)\big)
+C\int_{I}\big(p_\al^2+p_\beta^2\big)dr.
\end{align}
\begin{proof}
We get by integration by parts, Young's inequality and Corollary \ref{cor:p(pm1)} that
\begin{align*}
\int_{I}\big(p_{\al} \pa_rp_{\beta} -p_{\beta}\pa_rp_{\al} \big)a(r) dr&=-2\int_{I} p_{\beta}\pa_rp_\al a(r) dr-\int_{I}p_\al p_\beta a'(r)dr+p_\al p_\beta a(r)|_{-1}^{1}\nonumber\\
&\leq \f{\nu_0}{3}\int_{I}|\pa_rp_\al |^2dr+C\int_{I}\big(p_\al ^2+p_\beta^2\big)dr+p_\al p_\beta a(r)|_{-1}^1\notag\\
&\leq \nu_0\mathcal{G}_{2}(p_\al)+\nu_0\mathcal{G}_{2}(p_{\beta})
+C\int_{I}\big(p_\al ^2+p_\beta^2\big)dr.\notag
\end{align*}
Indeed, we have used \eqref{lowerG2}. The lemma follows.
\end{proof}
\end{lemma}
\begin{lemma}\label{cross2}
For any $\nu_0\in(0,\f12)$, there exists $C(\nu_0,\|a(r)\|_{W^{1,\infty}})>0$ such that
\begin{align}\label{cross2maineq}
\int_{I}(p_\al \partial_rp_\beta-p_\beta \partial_rp_\al )a(r)dr\leq \nu_0\big(\mathcal{G}_1(p_\al )+\mathcal{G}_1(p_\beta )\big)+C\int_{I}(p_\al^2+p_\beta^2)dr.
\end{align}
\proof
Let
$$p_\al(r)=\xi_{1,\ve}(r)\bar{p}_\al(r),\ p_\beta(r)=\xi_{1,\ve}(r)\bar{p}_\beta(r).$$
Then one deduces by  Young's inequality that
\begin{align*}
\int_{I}(p_\al\partial_rp_\beta-\partial_rp_\al p_\beta)a(r) dr&=\int_{I}\xi_{1,\ve}^2(\bar{p}_\al\partial_r\bar{p}_\beta
-\partial_r\bar{p}_\al\bar{p}_\beta
)a(r)dr\\
&\leq \f{\nu_0}{4}\int_{I}\xi_{1,\ve}^2\left((\partial_r\bar{p}_\al)^2
+(\partial_r\bar{p}_\beta)^2\right)dr+C\int_I\big(
p_\al^2+p_\beta^2\big)dr.
\end{align*}
As Lemma \ref{le:1} gives $\mathcal{G}_1(p_{\gamma})\geq \f14 \int_{I}\xi_{1,\ve}^2(\partial_r\bar{p}_{\gamma})^2
dr-C\int_{I}p_{\gamma}^2dr$ for $\gamma\in\{\al,\beta\}$,  then \eqref{cross2maineq} follows.
\qed
\end{lemma}
\begin{lemma}\label{cross3}
For any $\delta_0,\nu_0\in(0,\f12)$, there exists $C(\nu_0,\delta_0,\|a(r)\|_{W^{1,\infty}})$ such that
\begin{align}
\int_{I}(p_0 \partial_rp_\beta-p_\beta \partial_rp_0 )a(r) dr\leq \f{\delta_0}{4}\int_{I}(\pa_rp_\beta)^2dr+\nu_0\mathcal{G}_0(p_0)
+C\int_{I}\big(p_0^2+p_\beta^2\big)dr.
\end{align}
\begin{proof}
 It follows from Corollary \ref{cor:p(pm1)} that
 \begin{align*}
|p_0(\pm1)|^2\leq \nu_0\mathcal{G}_0(p_0)+C\int_{I}p_0^2dr.
\end{align*}
It follows from  Gagliardo-Nirenberg inequality and Young's inequality that
\begin{align*}
|p_{\beta}(\pm1)|^2\leq\f{\delta_0}{8}
\int_{I}(\pa_rp_{\beta})^2dr+C(\delta_0)\int_{I}p_{\beta}^2dr.
\end{align*}
Thus, we integrate by parts to get
\begin{align}
\int_{I}&(p_0 \partial_rp_\beta-p_\beta \partial_rp_0 )a(r) dr\notag\\
&=2\int_{I}p_0\pa_rp_\beta a(r)dr-\int_{I}p_0p_\beta a'(r)dr-p_0p_\beta a(r) |_{-1}^1\notag\\
&\leq \f{\delta_0}{4}\int_{I}|\pa_rp_\beta|^2dr+\nu_0\mathcal{G}_0(p_0)
+C(\nu_0,\delta_0,\|a(r)\|_{W^{1,\infty}})\int_{I}\big(p_\beta^2+p_0^2\big)dr\notag.
\end{align}
\end{proof}
\end{lemma}

\subsection{Estimate for correction terms}\label{6.6}
 By \eqref{transformation} and \eqref{dK=r}, the functions $s_{1,i}(\f{d^K(x,t)}{\ve},x,t)\ (i=0,\cdots,4)$ defined in \eqref{q1near} can be rewritten  as $s_{1,i}(\f{r}{\ve},x(\sigma,r),t)$.
We recall correction terms as the following:
\begin{align*}
\frac{1}{\ve}\int_I&f''_{Q_0}(P, P):Q_1dr\\
=&\frac{1}{\ve}\int_I\Big(4(3+2s)s_{1, 0}(p_3^2+p_4^2)-36s_{1, 1}(p_1p_4+p_3p_2)-36s_{1, 2}(p_1p_3-p_4p_2)\\
&+8(3+2s)p_0(s_{1, 3}p_3+s_{1, 4}p_4)\Big)dr\\
&+\frac{1}{\ve}\int_I\Big(
2(-3+4\xi_{1, \ve})(p_1^2+p_2^2)s_{1, 0}
+\f23(-6+12\xi_{1, \ve})p_0^2s_{1, 0}\\
&+18(p_2^2-p_1^2)s_{1, 4}-36p_1p_2s_{1, 3}\Big)dr+\frac{4}{\ve}\int_I\frac{d}{ds}
\kappa(\xi_{1, \ve})p_0(s_{1, 1}p_1+s_{1, 2}p_2)dr\\
:=&I_{G_{1}}+I_{G_{2}}+I_B.
\end{align*}
We have by \eqref{lowerG2} that
$
\ve^{-2}\int_Ip_{i}^2dr\leq \mathcal{G}_{2}(p_{i})$ for $i=3,4$,
which gives:
\begin{align}\label{Correct-G1}
|I_{G_1}| &
\leq \nu_0\sum_{i=3,4}\mathcal{G}_2(p_{i})
+C(\nu_0)\sum_{i=0}^4\int_{I}p_{i}^2dr,
\end{align}
for small $\ve$ and any $\nu_0\in(0, 1/2)$.

On the other hand, by \eqref{s1-0:estimate}--\eqref{s1-12:estimate} and the fact $d_0=r+O(\ve)$,  we have \begin{align*}
\frac{1}{\ve}\Big|s_{1,i}(\frac{r}{\ve},x,t)\Big|\leq C\frac{|d_0|}{\ve}e^{-\gamma\frac{|r|}{\ve}}\leq C,\;\;i=0,3,4.
\end{align*}
Thus,
\begin{align}
|I_{G_2}|&\leq C\sum_{i=0}^4\int_{I}p_{i}^2dr.\label{1vepartial}
\end{align}

For the last term $I_B$, we have the following estimate.
\begin{lemma}\label{correction}
For any $0<\delta_0,\nu_0<\f12$, there exists $C>0$ such that
\begin{align*}
|I_B|\leq \f12\delta_0\int_{I}|\pa_rp_1|^2+|\pa_rp_2|^2dr+\nu_0\Big(\mathcal{G}_0(p_0)+\sum_{i=1,2}\mathcal{G}_1(p_{i})\Big)+C\sum_{i=0}^4\int_{I}p_i^2dr.
\end{align*}
\end{lemma}
The key of the proof is the following cancellation structure.
\begin{lemma}\label{Ove}
For $i=1,2$, it holds that
\begin{align*}
\int_{I}\frac{d}{ds}\kappa(  \xi_{1,\ve})s_{1,i}\xi_{0,\ve}\xi_{1,\ve}dr \le C\ve.
\end{align*}
\begin{proof}It suffices to consider $i=1$. We have by changing variable $z=\frac{r}{\ve}$ that
\begin{align*}
\int_{I}\frac{d}{ds}\kappa(  \xi_{1,\ve})s_{1,1}\xi_{0,\ve}\xi_{1,\ve}dr&=\ve\int_{-\frac{1}{\ve}}^{\frac{1}{\ve}}(-3+4s)ss's_{1,1}dz\\
&=\ve\int_{|z|\leq 1}(-3+4s)ss's_{1,1}dz+\ve\int_{1\leq|z|\leq \frac{1}{\ve}}(-3+4s)ss's_{1,1}dz\\
&:=I_1+I_2.
\end{align*}
Thanks to Lemma \ref{lem2}, we have $|s_{1,1}|\leq C$. Then it follows $|I_1|\leq C\ve$.
Furthermore, using the exponential decay of $s'(z)$ in \eqref{est:s(z)}, we have
$$|I_2|\leq C\ve \int_{1\leq |z|\leq \frac{1}{\ve}}e^{-\gamma |z|}dz\leq C\ve.$$
Summing up $I_1$ and $I_2$, the estimate follows.
\end{proof}
\end{lemma}
\begin{proof}[Proof of Lemma \ref{correction}]
Set
\begin{align*}
p_0=\mu \xi_{0,\ve}+\hat{p}_0\text{ and }p_1=\kappa_1 \xi_{1,\ve}+\hat{p}_1\text{ with }\int_{I}\xi_{0,\ve}\hat{p}_0dr=0\text{ and } \int_{I}\xi_{1,\ve}\hat{p}_1\omega_\ve(r)dr=0.
\end{align*}
Let $$\|p\|_{\omega}=\Big(\int_{I}{\omega_\ve} p^2dr\Big)^{\frac{1}{2}}\text{ and }\|p\|_{L^2}=\Big(\int_{I}p^2dr\Big)^{\frac{1}{2}}.$$
We directly get
\begin{align}
\frac{4}{\ve}\int_{I}\frac{d}{ds}\kappa(  \xi_{1,\ve})p_0s_{1,1}p_1dr=&\frac{4\kappa_1\mu}{\ve}\int_{I}\frac{d}{ds}\kappa(  \xi_{1,\ve})s_{1,1}\xi_{0,\ve}\xi_{1,\ve}dr
+\frac{4\mu}{\ve}\int_{I}\frac{d}{ds}\kappa(  \xi_{1,\ve})s_{1,1}\xi_{0,\ve}\hat{p}_1dr\nonumber\\
&+\frac{4}{\ve}\int_{I}\frac{d}{ds}\kappa(  \xi_{1,\ve})s_{1,1}\hat
{p}_0
p_1dr\nonumber\\
:=&I_1+I_2+I_3.\nonumber
\end{align}
We deduce from $p_0=\mu \xi_{0,\ve}+\hat{p}_0$ with $\int_I\xi_{0,\ve}\hat{p}_0dr=0$ that
\begin{align*}
\int_{I}p_0^2dr=&\int_{I}\Big(\mu\xi_{0,\ve}+\hat{p}_0\Big)^2dr\\
=&\mu^2\int_{I}\xi_{0,\ve}^2dr+\int_{I}|\hat{p}_0|^2dr.
\end{align*}
Then we have
$$|\mu|\leq \frac{C}{\sqrt{\ve}}\|p_0\|_{L^2}.$$
By applying
 $$\int_{I}\omega_\ve|\pa_r p_1|^2dr=\kappa_1^2\int_{I}\omega_\ve|\pa_r \xi_{1,\ve}|^2dr+\int_{I}\omega_\ve|\pa_r \hat{p}_1|^2dr+2\kappa_1\int_{I}\omega_\ve \pa_r\hat{p}_1\pa_r\xi_{1,\ve}dr$$
and the Cauchy-Schwarz inequality, we have
$$|\kappa_1|\leq C\sqrt{\ve}(\|\pa_r p_1\|_\omega+\|\pa_r \hat{p}_1\|_\omega).$$
By applying Lemma \ref{g1second} and Lemma \ref{Ove}, one has
\begin{align}
|I_1|&=\frac{4}{\ve}\Big|\kappa_1\mu\int_{I}\frac{d}{ds}\kappa(  \xi_{1,\ve})s_{1,1}\xi_{0,\ve}\xi_{1,\ve}dr\Big|\notag\\
&\leq C\big(\|p_0\|_{L^2}\|\pa_r p_1\|_{\omega}+C\|p_0\|_{L^2}\|\pa_r \hat{p}_1\|_{\omega}\big)\frac{1}{\ve}\Big|\int_{I}\frac{d}{ds}\kappa(  \xi_{1,\ve})s_{1,1}\xi_{0,\ve}\xi_{1,\ve}dr\Big|\nonumber\\
&\leq \f12\delta_0\int_{I}|\pa_r p_1|^2dr+\|\pa_r\hat{p}_1\|_{\omega}^2+C\int_{I}p_0^2+p_1^2dr\nonumber\\
&\leq \f12\delta_0\int_{I}|\pa_r p_1|^2dr+\frac{1}{4}\mathcal{G}_1(p_1)+C\int_{I}p_0^2+p_1^2dr.\label{s112}
\end{align}
From Young's inequality,  H$\ddot{o}$lder's inequality and Lemma \ref{g1second}, we deduce that
\begin{align}
|I_2|&\leq C\frac{\mu}{\ve}\Big(\int_I|\xi_{0,\ve}|dr\Big)^{\f12}\|\hat{p}_1\|_\omega\notag\\
&\leq C\mu^2 \int_I \xi_{0,\ve}^2dr+\frac{1}{\ve^2}\|\hat{p}_1\|_\omega\notag\\
&\leq C\int_I p_0^2+p_1^2dr+\nu_0\mathcal{G}_1(p_1).\label{ibi2}
\end{align}
Indeed, we take $\omega_\ve(r)=|\big(\frac{d}{ds}\kappa(\xi_{1,\ve})\big)(r)s_{1,1}
(\f{r}{\ve},x,t)|^2_\G|
\xi_{0,\ve}(r)|$, which has at most finite zeros on $I$, as the non-zero function $s_{1,1}|_{\Gamma}(z)$ is real analytic on any finite interval by the Remark \ref{remark:analyitical}.
Moreover, it follows by Young's inequality and Lemma \ref{secondeigenvalue} that
\begin{align}\label{ibi3}
|I_3|\leq \frac{1}{\ve^2}\int_I\hat{p}_0^2dr
+C\int_{I}p_1^2dr
\leq \frac{\nu_0}{2}\mathcal{G}_0(p_0)+C\int_{I}p_0^2+p_1^2dr.
\end{align}
Combining \eqref{s112}, \eqref{ibi2} and \eqref{ibi3}, we have
$$\Big|\frac{4}{\ve}\int_{I}\frac{d}{ds}\kappa(  \xi_{1,\ve})p_0s_{1,1}p_1dr\Big|\leq \f12\delta_0\int_{I}|\pa_r p_1|^2dr+\nu_0\big(\mathcal{G}_0(p_0)+\mathcal{G}_1(p_1)\big)+C\sum_{i=0}^4\int_{I}p_i^2dr.$$
The estimate of $\frac{1}{\ve}\Big|\int_{I}\frac{d}{ds}\kappa(  \xi_{1,\ve})p_0s_{1,2}p_2dr\Big|$ follows in a similar manner.
\end{proof}

\bigskip
\bigskip

\section{Uniform error estimates}\label{section7}
Define
$$Q_R:=\frac{Q^\ve-Q^K}{\ve^k},$$
 where $Q^\ve$ is the solution of \eqref{equation:main} and $Q^K$ is the approximate solution constructed in Theorem \ref{th1}. Using \eqref{equation:main} and \eqref{equm}, we derive that $Q_R$ is the solution of
\begin{equation}\label{error}
\pa_t Q_R=\mathcal{L}Q_R-\frac{1}{\ve^2}\left(\mathcal{H}_{Q^K}Q_R+\frac{\ve^k}{2}f''_{Q^K+\frac{1}{3}\ve^kQ_R}(Q_R,Q_R)\right)+\mathcal{R}_k,
\end{equation}
where
\begin{equation}\label{pamrk}
\pa^i\mathcal{R}_k=O(\ve^{K-k-1-3i}) \text{ for } i\geq 0.
\end{equation}
We choose $k=9$ and $K\geq 10$. Set
$$\mathcal{E}(Q)=\sum_{i=0}^{2}\ve^{6i}\int_{\Omega}\|\pa^iQ\|^2dx.$$
It is obvious  that
\begin{align}\label{qinfty}
\|Q\|_{L^\infty(\Omega)}\leq \|Q\|_{H^2(\Omega)}\leq \ve^{-6}\mathcal{E}(Q)^{\frac{1}{2}},\ \ \|\nabla Q\|_{L^4(\Omega)}\leq \|Q\|_{H^2(\Omega)}\leq \ve^{-6}\mathcal{E}(Q)^{\frac{1}{2}}.
\end{align}

 Multiplying \eqref{error} by $Q_R$ and integrating over $\Omega$, one has
\begin{align}
\frac{1}{2}\pa_t\int_{\Omega}Q_R:Q_Rdx=&-\int_{\Omega}|\nabla Q_R|^2dx-L\int_{\Omega}|\nabla\cdot Q_R|^2dx\nonumber\\
&-\frac{1}{\ve^2}\int_{\Omega}\left(\mathcal{H}_{Q^K}Q_R+\frac{1}{2}\ve^kf''_{Q^K+\frac{1}{3}\ve^9Q_R}(Q_R, Q_R)\right):Q_Rdx\notag\\
&+\int_{\Omega}\mathcal{R}_9:Q_Rdx.
\end{align}
 Based on the spectral estimate Theorem \ref{th:uA}, \eqref{pamrk} and \eqref{qinfty}, we can establish the inequality:
\begin{align*}
\pa_t\int_{\Omega}|Q_R|^2dx\leq& -\int_{\Omega}|\nabla Q_R|^2dx-L\int_{\Omega}|\nabla\cdot Q_R|^2dx-\frac{1}{\ve^2}\int_{\Omega}\mathcal{H}_{Q^K}Q_R:Q_Rdx\\
&-\frac{1}{2}\ve^{7}\int_{\Omega}f''_{Q^K+\frac{1}{3}\ve^9Q_R}(Q_R,Q_R):Q_Rdx+C\ve^{K-10}\int_{\Omega}|Q_R|dx\\
\leq &C(1+\ve\mathcal{E}(Q_R)^{\frac{1}{2}})\int_{\Omega}|Q_R|^2dx+C\ve^{K-10}\mathcal{E}(Q_R)^{\frac{1}{2}}\\
\leq &C(1+\mathcal{E}(Q_R)+\ve\mathcal{E}(Q_R)^{\frac{3}{2}}).
\end{align*}
For $i=1$, differentiating both sides of \eqref{error} with respect to $x_m$ for $m\in \{1,2,3\}$, we have
\begin{align}
\pa_t\pa_mQ_R=&\pa_m\mathcal{L}Q_R-\frac{1}{\ve^2}\left(\mathcal{H}_{Q^K}\pa_mQ_R+\ve^9f''_{Q^K+\frac{1}{3}\ve^9Q_R}(\pa_mQ_R,Q_
R)\right)+\pa_m\mathcal{R}_9\nonumber\\
&-\ve^{-2}f''_{Q^K}(\pa_m Q^K,Q_R)-3\ve^7\Big(2Q_R(\pa_mQ^K+\frac{1}{3}\ve^9\pa_mQ_R):Q_R\nonumber\\
&+(\pa_mQ^K+\frac{1}{3}\ve^9\pa_mQ_R)|Q_R|^2\Big).
\end{align}
Integrating the above equation over $\Omega$ after contracting with $\pa_mQ_R,$ we have
\begin{align}
\frac{1}{2}\pa_t\int_{\Omega}|\pa_mQ_R|^2dx=&-\int_{\Omega}|\nabla \pa_m Q_R|^2dx-L\int_{\Omega}|\nabla\cdot \pa_m Q_R|^2dx\nonumber\\
&-\int_{\Omega}\frac{1}{\ve^2}\Big(\mathcal{H}_{Q^K}\pa_mQ_R+\ve^9 f''_{Q^K}(\pa_mQ_R, Q_
R)\Big):\pa_m Q_Rdx\nonumber\\
&-\int_{\Omega}\Big(\ve^{-2}\mathcal{H}_{Q^K}(\pa_m Q^K, Q_R)\notag\\
&+3\ve^{7}\big(2Q_R(\pa_mQ^K+\frac{1}{3}\ve^9\pa_mQ_R):Q_R\nonumber\\
&+(\pa_mQ^K+\frac{1}{3}\ve^9\pa_mQ_R)|Q_R|^2\Big):\pa_mQ_Rdx\notag\\
&+\int_{\Omega}\pa_m\mathcal{R}_9:\pa_mQ_Rdx\nonumber\\
&-\int_{\Omega}\ve^{7}f''_{Q^K+\frac{1}{3}\ve^9Q_R}(\pa_mQ_R, Q_
R):\pa_m Q_Rdx.
\end{align}
According to Theorem \ref{th:uA}, \eqref{pamrk} and \eqref{qinfty}, we deduce that
\begin{align}
\ve^6\pa_t\int_{\Omega}|\nabla Q_R|^2dx\leq & C\ve^6\int_{\Omega}|\nabla Q_R|^2dx+\ve^{K-5} \left(\int_{\Omega}|\nabla Q_R|^2dx\right)^{\frac{1}{2}}\nonumber\\
& +C\ve^3\int_{\Omega}|Q_R||\nabla Q_R|dx+C\ve^{13}\int_{\Omega}|\nabla Q_R|^2|Q_R|dx\nonumber\\
&+C\ve^{12}\int_{\Omega}|Q_R|^2|\nabla Q_R|dx+C\ve^{22}\int_{\Omega}|Q_R|^2|\nabla Q_R|^2dx\nonumber\\
\leq & C\ve^3\left(\int_{\Omega}|\nabla Q_R|^2dx\right)^{\frac{1}{2}}\Big(\mathcal{E}^{\frac{1}{2}}+\ve\mathcal{E}+\ve^3\mathcal{E}+\ve^4\mathcal{E}^{\frac{3}{2}}\Big)\nonumber\\
\leq & C(1+\mathcal{E}+\ve\mathcal{E}^{\frac{3}{2}}+\ve^4\mathcal{E}^2).
\end{align}
Similarly, we obtain
\begin{align}
\frac{\ve^{12}}{2}\pa_t\int_{\Omega}|\Delta Q_R|^2dx\leq &C\ve^{12}\int_{\Omega}|\Delta Q_R|^2dx+\ve^{K}(\int_{\Omega}|\Delta Q_R|^2dx)^{\frac{1}{2}}\nonumber\\
&+C\ve^8\int_{\Omega}|\Delta Q_R|  |Q_R|dx+C\ve^9\int_{\Omega}|\nabla Q_R| |\Delta Q_R|dx\nonumber\\
&+C\ve^{19}\int_{\Omega}|\Delta Q_R|^2|Q_R|dx+C\ve^{18}\int_{\Omega}|Q_R| |\nabla Q_R| |\Delta Q_R|dx\nonumber\\
&+C\ve^{17}\int_{\Omega}|Q_R|^2|\Delta Q_R|dx+C\ve^{19}\int_{\Omega}|\nabla Q_R|^2|\Delta Q_R|dx\nonumber\\
&+C\ve^{28}\int_{\Omega}|Q_R|^2|\Delta Q_R|^2dx+C\ve^{28}\int_{\Omega}|\nabla Q_R|^2|\Delta Q_R| |Q_R|dx.
\end{align}
We can show that
\begin{align}
\frac{\ve^{12}}{2}\pa_t\int_{\Omega}|\Delta Q_R|^2dx&\leq C\ve^6\left(\int_{\Omega}|\Delta Q_R|^2dx\right)^{\frac{1}{2}}(1+\mathcal{E}^{\frac{1}{2}}+\ve\mathcal{E}+\ve^{4}\mathcal{E}^{\frac{3}{2}})\nonumber\\
&\leq C(1+\mathcal{E}+\ve\mathcal{E}^{\frac{3}{2}}+\ve^4\mathcal{E}^2).
\end{align}
Summing $i$ from 0 to 2, we have
\begin{equation}
\pa_t \mathcal{E}\leq C(1+\mathcal{E}+\ve\mathcal{E}^{\frac{3}{2}}+\ve^4\mathcal{E}^2).
\end{equation}
Therefore, Theorem \ref{th:error} can be obtained by a standard continuity argument.
\bigskip

\appendix
\section{Properties of the reduced ODEs}
\subsection{Solvability of \eqref{modele}}\label{odesolveappendix}

We consider the solvability of the following system:
 \begin{equation}\label{modele}
 \begin{cases}
\begin{array}{ll}
-(1+A)\pa_y^2u+(1-\hat{s})(1-2\hat{s})u=f,\\
u(-\infty,x,t)=u^-(x,t),\quad u(+\infty,x,t)=u^+(x,t),
\end{array}
\end{cases}
\end{equation}
where $A>0$, $f(y,x,t)$ is a smooth function in $(y,x,t)\in\mathbb{R}\times \G(\delta)$, $u^{\pm}(x,t)$ are  smooth functions in $(x,t)\in\G(\delta)$, and $\hat{s}(y)=\frac{1}{1+{e}^{-y}} $ satisfies
 \begin{equation*}
\left\{
\begin{array}{ll}
\hat{s}''=\hat{s}(1-\hat{s})(1-2\hat{s}),\\
\hat{s}(-\infty)=0,\ \hat{s}(+\infty)=1.
\end{array}
\right.
\end{equation*}
The proofs of the following two lemmas are given at the end of this subsection.
Let $\hat{u}_1,\cdots ,\hat{u}_4$ be the smooth solutions of the following four equations:
\begin{equation}\label{u1}
\left\{
\begin{array}{ll}
(1+A)\pa_y^2\hat{u}_1-(1-\hat{s})(1-2\hat{s})\hat{u}_1=0,\\
\hat{u}_1(y)\sim e^{\f{y}{{\sqrt{1+A}}}},\;\;y\to-\infty.
\end{array}
\right.
\end{equation}
\begin{equation}\label{u2}
\left\{
\begin{array}{ll}
(1+A)\pa_y^2\hat{u}_2-(1-\hat{s})(1-2\hat{s})\hat{u}_2=0,\\
\hat{u}_2(y)\sim e^{-\f{y}{{\sqrt{1+A}}} },\;\;y\to-\infty.
\end{array}
\right.
\end{equation}
 \begin{equation}\label{u3}
\left\{
\begin{array}{ll}
(1+A)\pa_y^2\hat{u}_3-(1-\hat{s})(1-2\hat{s})\hat{u}_3=0,\\
\hat{u}_3(y)\sim 1,\;\;y\to+\infty.
\end{array}
\right.
\end{equation}
 \begin{equation}\label{u4}
\left\{
\begin{array}{ll}
(1+A)\pa_y^2\hat{u}_4-(1-\hat{s})(1-2\hat{s})\hat{u}_4=0,\\
\hat{u}_4(y)\sim y,\;\;y\to+\infty.
\end{array}
\right.
\end{equation}
 \begin{lemma}\label{u1solve}
Each system in \eqref{u1}-\eqref{u4} has a unique smooth solution, as follows: \begin{equation}\label{hatu1solution}
\hat{u}_1(y)=e^{ \f{y}{\sqrt{1+A}}}\Big(1+O(e^{y})\Big),
\end{equation}
and
\begin{equation}
\begin{split}\label{hatu2-4solution}
\hat{u}_2(y)&=e^{ -\f{y}{\sqrt{1+A}}}\Big(1+O(e^{y})\Big),\\
\hat{u}_3(y)&=1+O(e^{-y}),\\
\hat{u}_4(y)&=y\Big(1+O(e^{-y})\Big).
\end{split}
\end{equation}
\end{lemma}
$W(u,v)=u'v-v'u$ denotes the Wronskian of the second order ODE.
 \begin{lemma}\label{omega}
It follows that $W(\hat{u}_1,\hat{u}_2),\ W(\hat{u}_3,\hat{u}_4),\ W(\hat{u}_1,\hat{u}_3)\neq 0$.
\end{lemma}
Recall the definition
\begin{align*}S(\gamma,k)=\{&f\in C^{\infty}(\mathbb{R}): f^{\pm}=\lim_{y\to\pm\infty}f(y)\;\text{exist},\\
&\text{and}\;
\text{for}\;
 \forall j\geq 0, \exists k\geq 0,\;s.t. \;|\pa_z^j(f(y)-f^{\pm})|\lesssim |y|^ke^{-\gamma |y|} \}.\end{align*}
 The solvability of the equation \eqref{modele} is stated  as follows.
\begin{lemma}\label{le:mmodelsolve}
Let $f(y)\in S(1,k)$  with $f^{\pm}=0$.
 Then \eqref{modele} has a unique bounded solution if and only if
 \begin{equation}\label{solvemodel}
u^-=0,\;\;\;\text{and}\;\;\;
\Big((1+A)W(\hat{u}_1,  \hat{u}_3)\Big)^{-1}\int_{\mathbb{R}}\hat{u}_1(y)f(y)dy=u^+,
\end{equation}
where  $\hat{u}_1$  and $\hat{u}_3$ are as in  \eqref{u1} and \eqref{u3}. Moreover, the solution takes the form
\begin{align}\label{lemmodelsolve}
u(y)=&\Big((1+A)W(\hat{u}_1,  \hat{u}_3)\Big)^{-1}\Big(\int_{-\infty}^{y}
\hat{u}_3(y)\hat{u}_1(z)f(z)dz
+\int_{y}^{+\infty}\hat{u}_1(y)\hat{u}_3(z)f(z)dz\Big),
\end{align}
satisfying $u(y)\in S(1,k+1)$ with $u(+\infty)=u^{+}$.
\end{lemma}
\begin{proof}
It is easy to verify that
\begin{align}\label{ustar}
u_*(y)=&\Big((1+A)W(\hat{u}_1,  \hat{u}_3)\Big)^{-1}\Big(
\hat{u}_3(y)\int_{-\infty}^{y}\hat{u}_1(z)f(z)dz +\hat{u}_1(y)\int_{y}^{+\infty}\hat{u}_3(z)f(z)dz\Big)
\end{align}
is a $C^{\infty } (\mathbb{R})$ special solution of the equation
\begin{align}\label{unique}
-(1+A)\pa_y^2u+(1-\hat{s})(1-2\hat{s})u=f.
\end{align}
It follows directly by \eqref{ustar}, $f\in S(1,k)$, $\hat{u}_1,\hat{u}_3\in S(-\f{1}{\sqrt{A+1}},1)$ ($A>0$) that
 \begin{align}\label{limit-upm}
 \begin{split}
u_*(+\infty)&
=\Big((1+A)W(\hat{u}_1,  \hat{u}_3)\Big)^{-1}\int_{\mathbb{R}}\hat{u}_1(y)f(y)dy<\infty
,\\
u_*(-\infty)&=\lim_{z\rightarrow -\infty}e^{\f{ z}{\sqrt{1+A}}}\Big((1+A)W(\hat{u}_1,  \hat{u}_3)\Big)^{-1} \int_{\mathbb{R}}\hat{u}_3(y)f(y)dy=0,
\end{split}
\end{align}
which implies that  $u_*(y)$ is a bounded function.
Indeed, we claim that $u_*(y)$ is the unique bounded  solution of the equation \eqref{unique}. Therefore,  if \eqref{solvemodel} holds, then $u_*(y)$ exactly solves \eqref{modele}; if \eqref{solvemodel} does not hold, then the only possible bounded solution $u_*(y)$ fails to solve \eqref{modele}.

It remains to prove the uniqueness of the  bounded  solution of the equation \eqref{unique}. It suffices to show that the bounded solution for the homogeneous solution \eqref{unique} ($f\equiv 0$) must be zero. Indeed, if  $u(y)$ is a non-zero bounded  solution, then it must be the linear combination of $\hat{u}_1(y)$  and $\hat{u}_2(y)$:
\begin{align}\label{def:u-homo}
u(y)=C_1\hat{u}_1(y)+C_2\hat{u}_2(y).
\end{align}
We show $C_1=C_2=0$. From $\lim_{z\rightarrow -\infty}\hat{u}_1(z)=0 \text{ and } \lim_{y\rightarrow-\infty}\hat{u}_2(y)=+\infty$, we have $C_2\equiv 0$. Using Lemma \ref{omega}, there exist nonzero constants $C_3,\ C_4$ such that $\hat{u}_1(y)=C_3\hat{u}_3(y)+C_4\hat{u}_4(y),$
which together with \eqref{u3} and \eqref{u4} gives that \begin{align}\label{limitu1z}
\hat{u}_1(y)\sim C_3+C_4y\to+\infty,\;\;y\to+\infty.
\end{align}
Therefore, we have $C_1\equiv 0$. The proof of the uniqueness is finished.

Now we show $u_*(y)\in S(1,k+1)$, i.e., the exponential decay of $u_*(y)-u^{\pm}$. It follows by \eqref{solvemodel} and the condition $(1+A)W(\hat{u}_1,  \hat{u}_3)
u^+=\int_{\mathbb{R}}\hat{u}_1(y)f(y)dy$  that
\begin{align*}
(1+A)W(\hat{u}_1,  \hat{u}_3)&\big(u_*(y)-u^+(x,  t)\big)\\
=&\big(\hat{u}_3(y)-1\big)\int_{\mathbb{R}}\hat{u}_1(z)f(z)dz
-\hat{u}_3(y)\int_{y}^{+\infty}\hat{u}_1(z)f(z)dz
+\hat{u}_1(y)\int_{y}^{+\infty}\hat{u}_3(z)f(z)dz.
\end{align*}
 Then using $\hat{u}_3(y)-1=O(e^{-y})$(by \eqref{hatu2-4solution}), $|\hat{u}_1(y)|\lesssim y$(by \eqref{limitu1z}) and $ f\in S(1,k)$, that
 $$|u_*(y)- u^+|\lesssim
 e^{-y}+\int_{y}^{+\infty}zz^ke^{-z}dz+y\int_{y}^{+\infty}
z^ke^{-z}dz\lesssim y^{k+1}e^{-y},\;\;y\rightarrow +\infty.$$
On another hand, we write by the condition $u^-=\int_{\mathbb{R}}\hat{u}_3(y)f(y)dy=0$ that
\begin{align*}
&\quad(1+A)W(\hat{u}_1, \hat{u}_3)u_*(y)\\
=&\big(\hat{u}_1(y)
-e^{\f{y}{\sqrt{1+A}}}\big)\int_{y}^{+\infty}\hat{u}_3(z)f(z)dz
-e^{\f{y}{\sqrt{1+A}}}\int_{-\infty}^{y}\hat{u}_3(z)f(z)dz
+\hat{u}_3(y)\int_{-\infty}^y\hat{u}_1(z)f(z)dz.
\end{align*}
 Then using $\hat{u}_1(y)-e^{\f{y}{\sqrt{1+A}}}=O(e^{\f{y}{\sqrt{1+A}}+y})$ (by \eqref{u1solve}), $|\hat{u}_3(y)|\lesssim e^{-\f{y}{\sqrt{1+A}}}$ and $ f\in S(1,k)$, we get that
 $$|u_*(y)|\lesssim
e^{\f{y}{\sqrt{1+A}}+y}
+e^{\f{y}{\sqrt{1+A}}}\int_{-\infty}^yz^ke^{-\f{z}{\sqrt{1+A}}}e^{z}dz
+e^{-\f{y}{\sqrt{1+A}}}\int^{y}_{-\infty}
e^{\f{z}{\sqrt{1+A}}}z^ke^{-z}dz\lesssim y^{k+1}e^y,\;\;y\rightarrow -\infty.$$
The exponential decay rate for the derivatives follows in a similar manner.
\end{proof}
 \begin{proof}[The proof of Lemma \ref{u1solve}]
We first  need a lemma from  \cite[Lemma 2.4]{WAW}.
\begin{lemma}\label{lem:vol}
Let $a\in\mathbb{R}$, $g\in L^\infty(\pm\infty,a)$ and
$$\mu:=\Big|\int_{\pm\infty}^{a}\sup_{\tilde{y}\in[z,a]}
|K(z,\tilde{y})|dz\Big|<+\infty.$$
Then  the integral equation
$$u(y)=g(y)+\int_{\pm\infty}^{y}K(z,y)u(z)dz$$
admits a unique solution with the bound
$$\|u\|_{L^\infty(\pm \infty,a)}\leq e^{\mu} \|g\|_{L^\infty(\pm\infty,a)}.$$
\end{lemma}
Then we give the existence and behavior of $\hat{u}_i$ for $i=1,\cdots,4$. According to
$$(1-\hat{s})(1-2\hat{s})-1=O(e^{-|y|}) \text{ as } y\rightarrow -\infty,$$
we derive that the equation \eqref{u1} is equivalent  to
$$(1+A)\pa_y^2\hat{u}_1-\hat{u}_1=O(e^{-|y|})\hat{u}_1.$$
Then we have
\begin{equation}\label{hatu1}
\hat{u}_1(y)=e^{\f{y}{\sqrt{1+A}}}+\frac{1}{2\sqrt{1+A}}
\int_{-\infty}^{y}\Big(e^{\f{(y-z)}{\sqrt{1+A}}}
-e^{\f{(z-y)}{\sqrt{1+A}}}\Big)O(e^{-|z|})\hat{u}_1(z)dz.
\end{equation}
Let $\tilde{u}_1(y)=\hat{u}_1(y)e^{-\f{y}{\sqrt{1+A}}}$ and
\begin{equation}\label{kyz}
K(z,y)=\frac{1}{2\sqrt{1+A}}
\left(1-e^{\f{2(z-y)}{\sqrt{1+A}}}\right)O(e^{-|z|}).
\end{equation}
Then by \eqref{hatu1}, we have
\begin{equation*}
\tilde{u}_1(y)=1+\int_{-\infty}^{y}K(z,y)\tilde{u}_1(z)dz.
\end{equation*}
By \eqref{kyz}, for any $a\in\mathbb{R}$, we have
\begin{equation*}
\int_{-\infty}^{a}\sup_{z<\tilde{y}<a}|K(z,\tilde{y})|dz\leq C\max\{e^{a},1\}\leq C.
\end{equation*}
We deduce from Lemma \ref{lem:vol} that
$$\tilde{u}_1(y)=1+\int_{-\infty}^{y}K(z,y)\tilde{u}_1(z)dz$$
has a unique solution $\tilde{u}_1(y)$ with the bound $\|\tilde{u}_1\|_{L^\infty(- \infty,a)}\leq e^{C} \|1\|_{L^\infty(-\infty,a)}\leq C$ for  $a\in \mathbb{R}$ and satisfies $\lim_{y\rightarrow -\infty}\tilde{u}_1(y)=1$. Using \eqref{kyz} and  a more precise conclusion (cf. Lemma 4.1 in \cite{MRZ}), we can prove the asymptotic expansion as $y\to-\infty$:
\begin{align*}
\tilde{u}_1(y)=1+O(e^y).
\end{align*}
Thus we conclude that $\hat{u}_1(y)=e^{\f{y}{\sqrt{1+A}}}\tilde{u}_1(y)$ satisfies \eqref{u1}. The existence and behavior \eqref{hatu2-4solution} of $\hat{u}_2,\ \hat{u}_3$ and $\hat{u}_4$ follows analogously. We only list the integral equations:
\begin{align}
\begin{split}
\hat{u}_2(y)&=e^{-\f{y}{\sqrt{1+A}}} +\frac{1}{2\sqrt{1+A}}\int_{-\infty}^{y}\Big(e^{\f{(y-z)}{\sqrt{1+A}}}
-e^{\f{(z-y)}{\sqrt{1+A}}}\Big)O(e^{-|z|})\hat{u}_2(z)dz,\label{u2eq}\\
\hat{u}_3(y)&=1-\int_{y}^{+\infty}(y-z)O(e^{-|z|})\hat{u}_3(z)dz,\\
\hat{u}_4(y)&=y-\int_{y}^{+\infty}(y-z)O(e^{-|z|})\hat{u}_4(z)dz.
\end{split}
\end{align}
\end{proof}
 \begin{proof}[Proof of Lemma \ref{omega}]
$e^{\pm\f{y}{\sqrt{1+A}}}$ are the fundamental solutions of
$(1+A)\pa_y^2\hat{u}-\hat{u}=0,$
which gives \begin{align*}
W(\hat{u}_1,\hat{u}_2)&=\lim_{y\rightarrow -\infty}(e ^{\f{ y}{\sqrt{1+A}}})'e^{-\f{ y}{\sqrt{1+A}}}-e^{\f{ y}{\sqrt{1+A}}}(e^{-\f{ y}{\sqrt{1+A}}})'=\f{2}{\sqrt{1+A}} \neq 0.
\end{align*}
 $\hat{u}=y$ and $\hat{u}=1$ are the fundamental solutions of
$(1+A)\pa_y^2\hat{u}=0,$ which gives $$W(\hat{u}_3,\hat{u}_4)=-1\neq 0.$$
 $W(\hat{u}_1,\hat{u}_3)\neq 0$ follows from a contradiction argument.
If we have $W(\hat{u}_1,\hat{u}_3)= 0$, then there exists $C\neq 0$ such that $\hat{u}_1(y)=C\hat{u}_3(y)$, which gives a simple system for $\hat{u}_1$:
\begin{equation}\label{u4s}
\begin{cases}
\begin{array}{ll}
(1+A)\pa_y^2\hat{u}_1-(1-\hat{s})(1-2\hat{s})\hat{u}_1=0,\\
\lim_{y\rightarrow -\infty} \hat{u}_1(y)=e^{\f{y}{\sqrt{1+A}}},\\
\lim_{y\rightarrow +\infty} \hat{u}_1(y)=C.
\end{array}
\end{cases}
\end{equation}
Using
\begin{align*}
(1+A)\pa_y^2\hat{u}_1-(1-\hat{s})(1-2\hat{s})\hat{u}_1=
A\pa_y^2\hat{u}_1+\hat{s}^{-1}\big(\hat{s}^2\big(\frac{\hat{u}_1}{\hat{s}}\big)'\big)',
\end{align*}
one has
\begin{align*}
0&=\int_{\mathbb{R}}A\pa_y^2\hat{u}_1\hat{u}_1+s^{-1}\big(\hat{s}^2
\big(\frac{\hat{u}_1}{\hat{s}}\big)'\big)'\hat{u}_1dy\\
&=\int_{\mathbb{R}}-A|\pa_y\hat{u}_1|^2dy-A\hat{u}_1'\hat{u}_1
\Big|_{-\infty}^{+\infty}-\int_{\mathbb{R}}\hat{s}^2\big(
\big(\frac{\hat{u}_1}{\hat{s}}\big)'\big)^2dy+\hat{s}^2
\big(\frac{\hat{u}_1}{\hat{s}}\big)'\frac{\hat{u}_1}{\hat{s}}\Big|_{-\infty}^{+\infty}\\
&=\int_{\mathbb{R}}-A|\pa_y\hat{u}_1|^2-\hat{s}^2\big
(\big(\frac{\hat{u}_1}{\hat{s}}\big)'\big)^2dy+(\hat{u}_1
\hat{u}_1'-\frac{\hat{u}_1^2\hat{s}'}{\hat{s}})\Big|_{-\infty}^{+\infty}.
\end{align*}
Then we obtain
$$0=\int_{\mathbb{R}}A|\pa_y\hat{u}_1|^2+\hat{s}^2\big(\big(\frac{\hat{u}_1}{\hat{s}}\big)'\big)^2dy.$$
We infer from $A>0$ that $\hat{u}_1(y)\equiv 0$, which contradicts the fact $\hat{u}_1(y)=C\hat{u}_3(y)$.
\end{proof}
\subsection{Two integral identities}
\begin{lemma}\label{lem:energys0}
  Let $s_0$ be the bounded solution to \eqref{ode:s0} given in Lemma \ref{lem1}, then we have:
   \begin{align}\label{energy:s0}
2(1+\f{2L}{3})\int_{\mathbb{R}}s'(z)s_0'(z)dz=\int_{\mathbb{R}}zs'(z)f_0(z)dz.
\end{align}
\end{lemma}
\begin{proof}
The equation \eqref{ode:s0} can be equivalently written as
\begin{align*}
  -(1+\f{2L}{3})\frac{1}{s'}\partial_z\Big((s')^2\partial_z\Big(\frac{s_0}{s'}\Big)\Big)=f_0(z).
\end{align*}
Therefore we get
\begin{align*}
\int_{\mathbb{R}}zs'(z)f_0(z)\ud z=&-(1+\f{2L}{3})\int_{\mathbb{R}}z\partial_z\Big((s')^2\partial_z\Big(\frac{s_0}{s'}\Big)\Big)\ud z\\
=&-(1+\f{2L}{3})\int_{\mathbb{R}}z\partial_z(s's_0'-s''s_0)\ud z\\
=&~(1+\f{2L}{3})\int_{\mathbb{R}}(s's_0'-s''s_0)\ud z\\
=&~2(1+\f{2L}{3})\int_{\mathbb{R}}s's_0'\ud z,
\end{align*}
where we used the fact that $\int_{\mathbb{R}}(s''s_0+s's_0')\ud z=\int_{\mathbb{R}}\partial_z(s's_0)\ud z=0$.
\end{proof}

\begin{lemma}\label{lem:energys1}
  Let $s_1$ be the bounded solution to \eqref{ode:s1} given in Lemma \ref{lem2}, then we have:
   \begin{align}\label{energy:s1}
\int_{\mathbb{R}}s'(z)s_1'(z)dz=-\frac{6}{L}\int_{\mathbb{R}}s(z)f_1(z)dz.
\end{align}
\end{lemma}
\begin{proof}
The equation \eqref{ode:s1} can be equivalently written as
\begin{align*}
\frac{L}{6}s_1'' - (1+\f{2L}{3})\frac{1}{s}\partial_z\Big(s^2\partial_z\Big(\frac{s_1}{s}\Big)\Big)=f_1(z).
\end{align*}
Therefore
\begin{align*}
\int_{\mathbb{R}}s(z)f_1(z)\ud z=&\frac{L}{6}\int_{\mathbb{R}}ss_1''\ud z-(1+\f{2L}{3})\int_{\mathbb{R}}\partial_z\Big(s^2\partial_z\Big(\frac{s_1}{s}\Big)\Big)\ud z
=-\frac{L}{6}\int_{\mathbb{R}}s's_1'\ud z.
\end{align*}
The proof is completed.
\end{proof}
\section{Tensor calculations via exploitation of essential structures}
\subsection{Preliminary identities}
We derive some relations which hold in $\Gamma(\delta)$.
Using the fact $$ (n+\nabla d_0)\cdot(n-\nabla d_0)=|n|^2-|\nabla d_0|^2\equiv 0,$$ we get
\begin{align*}
   (n+\nabla d_0)\cdot \bar{h}=\frac{1}{d_0}(n+\nabla d_0)\cdot(n-\nabla d_0)\equiv 0.
\end{align*}
Thus
\begin{align}\label{eq:n-h}
&  n\cdot \bar{h}=\frac12(n+\nabla d_0)\cdot \bar{h}+\frac12(n-\nabla d_0)\cdot \bar{h}=\frac12 d_0|\bar{h}|^2,\\
&   \nabla d_0\cdot \bar{h}=-\frac12d_0|\bar{h}|^2,\qquad
  n\cdot\nabla d_0=1+\frac{1}{2}d_0^2|\bar{h}|^2. \label{eq:n-d0}
\end{align}
In addition, it holds that
\begin{align}\label{eq:divn}
  \nabla\cdot n=&\Delta d_0+d_0\nabla\cdot \bar{h}+\bar{h}\cdot\nabla d_0=\Delta d_0+d_0\nabla\cdot \bar{h}-\frac12d_0|\bar{h}|^2,\\  \nonumber
  (n\cdot\nabla) n=&~n\cdot\nabla(\nabla d_0+d_0h_{0})=\nabla^2d_0\cdot n+n\cdot\nabla d_0 \bar{h}+d_0(n\cdot\nabla)h_{0}\\ \nonumber
  =&~d_0\nabla^2d_0\cdot \bar{h}+n\cdot\nabla d_0 h_{0}+d_0(n\cdot\nabla)h_{0}\\
  =&~\bar{h}+d_0\big(\nabla^2d_0\cdot \bar{h}+(n\cdot\nabla)h_{0}\big)+\frac12d_0^2|\bar{h}|^2\bar{h}. \label{eq:dn}
  \end{align}
In particular, we have on $\Gamma$ that
\begin{align}
   \nabla\cdot n=\Delta d_0,\quad  (n\cdot\nabla) n=\bar{h}. \label{eq:div-n-Gamma}
\end{align}
\subsection{Proof of Lemmas \ref{lem:F1E0-1}, \ref{lem:F1Ei} and \ref{lem:boudaryFk+1Ei}}\label{f1123proof}
Substituting  \eqref{q0} and \eqref{G0} into \eqref{f1}, we have
\begin{align}\nonumber
F_1=&(\partial_t d_0-\Delta d_0)E^0s'-2(\nabla d_0\cdot\nabla) E^0s'-\frac{L}{2}\mathcal{N}_1(E^0,\nabla d_0)s'\\
&-\frac{L}{2}\Big(\mathcal{N}_2(E^0,\nabla d_0,\nabla d_1)+\mathcal{N}_2(E^0,\nabla d_1,\nabla d_0)\Big)s''
+\frac{L}{2}g_0(d_1-z)s''\nonumber\\
\triangleq & \sum_{j=1}^5I_j.\label{F1:fm'}
\end{align}

\begin{proof}[Proof of Lemma \ref{lem:F1E0-1}]
 Using \eqref{orth:E}, we have
 \begin{align}\label{F1E0-1}
 I_1:E^0&~=\f23(\partial_t d_0-\Delta d_0)s',\\ \label{F1E0-2}
I_2:E^0&~=-2s'(\nabla d_0\cdot \nabla) E^0:E^0=-s'(\nabla d_0\cdot \nabla) |E^0|^2=0,
\end{align}
Using the identity \eqref{N1:QQ} and the fact that $(E^0)^2=\frac19(3nn+I)$, we get
\begin{align*}
I_3:E^0=&-\frac{Ls'}{2}\mathcal{N}_1(E^0,\nabla d_0):E^0\\
=&-\frac{Ls'}{9}\nabla\cdot((3nn+I)\cdot\nabla d_0)\\
=&\f{-L}{9}\nabla \cdot \Big(4n-3(n\cdot \bar{h})n d_0-\bar{h}d_0\Big)s'\\
=&\f{-L}{9}\nabla \cdot \Big(4\nabla d_0-3(n\cdot \bar{h})n d_0+3\bar{h}d_0\Big)s'\\ \mynumber{F1E0-3}
=&\f{-4Ls'}{9}\Delta d_0+\frac{Ls'}{3}\nabla\cdot\big((n\cdot \bar{h})n d_0-\bar{h}d_0\big).
\end{align*}
Furthermore, we have  by \eqref{eq:n-h} that
\begin{align*}
  \nabla\cdot\big((n\cdot \bar{h})n d_0-\bar{h}d_0\big)&=\nabla\cdot\big(\frac12 d_0^2|\bar{h}|^2n-\bar{h}d_0\big)\\
  &=-h\cdot\nabla d_0-d_0\nabla\cdot \bar{h}+d_0|\bar{h}|^2n\cdot\nabla d_0+\frac12d_0^2\nabla\cdot(|\bar{h}|^2n)\\ \mynumber{J3-2}
  &=d_0\big(\frac12|\bar{h}|^2-\nabla\cdot \bar{h}+|\bar{h}|^2n\cdot\nabla d_0\big)+\frac12d_0^2\nabla\cdot(|\bar{h}|^2n),
\end{align*}
which vanishes on $\Gamma$.

In addition, we have
\begin{align*}
I_4:E^0=&-2LE_{kj}^0\pa_kd_0\pa_id_1E_{ij}^0=-\frac{2L}{9} (3nn+I):\nabla d_0\nabla d_1 s''\\ \mynumber{F1E0-4}
=&-\frac{2L}{3} (n\cdot\nabla d_0)(n\cdot\nabla d_1) s''=-\frac{2L}{3} d_0 (n\cdot\nabla d_0)(\bar{h}\cdot\nabla d_1) s'',
\end{align*}
and
\begin{align*}
I_5:E^0=&~(d_1-z)\f{Ls''}{2}(g_0:E^0)\\
=&~\f{Ls''}{2}(d_1-z)\Big\{-\frac{1}{3}(n\bar{h}+\bar{h}n)-2nn(n\cdot \bar{h})+\frac{8}{9}(\bar{h}\cdot n) I\nonumber\\ \nonumber
 &+d_0\Big((n\cdot \bar{h})(n\bar{h}+\bar{h}n) -\frac{2}{3}(\bar{h}\cdot n)^2I-\frac{2}{3}\bar{h}\bar{h}+\frac{2}{9}|\bar{h}|^2I\Big)\Big\}:E^0\\
=&~\f{Ls''}{2}(d_1-z)\Big(-\frac{4}{9}(\bar{h}\cdot n)-\frac43(\bar{h}\cdot n)+\frac{4d_0}{3}(\bar{h}\cdot n)^2-\frac{2d_0}{3}(\bar{h}\cdot n)^2+\frac{2d_0}{9}|\bar{h}|^2\Big)\\ \mynumber{F1E0-5}
=&~\f{Ls''}{3}(d_1-z)d_0\Big(-|\bar{h}|^2+(\bar{h}\cdot n)^2\Big).
 \end{align*}
Apparently, $I_4:E^0=I_5:E^0=0$ on $\Gamma$.  Thus, Lemma \ref{lem:F1E0-1} is proved.
\end{proof}

\begin{proof}[Proof of Lemma \ref{lem:F1Ei}]
First, we have on $\Gamma$ that
\begin{align*}
   I_1:E^i&~=0\quad \text{ for }i=1,\cdots,4;\\
I_2:E^1&~=-4s'(n\cdot \nabla )n\cdot l=-4s'\bar{h}\cdot l,\quad I_2:E^2=-4s'\bar{h}\cdot m,\\
 I_2:E^3&~=0,\quad I_2:E^4=0.
\end{align*}
Using \eqref{N1:def}, we have
\begin{align*}
I_3:E=&-Ls'\Big\{\pa_{ik}d_0(n_kn_j-\frac{1}{3}\delta_{kj})+n_i\pa_k(n_kn_j)+n_k\pa_i(n_kn_j)\Big\}_\G E_{ij}.
\end{align*}
Then, by the fact that $\partial_{ij}d_0 n_i =\partial_j(|\nabla d_0|^2)=0$ on $\Gamma$, we get on $\Gamma$ that
\begin{align*}
I_3:E^1=&-Ls'\Big\{-\frac13\pa_{ij}d_0+\pa_kn_k(n_in_j)+n_i(n_k\partial_k)n_j+\pa_i n_j\Big\}(n_il_j+l_in_j)\\
=&-2Ls'(n\cdot \nabla )n\cdot l=-2Ls'\bar{h}\cdot l,\\
I_3:E^2=&-2Ls'\bar{h}\cdot m,\quad \ I_3:E^3|_\G=0,\quad \ I_3:E^4|_\G=0.
\end{align*}
By \eqref{N2:def}, we get on $\Gamma$ that
\begin{align*}
 &\mathcal{N}_2(E^0,\nabla d_0,\nabla d_1):E= 2E^0_{kj}\pa_kd_0\pa_id_1E_{ij}= \frac{4}{3}\pa_jd_0\pa_id_1E_{ij},\\
 &\mathcal{N}_2(E^0,\nabla d_1,\nabla d_0):E= 2E^0_{kj}\pa_kd_1\pa_id_0E_{ij}= -\frac{2}{3}\pa_jd_1\pa_id_0E_{ij}.
\end{align*}
Thus it holds on $\Gamma$ that
\begin{align*}
I_4:E^1&=-\frac{L}{3}s''\pa_jd_0\pa_id_1(n_il_j+l_in_j)=-\frac{L}{3}s''(l\cdot \nabla d_1),\\
I_4:E^2&=-\frac{L}{3}s''(m\cdot \nabla d_1),\quad \ I_4:E^3=0,\quad  I_4:E^4=0.
\end{align*}
Finally, by \eqref{g0e}, we have
\begin{align*}
&I_5:E^1=-\frac{L}{3}(d_1-z)s''(\bar{h}\cdot l),\quad I_5:E^2=-\frac{L}{3}(d_1-z)s''(\bar{h}\cdot m),\\
&I_5:E^3=0,\qquad \ I_5:E^4=0.
\end{align*}
Combining the above identities, we obtain \eqref{f1e1}--\eqref{f1e34}.
\end{proof}
\begin{proof}[Proof of Lemma \ref{lem:boudaryFk+1Ei}]
We observe by the definition \eqref{Fm} that
 \begin{align*}
 F_{k+1}&=-\f{L}{2}\Big(\mathcal{N}_2(\pa_z^2Q_0,\nabla d_{k+1},\nabla d_0)
+\mathcal{N}_2(\pa_z^2Q_0,\nabla d_{0},\nabla d_{k+1})\Big)+G_0d_{k+1}\\
&\quad+\text{terms}\;\text{depending}
\;\text{on}\;d_{p},\;s_{p,j}(p\leq k,0\leq j\leq 4).
\end{align*}
On $\G$, for $i=1,2$, we denote  $I_4=-\f{L}{2}\Big(\mathcal{N}_2(\pa_z^2Q_0,\nabla d_{k+1},\nabla d_0)
+\mathcal{N}_2(\pa_z^2Q_0,\nabla d_{0},\nabla d_{k+1})\Big)$ and $I_5=G_0d_{k+1}$ similarly as in the proof of Lemma \ref{lem:F1Ei}, and calculate  $I_4:E^i$ and $I_5:E^i$ . Using \eqref{N2:def}, we get on $\Gamma$ that
\begin{align*}
 &\mathcal{N}_2(E^0,\nabla d_0,\nabla d_{k+1}):E= 2E^0_{lj}\pa_ld_0\pa_id_{k+1}E_{ij}= \frac{4}{3}\pa_jd_0\pa_id_{k+1}E_{ij},\\
 &\mathcal{N}_2(E^0,\nabla d_{k+1},\nabla d_0):E= 2E^0_{lj}\pa_ld_{k+1}\pa_id_0E_{ij}.
\end{align*}
Then we have
\begin{align*}
I_4:E^1&=-\frac{L}{3}s''\pa_jd_0\pa_id_{k+1}(n_il_j+l_in_j)
=-\frac{L}{3}s''(l\cdot \nabla d_{k+1}),\quad I_4:E^2=-\frac{L}{3}s''(m\cdot \nabla d_{k+1}),
\end{align*}
 and
 \begin{align*}
&I_5:E^1=\frac{2L}{3}d_{k+1}s''(\bar{h}\cdot l),\quad I_5:E^2=\frac{2L}{3}d_{k+1}s''(\bar{h}\cdot m).
 \end{align*}
Therefore, we arrive at the first two lines in \eqref{fm:boudaryFk+1Ei}.

 For $i=0,3,4$, we have on $\G$ that $G_0:E^i=0$.
Then we have by \eqref{didj} that
\begin{align*}
    &I_4:E^0=-\frac{8}{9}Ln\cdot\nabla d_{k+1}=\frac{4}{9}L\sum_{j=1}^{k}\nabla d_j\cdot\nabla d_{k+1-j},\\
    &I_4:E^i=0,\text{ for }\;i=3,4.
\end{align*}
Thus we conclude that on $\Gamma$, $F_{k+1}:E^i$($i=0,3,4$) only depends on $d_p,s_{p,j}(0\leq p\leq k,\;0\leq j\leq 4)$, which is the last line in \eqref{fm:boudaryFk+1Ei}.
\end{proof}
\subsection{Proof of Lemma \ref{lem:d1}}\label{lem:d1proof}
For simplicity, we denote $\langle \phi_1,\phi_2 \rangle$ by $\int_{\mathbb{R}}\phi_1(z)\phi_2(z)\ud z$ for scalar functions $\phi_1, \phi_2$ on $\mathbb{R}$, and $\langle P, Q \rangle=\int_{\mathbb{R}}P_{ij}(z)Q_{ij}(z)\ud z$ for matrix-valued functions $P$ and $Q$.

First, we have
\begin{lemma}
It hold on $\Gamma$ that
\begin{align}\label{dF1E0} \nonumber
\pa_{\nu}(F_1:E^0)
&=\Big\{\f23\pa_{\nu}\big(\partial_t d_0-(1+\f{2L}{3})\Delta d_0\big)+\f{L}{3}\Big(\f32|\bar{h}|^2-(\nabla\cdot \bar{h})\Big)\Big\}s'\\
&\quad-\f{L}{3}\Big(2(\bar{h}\cdot \nabla d_1) +|\bar{h}|^2(d_1-z)\Big)s'',
\end{align}
and
\begin{align}\label{fm:g10}
 g_{1,0}
&=\f23\pa_{\nu}\Big[\partial_t d_0-(1+\f{2L}{3})\Delta d_0\Big]+\f{L}{3}\Big(|\bar{h}|^2-(\nabla\cdot \bar{h})\Big)\f{\langle s', s'\rangle}{\langle \eta', s'\rangle}.
\end{align}
Thus
\begin{align}\label{dF1-integral}
\int_{\mathbb{R}}zs'\pa_{\nu}(F_1:E^0)\ud z=&\f{L}{6}\Big(2(\bar{h}\cdot \nabla d_1) +|\bar{h}|^2d_1\Big)\langle s', s'\rangle.
\end{align}
\end{lemma}
\begin{proof}
Since $\partial_\nu(d_0 f)|_{\Gamma}=f_{\Gamma}$ for a $C^1$ function $f$ defined on $\Gamma(\delta)$, equation \eqref{dF1E0} is
a consequence of  identities \eqref{F1E0-1}-\eqref{F1E0-5}.
Then  due to the facts
\begin{align*}
  \langle s'', s'\rangle=0,\quad \text{ and } \langle zs'', s'\rangle=\langle z, \big(\frac12(s')^2\big)'\rangle=-\frac12\langle s', s'\rangle,
\end{align*}
 \eqref{fm:g10} follows.
In addition with
\begin{align*}
  \langle zs'', zs'\rangle=\langle z^2, \partial_z(\frac12(s')^2)\rangle=-\langle z, (s')^2\rangle =0,
\end{align*}
we obtain \eqref{dF1-integral}.
\end{proof}
\begin{proof}[Proof of Lemma \ref{lem:d1} ]
By the definition of $F_2$ and $F_1$, we have:
\begin{align*}
\langle F_2, ~~\pa_zQ_0 \rangle =&
\langle \pa_zQ_0 , ~\pa_zQ_0 \rangle(\pa_t-\Delta )d_1 +\langle \pa_zQ_1, ~\pa_zQ_0 \rangle (\pa_t-\Delta )d_0\\
&-2\langle(\nabla d_1\cdot \nabla)\pa_zQ_0, ~\pa_zQ_0\rangle
-2\langle (\nabla d_0\cdot \nabla)\pa_zQ_1, ~\pa_zQ_0\rangle\\
&-\frac{L}{2}\langle\mathcal{N}_1(\pa_zQ_1,\nabla d_0),~\pa_zQ_0\rangle
-\frac{L}{2}\langle\mathcal{N}_1(\pa_zQ_0,\nabla d_1),~\pa_zQ_0\rangle\\
&-\frac{L}{2}\sum_{\substack{q+r=2}}
\langle \mathcal{N}_2(\pa_z^2Q_0,\nabla d_q,\nabla d_r),~\pa_zQ_0\rangle\\
&-\frac{L}{2}\langle
\mathcal{N}_2(\pa_z^2Q_1,\nabla d_0,\nabla d_1)+\mathcal{N}_2(\pa_z^2Q_1,\nabla d_1,\nabla d_0),~\pa_zQ_0\rangle\\
&+\langle (\pa_t -\mathcal{L})Q_{0},~\pa_zQ_0\rangle-\langle B_1+C_1,~\pa_zQ_0\rangle\\
&+d_2\langle G_0,~ \pa_zQ_0\rangle
+\langle (d_1-z)G_1,~\pa_zQ_0\rangle,
\end{align*}
and
\begin{align*}
\langle F_1,~\pa_zQ_1\rangle=&~
\langle \pa_zQ_0,~\pa_zQ_1 \rangle(\pa_t-\Delta )d_0-2\langle \nabla d_0\cdot\nabla \partial_zQ_0,~\pa_zQ_1\rangle\\
&-\frac{L}{2}\langle \mathcal{N}_1(\pa_zQ_0,\nabla d_0),~\pa_zQ_1\rangle+\langle G_0(d_1-z),~\pa_zQ_1\rangle\\
&-\frac{L}{2}\langle\mathcal{N}_2(\pa_z^2Q_0,\nabla d_0,\nabla d_1)+\mathcal{N}_2(\pa_z^2Q_0,\nabla d_1,\nabla d_0),~\pa_zQ_1\rangle.
\end{align*}
By the equality \eqref{N2:PQ}, we have
\begin{align*}
&  \langle\mathcal{N}_2(\pa_z^2Q_1,\nabla d_0,\nabla d_1)+\mathcal{N}_2(\pa_z^2Q_1,\nabla d_1,\nabla d_0),~\pa_zQ_0\rangle\\
&\qquad+\langle\mathcal{N}_2(\pa_z^2Q_0,\nabla d_0,\nabla d_1)+\mathcal{N}_2(\pa_z^2Q_0,\nabla d_1,\nabla d_0),~\pa_zQ_1\rangle\\
&=\int_{\mathbb{R}}\partial_z\Big(\big(\mathcal{N}_2(\pa_zQ_1,\nabla d_0,\nabla d_1)+\mathcal{N}_2(\pa_zQ_1,\nabla d_1,\nabla d_0)\big):\pa_zQ_0\Big)\ud z=0,
\end{align*}
and
\begin{align*}
& \langle \mathcal{N}_2(\pa_z^2Q_0,\nabla d_q,\nabla d_r)+ \mathcal{N}_2(\pa_z^2Q_0,\nabla d_r,\nabla d_q),~\pa_zQ_0\rangle\\
&=\frac12\int_{\mathbb{R}}\partial_z\Big(\big(\mathcal{N}_2(\pa_zQ_0,\nabla d_q,\nabla d_r)+\mathcal{N}_2(\pa_zQ_0,\nabla d_r,\nabla d_q)\big):\pa_zQ_0\Big)\ud z=0.
\end{align*}

Therefore, it holds on $\Gamma$ that
\begin{align*}
  \langle F_2, &~~\pa_zQ_0 \rangle +\langle F_1,~\pa_zQ_1\rangle\\
=&
\langle \pa_zQ_0 , ~\pa_zQ_0 \rangle(\pa_t-\Delta )d_1
-2\langle(\nabla d_1\cdot \nabla)\pa_zQ_0, ~\pa_zQ_0\rangle\\
&+2\langle \pa_zQ_1, ~\pa_zQ_0 \rangle (\pa_t-\Delta )d_0
-2(\nabla d_0\cdot \nabla)\langle \pa_zQ_1, ~\pa_zQ_0\rangle\\
&-\frac{L}{2}\Big(\langle\mathcal{N}_1(\pa_zQ_1,\nabla d_0),~\pa_zQ_0\rangle
+\langle\mathcal{N}_1(\pa_zQ_0,\nabla d_1),~\pa_zQ_0\rangle+\langle \mathcal{N}_1(\pa_zQ_0,\nabla d_0),~\pa_zQ_1\rangle\Big)\\
&+\langle (\pa_t -\mathcal{L})Q_{0},~\pa_zQ_0\rangle-\langle B_1+C_1,~\pa_zQ_0\rangle\\
&+d_2\langle G_0,~ \pa_zQ_0\rangle
+\langle (d_1-z)G_1,~\pa_zQ_0\rangle+\langle G_0(d_1-z),~\pa_zQ_1\rangle\\
\triangleq & J_1+J_2+\cdots +J_5. \mynumber{sum:F2F1}
\end{align*}
{\it Calculation of $J_1$:}  As $\langle \pa_zQ_0 , ~\pa_zQ_0 \rangle= \frac23\langle s',~s'\rangle$ and $-\nabla d_1\cdot \nabla(|E_0|^2)\langle s',s'\rangle=0$,  we have
\begin{align*}
  J_1=\frac23\langle s',~s'\rangle (\pa_t-\Delta )d_1.
\end{align*}
{\it Calculation of $J_2$:}  Since
\begin{align}\nonumber
&\int_{\mathbb{R}} (\pa_zQ_1\pa_zQ_0+\pa_zQ_0\pa_zQ_1)\ud z\\
&=\f29(I+3nn) \langle s',s_{1,0}'\rangle
+\f13(ln+nl)\langle s',s_{1,1}'\rangle
+\f13(mn+nm)\langle s',s_{1,2}'\rangle, \label{eq:Q1Q0}
\end{align}
we have
\begin{align*}
  \langle \pa_zQ_1, ~\pa_zQ_0 \rangle~&=\frac12 \tr \int_{\mathbb{R}} (\pa_zQ_1\pa_zQ_0+\pa_zQ_0\pa_zQ_1)\ud z=\frac23\langle s',s_{1,0}'\rangle.
\end{align*}
On $\Gamma$, $s_{1,0}$ and  $\partial_\nu s_{1,0}$ are the solutions to \eqref{ode:s0} with $f_0=-\frac32F_1:E^0=0$ and $f_0=-\frac32\partial_\nu(F_1:E^0)-\eta' g_{1,0}$ respectively. Thus, we get by Lemma \ref{lem:energys0} that on $\Gamma$
\begin{align}\label{eq:s10}
  \langle s',s_{1,0}'\rangle=&0,\\ \nonumber \label{eq:ds10}
  \partial_\nu\langle s',s_{1,0}'\rangle=&-\frac{1}{2(1+\frac{2L}{3})}\int_{\mathbb{R}}zs'(z)(\frac32\partial_\nu(F_1:E^0)+\eta' g_{1,0})dz\\
  =&-\frac{1}{2(1+\frac{2L}{3})}\frac{3}{2}\frac{L}{6}\Big(2(\bar{h}\cdot \nabla d_1) +|\bar{h}|^2d_1\Big)\langle s', s'\rangle,
\end{align}
where we have used \eqref{dF1-integral} and the fact that $\int_{\mathbb{R}}zs'\eta'\ud z=0$ (as $s'$, $\eta'$ are even functions on $\mathbb{R}$).
Therefore, we have on $\Gamma$ that
\begin{align*}
  \langle \pa_zQ_1, ~\pa_zQ_0 \rangle~ &=-\frac3{2(3+2L)}\int_{\mathbb{R}}zs'(z)F_1:E^0\ud z=0,
  \end{align*}
  and
  \begin{align*}
J_2&=2\nabla d_0\cdot\nabla \langle \pa_zQ_1, ~\pa_zQ_0 \rangle=\frac43\partial_\nu\langle s',s_{1,0}'\rangle\\
 &=-\frac{L}{3+{2L}}\Big((\bar{h}\cdot \nabla d_1) +\f12|\bar{h}|^2d_1\Big)\langle s', s'\rangle.
\end{align*}

{\it Calculation of $J_3$:} By \eqref{N1:PQ}, we have
\begin{align} \nonumber
&\langle \mathcal{N}_1(\pa_zQ_1,\nabla d_0), ~\pa_zQ_0\rangle+\langle \mathcal{N}_1(\pa_zQ_0,\nabla d_0),~\pa_zQ_1\rangle\\
=&2\nabla \cdot\Big(\int_{\mathbb{R}}\Big(\langle \pa_zQ_1\pa_zQ_0\rangle +\langle\pa_zQ_0\pa_zQ_1\rangle\Big)\ud z\cdot \nabla d_0\Big)\notag\\
=&\frac23\nabla \cdot\Big(\Big(\f23(I+3nn)\langle s',s_{1,0}'\rangle+(ln+nl)\langle s',s_{1,1}'\rangle
+(mn+nm)\langle s',s_{1,2}'\rangle\Big)\cdot(n-d_0\bar{h})\Big).
\end{align}
Direct calculations give us that
\begin{align*}
&\nabla \cdot\Big((I+3nn)\langle s',s_{1,0}'\rangle\cdot (n-d_0\bar{h})\Big)\\
&=\nabla \cdot\Big(\big(4n-d_0 \bar{h}-3d_0n(\bar{h}\cdot n)\big)\langle s',s_{1,0}'\rangle\Big)\\
&=4n\cdot\nabla \langle s',s_{1,0}'\rangle+(4\nabla\cdot n-\bar{h}\cdot \nabla d_0-3(n\cdot\nabla d_0)(\bar{h}\cdot n))\langle s',s_{1,0}'\rangle\\
&=4\partial_\nu\langle s',s_{1,0}'\rangle+4\nabla\cdot n\langle s',s_{1,0}'\rangle.
\end{align*}
Using the identities \eqref{eq:s10}-\eqref{eq:ds10}, we have on $\Gamma$ that
\begin{align*}
&\nabla \cdot\Big((I+3nn)\langle s',s_{1,0}'\rangle\cdot (n-d_0\bar{h})\Big)=-\frac{3L}{2(3+{2L})}\Big(2(\bar{h}\cdot \nabla d_1) +|\bar{h}|^2d_1\Big)\langle s', s'\rangle.
\end{align*}
Furthermore, we can directly get
\begin{align*}
&\nabla \cdot\Big((ln+nl)\langle s',s_{1,1}'\rangle\cdot (n-d_0\bar{h})\Big)\\
&=\nabla \cdot\Big(\big(l-d_0l(\bar{h}\cdot n)-d_0n(\bar{h}\cdot l)\big)\langle s',s_{1,1}'\rangle\Big)\\
&= l\cdot\nabla\langle s',s_{1,1}'\rangle+
\Big(\nabla \cdot l-(l\cdot\nabla d_0)(\bar{h}\cdot n)-(n\cdot\nabla d_0)(\bar{h}\cdot l)\Big)\langle s',s_{1,1}'\rangle \\
&= l\cdot\nabla\langle s',s_{1,1}'\rangle+
\big(\nabla \cdot l-(\bar{h}\cdot l)\big)\langle s',s_{1,1}'\rangle.
\end{align*}
As $l\cdot \nabla $ is tangential derivative on $\Gamma$, it suffices to calculate $\langle s',s_{1,1}'\rangle$ on the interface $\Gamma$.
Recalling that $s_{1,1}$ on $\Gamma$ satisfies the equation \eqref{ode:s1} with $f_1=-\frac{1}{2}F_1:E^1$, we get by Lemma \ref{lem:energys1} that:
\begin{align*}
\langle s',s_{1,1}'\rangle|_{\Gamma}=&\frac{3}{L}\int_{\mathbb{R}}s(z)F_1:E^1 \ud z\\
=& -\frac{3}{L}\int_{\mathbb{R}}s\Big\{(4+2L)s'(\bar{h}\cdot l)
   +\frac{L}{3}s''(l\cdot \nabla d_1)+\frac{L}{3}(d_1-z)s''(\bar{h}\cdot l)\Big\}\ud z\\
=&-\frac{12+7L}{2L}(\bar{h}\cdot l)+\big(l\cdot \nabla d_1+d_1(\bar{h}\cdot l)\big)\langle s', s'\rangle.
\end{align*}
Here we used the fact that $-\langle zs'', s \rangle= \langle zs', s' \rangle+\langle s', s \rangle=\langle s', s \rangle=\frac12$.
Similarly,
\begin{align*}
\langle s',s_{1,2}'\rangle|_{\Gamma}=-\frac{12+7L}{2L}(\bar{h}\cdot m)+\big(m\cdot \nabla d_1+d_1(\bar{h}\cdot m)\big)\langle s', s'\rangle.
\end{align*}
Therefore, using the relation $ll+mm=I-nn$, we get on $\G$ that
\begin{align*}
  l\langle s',s_{1,1}'\rangle+m\langle s',s_{1,2}'\rangle=&-\frac{12+7L}{2L}l(\bar{h}\cdot l)+l\big(l\cdot \nabla d_1+d_1(\bar{h}\cdot l)\big)\langle s', s'\rangle\\
  &-\frac{12+7L}{2L}m(\bar{h}\cdot m)+m\big(m\cdot \nabla d_1+d_1(\bar{h}\cdot m)\big)\langle s', s'\rangle\\
  =&-\frac{12+7L}{2L}(I-nn)\bar{h}+(I-nn)\big(\nabla d_1+d_1\bar{h}\big)\langle s', s'\rangle\\
  =&-\frac{12+7L}{2L}\bar{h}+\big(\nabla d_1+d_1\bar{h}\big)\langle s', s'\rangle, \mynumber{combine:s12}
\end{align*}
which yields that
\begin{align*}
  \nabla \cdot\Big(&(ln+nl)\langle s',s_{1,1}'\rangle\cdot (n-d_0\bar{h})\Big)
  +\nabla \cdot\Big((mn+nm)\langle s',s_{1,2}'\rangle\cdot (n-d_0\bar{h})\Big)\\
  =&\nabla\cdot\big(l\langle s',s_{1,1}'\rangle+m\langle s',s_{1,2}'\rangle\big)
  -\bar{h}\cdot\big(l\langle s',s_{1,1}'\rangle+m\langle s',s_{1,2}'\rangle\big)\\
=&-\frac{12+7L}{2L}\nabla\cdot\bar{h}+\big(\Delta d_1+d_1\nabla\cdot\bar{h}+\bar{h}\cdot\nabla d_1\big)\langle s', s'\rangle\\
&+\frac{12+7L}{2L}|\bar{h}|^2-\big(\bar{h}\cdot\nabla d_1+d_1|\bar{h}|^2\big)\langle s', s'\rangle\\
=&\frac{12+7L}{2L}(|\bar{h}|^2-\nabla\cdot\bar{h})+\big(\Delta d_1+d_1(\nabla\cdot\bar{h}-|\bar{h}|^2)\big)\langle s', s'\rangle.
\end{align*}
Thus
\begin{align*} 
&\langle \mathcal{N}_1(\pa_zQ_1,\nabla d_0),~\pa_zQ_0\rangle+\langle \mathcal{N}_1(\pa_zQ_0,\nabla d_0),~\pa_zQ_1\rangle\\
=&\frac23\nabla \cdot\Big(\Big(\f23(I+3nn)\langle s',s_{1,0}'\rangle+(ln+nl)\langle s',s_{1,1}'\rangle
+(mn+nm)\langle s',s_{1,2}'\rangle\Big)\cdot(n-d_0\bar{h})\Big)\\
=&-\frac{2L}{3(3+{2L})}\Big(2(\bar{h}\cdot \nabla d_1) +|\bar{h}|^2d_1\Big)\langle s', s'\rangle\\
&+\frac{12+7L}{3L}(|\bar{h}|^2-\nabla\cdot\bar{h})+\frac23\big(\Delta d_1+d_1(\nabla\cdot\bar{h}-|\bar{h}|^2)\big)\langle s', s'\rangle.
\end{align*}

By \eqref{N1:QQ}, we get
\begin{align*}
\langle \mathcal{N}_1(\pa_zQ_0,\nabla d_1),~\pa_zQ_0 \rangle
=&\mathcal{N}_1(E^0,\nabla d_1):E^0\langle s',~s' \rangle\\
=&\f{2}{9}\nabla\cdot\big((I+3nn)\cdot\nabla d_1\big)\langle s',~s' \rangle\\
=&\Big(\f{2}{9}\Delta d_1+\f{2}{3}\nabla\cdot\big(nn\cdot\nabla d_1\big)\Big)\langle s',~s' \rangle\\
=&\Big(\f{2}{9}\Delta d_1+\f{2}{3}(n\cdot\nabla )n\cdot\nabla d_1+\f23n_in_k\pa_{ik}d_1\Big)\langle s',~s' \rangle\\
=&\Big(\frac{2}{9}\Delta d_1+\f{2}{3} \bar{h} \cdot \nabla d_1\Big)\langle s',~s' \rangle.
\end{align*}
Combine the above two estimates, we get
\begin{align*}
J_3
=&\Big(-\frac{4L}{9}\Delta d_1
-\frac{L}{3+2L}\bar{h}\cdot \nabla d_1+
\big(-\f{L}{3}\nabla\cdot\bar{h}+\f{(1+L)L}{3+2L}|\bar{h}|^2\big)d_1
\Big)\langle s', s'\rangle\\
&-\f{12+7L}{6}(|\bar{h}|^2-\nabla\cdot\bar{h}).
\end{align*}

{\it Calculation of $J_5$:}
First, by the definition of $G_0$ \eqref{G0}, we have
\begin{align}
G_0(z,x,t)&=\frac{1}{d_0}(\mathcal{L}_\nu Q_0-\mathcal{L}_* Q_0)=\frac{L}{2}\Big(\partial_z^2Q_0 P_*+ P_*\partial_z^2 Q_0 -\frac{2I}{3} P_*:\partial_z^2Q_0\Big),
\end{align}
with
\begin{align*}
P_*=&~d_0^{-1}\big((\nabla d_0\nabla d_0)-nn\big)=-(n\bar{h}+\bar{h}n)+d_0\bar{h}\bar{h}.
\end{align*}
Thus,
\begin{align*}
d_2\langle G_0,~ \pa_z Q_0\rangle& = \frac{Ld_2}{2}\langle \partial_z^2Q_0 P_*+
 P_*\partial_z^2 Q_0 -\frac{2I}{3} P_*:\partial_z^2Q_0,~ \pa_z Q_0\rangle\\
 & = \frac{Ld_2}{2}\int_{\mathbb{R}} \partial_z\Big((\partial_zQ_0)^2: P_* \Big)\ud z=0.
\end{align*}

Recalling the definitions  $G_1$ \eqref{Gm}, we have
\begin{align}\label{G1}
G_1(z,x,t)&=\sum_{i=0}^{2}\eta_i'(z)g_{1,i}(x,t)E^i(x,t)+\frac{1}{d_0}(\mathcal{L}_\nu Q_1-\mathcal{L}_* Q_1)\\
&=\sum_{i=0}^{2}\eta_i'(z)g_{1,i}(x,t)E^i(x,t)+\frac{L}{2}\Big(\partial_z^2Q_1 P_*+ P_*\partial_z^2 Q_1 -\frac{2I}{3}P_*:\partial_z^2Q_1\Big).
\end{align}

Thus,
\begin{align*}
    J_{5}=&
    \langle(d_1-z)G_1,~ \pa_z Q_0\rangle+\langle (d_1-z)G_0,~\pa_zQ_1\rangle\\
    =&\langle(d_1-z)\eta'(z)g_{1,0}(x,t)E^0(x,t),\partial_zQ_0\rangle\\
    &+\frac{L}{2}\int_\mathbb{R}(d_1-z)\Big\{\Big(\partial_z^2Q_1 P_*+ P_*\partial_z^2 Q_1 -\frac{2I}{3}P_*:\partial_z^2Q_1\Big):\pa_z Q_0\\
    &+\Big(\partial_z^2Q_0 P_*+ P_*\partial_z^2 Q_0 -\frac{2I}{3} P_*:\partial_z^2Q_0\Big):\pa_zQ_1\Big\}\ud z\\
    =&\frac23g_{1,0}(x,t)\langle(d_1-z)\eta'(z), s'\rangle+L\int_\mathbb{R}(d_1-z)\partial_z\big(\partial_zQ_1 P_*:\partial_zQ_0\big)\ud z\\
    =&\frac23g_{1,0}(x,t)d_1\langle\eta'(z), s'\rangle+L\Big(\int_\mathbb{R} \partial_zQ_1\partial_zQ_0\ud z\Big): P_*\\
    =&-\f{L}3\bar{h}\cdot\nabla d_1 \langle s', s'\rangle-\f{L}{9}\Big(2(\nabla\cdot\bar{h})
    +|\bar{h}|^2\Big)d_1 \langle s', s'\rangle\\
    &\quad+\f49\langle\eta', s' \rangle\pa_{\nu}\Big[\partial_t d_0-(1+\f{2L}{3})\Delta d_0\Big] d_1+\frac{12+7L}{6}|\bar{h}|^2.
    \end{align*}
Indeed, it yields from \eqref{fm:g10} that
\begin{align*}
\frac23g_{1,0}(x,t)d_1\langle\eta', s' \rangle
  =\f49d_1\langle\eta', s' \rangle\pa_{\nu}\Big[\partial_t d_0-(1+\f{2L}{3})\Delta d_0\Big]+\f{2Ld_1}{9}\Big(|\bar{h}|^2-(\nabla\cdot \bar{h})\Big)\langle s', s'\rangle.
\end{align*}
By \eqref{eq:Q1Q0} and \eqref{combine:s12}, we have
\begin{align*}\nonumber
&2L\Big(\int_\mathbb{R} \partial_zQ_1\partial_zQ_0\ud z\Big): P_*\\
&=L\int_{\mathbb{R}} (\pa_zQ_1\pa_zQ_0+\pa_zQ_0\pa_zQ_1)\ud z:P_*\\
&=\f{2L}9(I+3nn):P_* \langle s',s_{1,0}'\rangle
+\f{L}3(ln+nl):P_*\langle s',s_{1,1}'\rangle \\
&\quad+\f{L}3(mn+nm):P_*\langle s',s_{1,2}'\rangle\nonumber \\
&=-\f{2L}3(\bar{h}\cdot l)\langle s',s_{1,1}'\rangle-\f23(\bar{h}\cdot m)\langle s',s_{1,2}'\rangle\\
&=-\f{2L}3 \bar{h}\cdot\Big(l\langle s',s_{1,1}'\rangle+m\langle s',s_{1,2}'\rangle\Big)\\
&=-\f{2L}3 \bar{h}\cdot\Big(-\frac{12+7L}{2L}\bar{h}+\big(\nabla d_1+d_1\bar{h}\big)\langle s', s'\rangle\Big)\\
&=\frac{12+7L}{3}|\bar{h}|^2-\f{2L}3\big(\bar{h}\cdot\nabla d_1+d_1|\bar{h}|^2\big)\langle s', s'\rangle.
\end{align*}

{\it Completion of the proof:} We have from \eqref{sum:F2F1} that
\begin{align}
\frac23 \langle s', s'\rangle& (\partial_td_1-\Delta d_1)+J_2+J_3+J_5\notag\\
 &=\langle F_1,~\pa_zQ_1\rangle+\langle B_1+C_1,~\pa_zQ_0\rangle-\langle (\pa_t -\mathcal{L})Q_{0},~\pa_zQ_0\rangle.\label{J_2-J5}
\end{align}
By the definition of $F_1$, $B_1$, $C_1$, we have
\begin{align*}
 & \langle F_1,~\pa_zQ_1\rangle+\langle B_1+C_1,~\pa_zQ_0\rangle\\
& = \langle  \mathcal{L}_{\nu}Q_{1}+\mathcal{H}_{Q_0}Q_{1}, ~\pa_zQ_1\rangle+\langle B_1+C_1,~\pa_zQ_0\rangle\\
& =\langle  \mathcal{L}_{\nu}Q_{1}, ~\pa_zQ_1\rangle+\langle\mathcal{H}_{Q_0}Q_{1}, \partial_zQ_1\rangle +\langle B_1+C_1,~\pa_zQ_0\rangle.
\end{align*}
From \eqref{op:L_nu}, we get
\begin{align*}
  \langle  \mathcal{L}_{\nu}Q_{1}, ~\pa_zQ_1\rangle=\frac12\int_{\mathbb{R}}\partial_z\Big(|\partial_zQ_1|^2+L(\partial_zQ_1)^2:(\nabla d_0\nabla d_0) \Big)\ud z=0.
\end{align*}
By the definition of $B_1$ and $C_1$, we have
\begin{align*}
  \langle\mathcal{H}_{Q_0}Q_{1}, \partial_zQ_1\rangle +\langle B_1+C_1,~\pa_zQ_0\rangle &=\int_{\mathbb{R}}\partial_z\Big({H}_{Q_0}Q_{1}:Q_1\Big)\ud z\\
 &= {H}_{Q_0^+}Q_{1}^+:Q_1^+ =0.
\end{align*}
Thus, combining all the above calculations in \eqref{J_2-J5},  we have  for some $\mu_2(x,t)$, $\mu_0(x,t)$ only depending on $d_0$ that
\begin{align}
&\quad\frac23 \langle s', s'\rangle \Big(\partial_td_1-(1+\f{2L}{3})\Delta d_1\Big)-\frac{9L+2L^2}{3(3+{2L})}\bar{h}\cdot \nabla d_1 \langle s', s'\rangle -\f23(\mu_{2}d_1+\mu_0)\langle s', s'\rangle=0,
\end{align}
which completes the proof.
\end{proof}

\subsection{Sketch of proof for Lemma \ref{lem:dk+1}}
Since the process is similar to the proof of Lemma \ref{lem:d1}, we only give a sketched proof. We deduce that
\begin{align*}
&  \langle\mathcal{N}_2(\pa_z^2Q_{k+1},\nabla d_0,\nabla d_1)+\mathcal{N}_2(\pa_z^2Q_{k+1},\nabla d_1,\nabla d_0),~\pa_zQ_0\rangle\\
&\qquad+\langle\mathcal{N}_2(\pa_z^2Q_0,\nabla d_0,\nabla d_1)+\mathcal{N}_2(\pa_z^2Q_0,\nabla d_1,\nabla d_0),~\pa_zQ_{k+1}\rangle\\
&=\int_{\mathbb{R}}\partial_z\Big(\big(\mathcal{N}_2(\pa_zQ_{k+1},\nabla d_0,\nabla d_1)+\mathcal{N}_2(\pa_zQ_{k+1},\nabla d_1,\nabla d_0)\big):\pa_zQ_0\Big)\ud z=0,
\end{align*}
and for $q+r=k+2$
\begin{align*}
& \langle \mathcal{N}_2(\pa_z^2Q_0,\nabla d_q,\nabla d_r)+ \mathcal{N}_2(\pa_z^2Q_0,\nabla d_r,\nabla d_q),~\pa_zQ_0\rangle\\
&=\frac12\int_{\mathbb{R}}\partial_z\Big(\big(\mathcal{N}_2(\pa_zQ_0,\nabla d_q,\nabla d_r)+\mathcal{N}_2(\pa_zQ_0,\nabla d_r,\nabla d_q)\big):\pa_zQ_0\Big)\ud z=0.
\end{align*}
Therefore, using the definition \eqref{Fm} and \eqref{f1}, it holds on $\Gamma$ that
\begin{align*}
 & \langle F_{k+2}, ~~\pa_zQ_0 \rangle +\langle F_1,~\pa_zQ_{k+1}\rangle\\
=&
\langle \pa_zQ_0 , ~\pa_zQ_0 \rangle(\pa_t-\Delta )d_{k+1}
-2\langle(\nabla d_{k+1}\cdot \nabla)\pa_zQ_0, ~\pa_zQ_0\rangle\\
&+2\langle \pa_zQ_{k+1}, ~\pa_zQ_0 \rangle (\pa_t-\Delta )d_0
-2(\nabla d_0\cdot \nabla)\langle \pa_zQ_{k+1}, ~\pa_zQ_0\rangle\\
&-\frac{L}{2}\Big(\langle\mathcal{N}_1(\pa_zQ_{k+1},\nabla d_0),~\pa_zQ_0\rangle
+\langle\mathcal{N}_1(\pa_zQ_0,\nabla d_{k+1}),~\pa_zQ_0\rangle+\langle \mathcal{N}_1(\pa_zQ_0,\nabla d_0),~\pa_zQ_{k+1}\rangle\Big)\\
&+\text{terms}\;\text{independent}\;
\text{of}\;d_i(i\geq k+1)-\langle B_{k+1}+C_{k+1},~\pa_zQ_0\rangle\\
&+d_{k+2}\langle G_0,~ \pa_zQ_0\rangle
+\langle (d_{1}-z)G_k,~\pa_zQ_0\rangle+\langle G_0(d_{1}-z),~\pa_zQ_{k+1}\rangle\\
\triangleq & J_1+J_2+\cdots +J_5. \mynumber{sum:Fk+2F1}
\end{align*}
Then similar with the calculation in the proof of Lemma \ref{lem:d1}, we arrive at equation \eqref{dk+1q11ex}.
\qed
\section*{Acknowledgments}
This work is supported by the National Key Research and Development Program of China (No. 2023YFA1008801). The first and the third author are also supported by NSF of China under Grant No. 12271476 and 11931010.

\end{document}